\documentclass[12pt]{article}
\usepackage{amsmath,amsfonts,amscd, amsthm,latexsym}

\newcommand{\ep}{\varepsilon}

\newcommand{\centr}[2]{\mathop{{\rm center}_{#1}}(#2)}
\newcommand{\cont}[1]{\mathop{\rm cont}_{#1}}

\newcommand{\LCS}[1]{\mathop{\rm LCS}#1}

\newcommand{\mult}[2]{\mathop{\rm mult}_{#1}#2}
\newcommand{\NE}[1]{\mathop{\rm NE}(#1)}
\newcommand{\NEc}[1]{\mathop{\overline{\rm NE}}(#1)}
\newcommand{\Pic}[1]{\mathop{\rm Pic}(#1)}
\newcommand{\pt}{\mathop{\rm pt.}}
\newcommand{\R}[1]{\mathop{\rm R}#1}
\newcommand{\reg}[1]{\mathop{\rm reg}#1}
\newcommand{\Supp}[1]{\mathop{\rm Supp}#1}

\newcommand{\Atd}[2]{\widetilde{\rm A}_{#1}^{#2}}
\newcommand{\Dtd}{\widetilde D}
\newcommand{\Xtd}{\widetilde X}

\newcounter{subsec}[section]
\theoremstyle{definition}
\newtheorem*{definition}{\thesection.\thesubsec. Definition}
\newtheorem*{conjecture}{\thesection.\thesubsec. Conjecture}
\newtheorem*{p-definition}{\thesection.\thesubsec.
 Proposition-Definition}
\theoremstyle{plain}
\newtheorem*{corollary}{\thesection.\thesubsec. Corollary}
\newtheorem*{lemma}{\thesection.\thesubsec. Lemma}
\newtheorem*{proposition}{\thesection.\thesubsec. Proposition}
\newtheorem*{statement1}{\thesection.\thesubsec.1}
\newtheorem*{statement2}{\thesection.\thesubsec.2}
\newtheorem*{statement3}{\thesection.\thesubsec.3}
\newtheorem*{theorem}{\thesection.\thesubsec. Theorem}
\theoremstyle{remark}
\newtheorem*{remark}{\thesection.\thesubsec. Remark}
\newtheorem*{remcor}{\thesection.\thesubsec.
Remark-Corollary}
\newtheorem*{example}{\thesection.\thesubsec. Example}
\newtheorem*{example-p}
{\thesection.\thesubsec. Example-Problem}
\newtheorem*{example-c}
{\thesection.\thesubsec. Example-Corollary}
\newtheorem*{question}{\thesection.\thesubsec. Question}

\input cyracc.def
\font\tencyr=wncyr10
\def\cyr{\tencyr\cyracc}

\date{November 16, 1997}

\title{COMPLEMENTS ON SURFACES}
\author{V.V. Shokurov}

\begin{document}

\maketitle

\begin{abstract}
The main result of the paper is
a boundedness theorem for $n$-complements
on algebraic surfaces.
In addition, applications of this theorem
to a classification of log Del Pezzo surfaces
and of birational contractions for $3$-folds
are formulated
\footnote[1]{1991 Mathematics Classification:
14E05, 14E30, 14J25, 14J30.}
\footnote[2]{Partial financial support was provided by the NSF under
the grant number DMS-9500971.}.
\end{abstract}

\tableofcontents

\section{Introduction}
\label{int}

An introduction to complements can be found
in \cite[\S 5]{Sh2}.
See also the end of this section.

\refstepcounter{subsec}
\label{contex}
\begin{example}~Let $S=\mathbb P_E(F_2)$ be
a ruled  surface over
a non-singular curve $E$ of genus $1$,
which corresponds to a non-splitting vector bundle
$F_2/E$ of rank 2 \cite[Theorem 5 (i)]{At}.
We denote the ruling by $f\colon S\to E$.
It has a single section $\mathbb P_E(F_1)$
in its linear and even numerical equivalence class
(by \cite[2.9.1]{Sh1}, cf. arguments below).
We identify the section with $E$.

Note that $E|_E\sim \det F_2\sim 0$,
where $\sim$ denotes {\it linear\/} equivalence
(cf. \cite[Proposition 2.9]{H}).
On the other hand $(K_S+E)|_E\sim K_E\sim 0$
according to the Adjunction Formula.
Thus $(K_S+2E)|_E\sim 0$ \cite[Lemma 2.10]{H}.
But $K_S+2E\equiv 0/E$,
where $\equiv$ denotes {\it numerical\/} equivalence.
So, $K_S+2E=f^*((K_S+2E)|_E)\sim 0$ \cite[2.9.1]{Sh1}.
Equivalently, we may choose $-K_S=2E$ as
an anti-canonical divisor.
In particular, $-K_S$ is nef because $E^2=0$.

The latter implies also that
the cone $\NEc{S}$ is two-dimensional and
has two extremal rays:
\begin{itemize}
\item the first one $R_1$ is generated by
a fibre $F$ of the ruling $f$, 
\item the second one $R_2$ is generated by
the section $E$.
\end{itemize}
We contend that $E$ is the only curve
in $R_2$.

Indeed if there is a curve $C\not= E$ in $R_2$,
then $C\equiv m E$ where $m=(C.F)$
is the degree of $C$ as a multi-section of $f$.
We may suppose that $m$ is minimal.
Then as above, $C\sim m E$ which
induces another fibering with projection
$g\colon S\to \mathbb P^1$.
All fibres of $g$ are non-singular curves
of genus $1$.
However, some of them may be multiple.
According to our assumptions,
any multiple fibre $G_i$ has
a multiplicity $m_i=m\ge 2$, and
we have at least two of them.
So, if $G_1=E$ then another
multiple fibre $G_2\equiv E$ gives
another section in $R_2$.
This is impossible.

Thus in general a nef anti-canonical divisor
$-K$ may not be semi-ample,
and we may have no complements when
$-K$ is just nef.
The same may occur in the log case.
\end{example}

\refstepcounter{subsec}
\label{posch}
\begin{example}
Now let $E$ be a fibre of
an elliptic fibering $f:S\to Z$ such that
\begin{itemize}
\item $E$ is a non-singular elliptic curve,
\item $m$ is the multiplicity of $f$ in $E$,
i.~e., $m$ is the minimal positive integer with $m E\sim 0$
locally$/Z$.
\end{itemize}
Then, for a canonical divisor $K$,
$K\sim a E$ locally$/Z$ with
a unique $0\le a\le m-1$.  
Moreover in the characteristic $0$,
$a=m-1$, or $K+E\sim 0$, which means that $K+bE$, 
with some real $0\le b\le 1$, has a $1$-complement locally$/Z$.  
However if
the characteristic $p$ is positive, one can construct such a fibering
with $m=p^k$ and any $0\le a\le m-1$.  
For instance, see
\cite[Example 8.2]{KU} for $m=p^2$ and $0\le a\le p-1$.  So,
$K+bE$, for $b$ close to $1$, has no $n$-complements with bounded
$n$ because such a complement would have the form $K+E\sim
(a+1)E$ and $n(K+E)\sim n(a+1)E\not\sim 0$ when $n,a\ll p$ and
$m\ge p$.
So, $K+E$ can be an $n$-complement having rather
high index $n|m$.
However, it is high in the Archimedean sense, 
but $n$ is small
$p$-adically when $m=p^k$.  \end{example}

So, everywhere below we suppose that
the characteristic is $0$.

\refstepcounter{subsec}
\theoremstyle{definition}
\newtheorem*{conjcs}{\thesection.\thesubsec.
Conjecture on complements}
\begin{conjcs}
\label{conjcs}
We consider log pairs $(X,B)$
with boundaries $B$
such that $K+B$ is log canonical.
Then complements in any given
dimension $d$ are {\it bounded}.
This means that there exists a finite set $N_d$ of
natural numbers such that
any contraction $f:X\to Z$, satisfying
certain conditions which we discuss below, and of $\dim X\le d$,
has locally$/Z$ an $n$-complement
for some {\it index\/} $n\in N_d$.

Of course, we need to assume that
\begin{description}
\item{(EC)} there exists at least some $n$-complement.
\end{description}
In particular, if $B=0$ then the linear system
$|-nK|\not=\emptyset$
(cf. Corollaries~\ref{int}.\ref{nonvan1}-\ref{nonvan2} 
and (NV) in Remark~\ref{indcomp}.\ref{rel})
and its generic element
has good singularities in terms of complements
(see explanations before 
Corollary~\ref{int}.\ref{nonvan1}).
Maybe it is enough 
the existence of complements (EC).
More realistic and important for applications
are the following conditions
\begin{description}
\item{(SM)} the multiplicities $b_i$ of $B$
are {\it standard}, i.~e.,
$b_i=(m-1)/m$ for a natural number $m$, or
$b_i=1$ (as for $m=\infty$); and
\item{(WLF)} $(X/Z,B)$ is a {\it weak log Fano\/} which means
that $-(K+B)$ is nef and big$/Z$.
\end{description}
Note that these conditions imply (EC)
according to the proof of \cite[Proposition 5.5]{Sh2}.
By Example~\ref{int}.\ref{contex}
a condition that $-(K+B)$ is just nef,
is not sufficient even for (EC).
\end{conjcs}

\refstepcounter{subsec}
\theoremstyle{plain}
\newtheorem*{mainth}{\thesection.\thesubsec.
Main Theorem}
\begin{mainth}
The complements in dimension $2$
are bounded under the condition (WLF)
and
\begin{description}
\item{{\rm (M)}} the multiplicities $b_i$ of $B$
are {\rm standard\/}, i.~e., $b_i=(m-1)/m$ 
for a natural number $m$, or
$b_i\ge 6/7$.
\end{description}
More precisely, for almost all
contractions and all contractions
of relative dimension $0$ and $1$,
we can take a complementary index in $RN_2=\{1,2,3,4,6\}$.
\end{mainth}

Here {``}almost all" means up to
a bounded family of contractions
in terms of moduli.
Really, this concerns moduli spaces
in the global case or that of log Del Pezzos.
Moduli itself may not
be the usual one according to
Remark~\ref{int}.\ref{mod}.
By the {\it global\/} case we mean that $Z=\pt$.
The other cases are {\it local\/}$/Z$.

\refstepcounter{subsec}
\begin{definition}
\label{dexc}
(Cf. \cite[Theorem 5.6]{Sh2}.)
We say that a complement
$K+B^+$ is {\it non-exceptional\/} if
it is not Kawamata log terminal
whenever $Z=\pt$, and it is not
purely log terminal
on a log terminal resolution whenever $Z\not=\pt$
(cf. \cite[Theorem~5.6]{Sh2}).
The former will be called {\it global\/}
and the latter {\it local\/}.
Respectively, the corresponding
log pair $(X/Z,B)$ is called {\it non-exceptional\/}
if it has a non-exceptional complement.
For instance, surface Du Val singularities
of types $\mathbb A_*$ and $\mathbb D_*$ are
non-exceptional even though they have trivial
complements that are canonical
(cf. Example~\ref{int}.\ref{tcomp}).
They have respectively other
$1$- and $2$-complements that
are non-exceptional.

On the other hand, the pair $(X/Z,B)$ is {\it exceptional\/}
if each of its complements $K+B^+$ is {\it exceptional\/}.
The latter means that $K+B^+$ is Kawamata log terminal
whenever $Z=\pt$, and it is
{\it exceptionally\/} log terminal whenever $Z\not=\pt$.
The exceptional log terminal property
means the purely log terminal property on
a log terminal resolution.
The Du Val singularities of exceptional types
$\mathbb E_6,\mathbb E_7$ and $\mathbb E_8$ are
really exceptional from this point of view
as well (cf.~\cite[Example 5.2.3]{Sh2}).

If the pair $(X/Z,B)$ is exceptional
and $Z\not=\pt$,
then for any complement $K+B^+$, 
any possible
divisor (at most one as one can prove)
with log discrepancy $0$ for $K+B^+$ has
the center over the given point in $Z$.
Otherwise we can find another complement
$B^{+'}>B^+$ that is non-exceptional
(see the Proof of  
Theorem~\ref{lcomp}.\ref{lmainth}: General case
in Section~\ref{lcomp}).

For instance, in the Main Theorem
the $n$-complements with $n$ not in $RN_2$
are over $Z=\pt$ and exceptional 
as we see later in Theorems~\ref{lcomp}.\ref{lmainth}
and \ref{gcomp}.\ref{gmainth}.
So, the theorem states that the global exceptions
are bounded.
Some of them will be discussed in Section~\ref{ecomp}.
In higher dimension, it is conjectured
that such complements and the corresponding pairs
$(X/\pt,B)$ under conditions (WLF) and (SM) are bounded.
Of course, this is not true if we drop (WLF)
(see Example~\ref{int}.\ref{tcomp}).
However, it may hold formally  and
even in the local case under certain conditions
as suggested in
Remark~\ref{int}.\ref{exc}.

$N_2$ in the Conjecture on complements
differs from $RN_2$ in the Main Theorem in
exceptional cases,
which are still not classified completely.
Nonetheless we define corresponding {\it exceptional\/} indexes
as $EN_2=N_2\setminus RN_2$.
We say that $RN_2$ is {\it regular\/} part
of $N_2$ in the Conjecture for dimension $2$.
Note that $RN_2=N_1=RN_1\cup EN_1$ where
$RN_1=\{1,2\}$ and $ER_1=\{3,4,6\}$.
\end{definition}

\refstepcounter{subsec}
\begin{conjecture}
\label{econj}
In general we may define
$EN_d=N_2\setminus N_{d-1}$ and
conjecture that $RN_d=N_{d-1}$,
or, equivalently, that
$EN_d$ really corresponds to exceptions.
Evidence supporting this is related
to our method in the proof
of the Main Theorem and Borisov-Alekseev's
conjecture~\cite{Al}.

Some extra conditions on multiplicities
are needed: (SM) or an appropriate
version of (M) (cf. (M)' in
Remark~\ref{gcomp}.\ref{mrem}.2,
(M)'' in \ref{csing}.\ref{singth}.1,
and the statement~\ref{indcomp}.\ref{indth}.1).
\end{conjecture}

Basic examples of complements may be found
in~\cite[5.2]{Sh2}.

\refstepcounter{subsec}
\begin{example-p}
\label{tcomp}
{\it Trivial complements.}
Let $X/Z$ be a contraction
having on $X$ only log canonical
singularities and with $K\equiv 0/Z$,
e.~g., Abelian variety, K3 surface,
Calabi-Yau 3-fold or a fibering of them.
Then it is known that locally$/Z$,
$K$ is semi-ample:
$K\sim_{\mathbb Q} 0$ or
$n K\sim 0$ for some natural number $n$
(at least in the log terminal case,
cf.~\cite[Remark~6-1-15(2)]{KMM}).
Such minimal $n$ is known as the {\it index\/}
of $X/Z$, to be more precise,
over a point $P\in Z$.
Indices are {\it global\/} when $Z=P=\pt$ and
{\it local\/} otherwise.
So, if we consider thought $X/Z$ as 
a log pair with $B=0$,
then (EC) is fulfilled for the above $n$.
In the global case
the Conjecture on complements suggests that
we can find such $n$ in $N_d$
for a given dimension $d=\dim X$.
Moreover, in that case,
we may replace $B=0$ by $B$ under (SM)
according to 
the Monotonicity~\ref{indcomp}.\ref{monl1} below,
and assuming that $X$ may be semi-normal.
In the local case we need an additional
assumption on a presence of a log canonical
singularity$/P$, i.~e.,
there exists an exceptional or non-exceptional divisor $E$
with the log discrepancy $0$ for $K+B$ and
with $\centr{X}{E}/P$.
Otherwise a complement may be {\it non-trivial\/}:
$B^+>B$.

This really holds in dimensions $1$ and $2$
by Corollary~\ref{int}.\ref{ind2} below.
In dimension $3$ it is know that
any global index $n$ divides
the Beauville number
$$b_3=2^5\cdot 3^3\cdot 5^2\cdot 7\cdot 11
\cdot 13\cdot 17\cdot 19$$
when $X/\pt$ has at most terminal singularities
and $B=0$.
Note that the similar $2$-dimensional number
is $b_2=12$.
Its proper divisors give $N_1=RN_2$.
So, we conjecture that $N_2$ consists of
(the proper) divisors of the Beauville number.
Perhaps we have something like this  
in higher dimensions.

Of course, if in a given dimension $d$ indices
are bounded, we have a universal index $I_d$
as their least common multiple.
The case of the Conjecture on complements
under consideration in dimension $d$
is equivalent to an existence of such $I_d$.
As always in mathematics, it is known only
that $I_1=12$
and that $I_2$ exists (see Corollary~\ref{int}.\ref{ind2}).
We suggest that $I_2\approx b_3$.
In particular, this means that $b_3$
corresponds to non-exceptional cases.
Note also that in the global case,
according to our definition,
contractions and their complements are
non-exceptional whenever they are Kawamata
log terminal, for instance, Abelian varieties.
These varieties have unbounded moduli
in any dimension $\ge 2$,
but their index is $1$.
So, from the formal viewpoint of
Remark~\ref{int}.\ref{exc} below,
we treat them as a regular case.

Note that we are still discussing 
the case that
is opposite to the assumption (WLF)
in the Main Theorem.
Indeed, if we consider a birational contraction
$X/Z$ then the indices are not bounded.
Take, for instance, $X=Z$ as
a neighborhood of a quotient singularity.
Nonetheless, the above should hold when
$X$ is really log canonical in $P$
as we suppose.
In the surface case such singularities
with $B=0$ are known as {\it elliptic\/}
and they have in a certain sense
a Kodaira dimension $0$.
Their universal index is $12$.
The Conjecture on complements implies
an existence of such a universal number in any
dimension.

An inductive explanation may be presented
in the log canonical case (cf. Remark~\ref{int}.\ref{exc}).
Using the LMMP it is possible to reduce
to the case when $B$ has a reduced component
$E$.
Then by the Adjunction formula the index
of $K+B$ near $E$ coincides
with the same for restriction $K_E+B_E=(K+B)|_E$.
So, it is not surprising
that indices for dimension $d-1=\dim E$
should divide indices in dimension $d$.
If the above restriction is epimorphic
for Cartier multiples of $K_E+B_E$
then both indices coincide over
a neighborhood of $P\in Z$.
This holds, for example, when $X/Z$ is birational.
We see in Section~\ref{lcomp} that this works
in fibre cases as well.

Finally, note that in the fibre case
the complements may be non-trivial
when  the above condition on a log canonicity
is not fulfilled.
Let $X$ be smooth and $B=0$.
Then the index of $K$ in a neighborhood
of a fibre $F\subset X$ over $P$ may
be arbitrarily high.
The point is that contraction $X/Z$ itself
is not smooth in those cases and
fibre $F$ is multiple.
For instance, $K$ and $F$ have
multiplicity $m$ near fibre $F$ of
type $mI_0$ in an elliptic fibering.
However, there exists $1$-complement
$K+F$ in a neighborhood of $F$.
This explains in particular why
complements are fruitful even
in a well-known case as
elliptic fiberings.
More details may be found in Section~\ref{lcomp}.
In the positive characteristic $p$, it does not work well
according to Example~\ref{int}.\ref{posch}.
Maybe it works modulo $p$ factors.

Each {\it trivial\/} complement, i.~e.,
a complement with $B^+=B$, defines
a cover $\Xtd\to X$ of degree $d$,
which is unramified in divisors with respect to $B$:
a covering is $B$-{\it unramified\/} if it
preserves the boundary $B$, i.~e.,
inequality $\mult{D_i}{B}=b_i\ge (m_j-1)/m_j$,
where $m_j$ is any ramification multiplicity$/D_i$ 
\cite[2.1.1]{Sh2}, holds for each prime $D_i$ in $X$.
The conjecture in this case states
a boundedness of such degrees.
In terms of a local or global algebraic fundamental group
$\pi (X/Z)$, they correspond to
normal subgroups of finite index and
with a cyclic quotient.
A general conjecture here states that
$\pi (X/Z)$ is quasi-Abelian, and
even finite in exceptional cases
(see Remark~\ref{int}.\ref{exc}).
The structure of the fundamental group is
interesting even for non-singular $X$.
\end{example-p}

In the previous example we know that
the condition (EC) is satisfied.
So, according to the general philosophy
and to conjectures there, 
the following results are not surprising.

\refstepcounter{subsec}
\begin{corollary}
\label{mcorol}
In the Main Theorem we may replace
the condition (WLF) by (EC).
\end{corollary}

This will be proven in
Sections~\ref{lcomp} and \ref{gcomp}.
More advanced results can be found
in Remark~\ref{gcomp}.\ref{mrem}.
Here we consider only a {\it trivial\/} case.

\begin{proof}[Proof: Numerically trivial case
and under (SM)]
So, we suppose that $K+B\equiv 0/Z$ and
it has a trivial complement $B^+=B$
(cf. Monotonicity~\ref{indcomp}.\ref{monl1}).
Moreover we consider the global case
when $Z=\pt$ (for the local case see
Section~\ref{lcomp}).
We need to bound the index of $K+B$.

If $B=0$ and $X$ has only canonical singularities
it is well-known that the index of $K$ divides $I_1=12$.
On the other hand,
pairs $(X,B)$ with $B\not=0$ and surfaces $X$
with non-canonical singularities
are bounded whenever they
are $\varepsilon$-log canonical for
any fixed $\varepsilon >0$, for instance, for
$\varepsilon\ge 1/7$ \cite[Theorem~7.7]{Al}.
So, their indices are bounded as well.
Note, that by (SM) the possible multiplicities
$b_i$ satisfy the d.c.c. (descending chain condition).

Thereby we need to bound the indices
in the non-$\varepsilon$-log canonical case
for an appropriate $\varepsilon$.
We take $\varepsilon=1/7$.
After a log terminal resolution we may suppose
that $B$ has a multiplicity $b=b_i>6/7$
in a prime divisor $D=D_i$.
The resolution has the same indices
(cf. \cite[Lemma~5.4]{Sh2}).
Moreover,
for the regular complements or
for almost all $(S,B)$, all $b_i=1$,
when $b_i\ge 6/7$, and
the index is in $RN_2$.
So, there exists a real $c>0$ such that
all $b_i=1$ whenever $B_i\ge c$,
and we have then a regular complement.
To reduce this to the Main Theorem
we apply the LMMP.
(Cf. Tsunoda's \cite[Theorem~2.1]{T},
which states that $K+B$ has the index $\le 66$,
when $B$ is reduced and $K+B\equiv 0$.)

We have an extremal contraction $f:X\to Z$
because $K+B-b D$ is negative on a covering
family of curves.
If $f$ is birational we replace $X$ by $Z$.
Then $D$ is not contracted since
$(K+B-b D.D)=-b D^2>0$.

If $Z$ is a curve, then $f$ is a ruling,
and, according to a 1-dimensional result
as our corollary, applied to the generic fibre ,
$b_i=1$ and $K+B$ has a regular index
near the generic fibre of $f$.
Note that $D$ is a multi-section of $f$ of
a multiplicity at most $2$.
If $D$ is a $1$-section we reduce the problem
to a $1$-dimensional case after an adjunction
$K_D+B_D=(K+B)|_D$.
In that case we use the same arguments
as in Example~\ref{int}.\ref{contex}
(cf. Lemma~\ref{indcomp}.\ref{indcomp1}).
We can take a complement with an even index
when $K+B$ has index $2$ near the generic fibre. 

Otherwise $D$ is an irreducible double section of $f$.
In that case the index is in $RN_2$ near each fibre. 
The same holds for entire $X$,
except for a case when $B=D$, 
$K_D+B_D=(K+B)_D\sim 0$ and $2(K+B)\sim 0$
(see Lemma~\ref{indcomp}.\ref{indcomp2}).

Finally, if $Z=\pt$, we may apply the Main Theorem
for $B-\delta D$ for some $\delta>0$.
The main difficulties are here.
\end{proof}

\refstepcounter{subsec}
\begin{corollary}
\label{ind2}
Under (SM),
$I_2$ exists.
Moreover, for the trivial complements,
each $b_i=1$ whenever $b_i\ge I_2/(I_2+1)$,
and , for the almost all trivial complements,
each $b_i=1$ whenever $b_i\ge 6/7$.

In the global case when $(X,B)$
is not Kawamata log terminal, and
in the local case we can replace
inexplicit $I_2$ by $I_1=12$,
(SM) by (M).
Then each $b_i=1$ whenever $b_i\ge 6/7$.

If we have an infinite number of
exceptional divisors with log discrepancies $0$
over $P$, then we can replace $I_1=12$ even
by $2$.
\end{corollary}

\begin{proof}[Proof: Global case]
We suppose that $Z=\pt$.
We consider other cases in
Section~\ref{lcomp}.
Then the results follow from the above
arguments and by \cite[Theorem~7.7]{Sh2}.
The non-Kawamata log terminal case
corresponds to regular complements by
the Main Theorem (cf.
the Inductive Theorem and
Theorem~\ref{gcomp}.\ref{gmainth} below).

For global and local cases, with
an infinite number of
exceptional divisors with log discrepancies $0$,
see \ref{indcomp}.\ref{indth}.2.
\end{proof}

\refstepcounter{subsec}
\label{ind3}
\begin{corollary}
Under (SM),
let $P\in (X,B)$ be a 3-fold
log canonical singularity as
in Example~\ref{int}.\ref{tcomp}
having the trivial complement.
Then its index divides $I_2$.
Moreover, each $b_i=1$ whenever $b_i\ge I_2/(I_2+1)$.

If $(X,B)$ has at least two exceptional
divisors$/P$ with the log discrepancy $0$, then
we can replace
inexplicit $I_2$ by $I_1=12$.
The same holds with one exceptional or
non-exceptional divisor with the center
passing through $P$ but $\not=P$.

If $(X,B)$ on some resolution has a triple point in an
exceptional locus of divisors$/P$ with the log discrepancy $0$,
then we can replace $I_1=12$ even by $2$.
\end{corollary}

\begin{proof}
According to the very definition of
the trivial complements, we have an exceptional divisor $E/P$
with the log discrepancy $0$.
By the LMMP (cf. \cite[Theorem~3.1]{Sh3}),
we can make a crepant resolution $g:Y\to X$ just of
this $E$ and $E$ is in the reduced part of boundary $B^Y$.
By the adjunction $K_E+B_E=(K_Y+B^Y)|_E\equiv 0$ and
$B_E$ has only standard coefficients.
So, it has a trivial complement.
Then by Corollary~\ref{int}.\ref{ind2}
and by the local inverse adjunction
$K_Y+B^Y$ has the same Cartier index and it divides $I_2$
whenever the restriction $K_E+B_E=(K_Y+B^Y)|_E$ is
Kawamata log terminal.
Here we may use a covering trick as well.
Then, by the Kawamata-Viehweg vanishing,
restrictions $m(K_E+B_E)=m(K_Y+B^Y)|_E$ are epimorphic
for $m$ dividing the index of $K_E+B_E$.
So, $K+B$ has the same index and it divides $I_2$.

Other special cases follow from
special cases in Corollary~\ref{int}.\ref{ind2}
(cf. \cite[Corollary~5.10]{Sh2}).
Then the index will divide $I_1=12$.
Here a difficulty is related to an application
of the above surjectivity to a semi-normal surface $E$.
However, its combinatoric is quite simple.
\end{proof}

We may simulate examples on well-known
varieties as in the next example.

\refstepcounter{subsec}
\label{simul}
\begin{example} Let $X=\mathbb P^d$
and $B=\sum b_i L_i$ where
prime divisors $L_i$ are hyperplanes
in a generic position and 
there are at most $l$ of them.
(The latter holds when $\min\{b_i\not=0\}>0$.)
Then for any boundary coefficients
$0\le b_i\le 1$, the log pair $(X,B)$
is log terminal.
If we have a reduced component $L_i$ in $B$
we can restrict our problem on $L_i$.
So, really new complements are related
to the case when all $b_i<1$ that
we assume below.
Note that if all $b_i$
are rational (or $\sum b_i<d+1$), the condition (EC)
follows from an inequality
$\sum b_i\le d+1$.
An existence of an $n$-complement
is equivalent to another inequality
$$\deg\lfloor (n+1)B\rfloor=
\sum \lfloor (n+1)b_i\rfloor/n\le d+1.$$
So, the Conjecture on complements
states that we may choose such $n$
in a finite set $N_d$.
The space
$$
T_d^l=\{(b_1,\dots,b_l)\mid b_i\in [0,1]
\text{ and } \sum b_i\le d+1\}
$$
is compact and 
the vectors $v=(b_1,\dots,b_l)$ without
$n$-complements form a union of
convex rational polyhedra:
$b_i\ge c_i$.
In the intersection they give 
a vector $v$ without any $n$-complements.
We can also assume that it is maximal
$\sum b_i=d+1$, and $b_i>0$.
But this is impossible,
because there exists an infinite set of 
approximations with natural numbers $n$ that
$(n+1)b_i=N_i+\delta_i, i=1,2,$
where $|\delta_i-b_i|<\varepsilon\ll 1$,
and $N_i$ is integer.
Indeed, then  each $\lfloor (n+1)b_i\rfloor=N_i$, 
and
$$
\sum \lfloor (n+1)b_i\rfloor/n=\sum N_i/n=
\sum b_i+\sum (b_i-\delta_i)/n\le
d+1+l\varepsilon /n.
$$
Hence $\sum \lfloor (n+1)b_i\rfloor/n\le d+1$,
when $l\varepsilon<1$, because $d$ is integer.
We can find the required approximation from
another one $nb_i=N_i+\varepsilon_i$,
where $\varepsilon_i=\delta_i-b_i$,
$|\varepsilon_i|<\varepsilon\ll 1$,
and $N_i$ is integer.
It has required solutions $n$,
according to the Kronecker Theorem
\cite[Theorem~IV in Section 5 of Chapter III]{Ca},
with $L_i=b_i x$ and $\alpha_i=0$.
Indeed, for every integer numbers $u_i$,
$\sum u_i\alpha_i=0$ is integer.

Explicitly, it is known up to $d=2$.
For $d=1$ \cite[Example 5.2.1]{Sh2}.
For $d=2$, under (SM) $n\le 66$ due to 
Prokhorov~\cite[Example~6.1]{P2}
(cf. \cite[ibid]{T}).

We may generalize this
and suppose that $d$ is any positive rational
(and, maybe, even real), 
the hyperplanes are in a non-generic position, and with
hypersurfaces rather than hyperplanes of fixed
degrees.
Moreover, we can consider complements with
$n$ divided by a given natural number $m$
(cf. Lemma~\ref{indcomp}.\ref{compl3}).
(Even for $l=\infty$, as for $d=1$.
However, the boundedness is
unknown for
$l=\infty$ and $d>1$.)

There exists a corresponding local case
when we consider a contraction
$X/Z$ with a fibre $F=\mathbb P^d$
and $B=\sum b_i L_i$ where
prime divisors $L_i$ intersect $F$ in
hyperplanes in a generic position.

Of course, more interesting cases
are related to a non-generic position and with
hypersurfaces rather than hyperplanes.
They also give non-trivial examples 
of trivial complements.
For instance, if $C\subset \mathbb P^2$
is a plane curve of degree $6$ with
one simple triple point singularity,
then $K+C/2$ is log terminal and
has a trivial $2$-complement.
The corresponding double cover
$\Xtd\to \mathbb P^2$ produces
a K3 surface $\Xtd$ with a single
Du Val singularity of type $\mathbb D_4$
over the triple singularity of $C$.
For a generic $C$, after a resolution of $\Xtd$
we get a non-singular K3 surface $Y$ with
the Picard lattice of rank $5$
and an involution on $Y$ which is
identical on the lattice.
Note that generic curve $C$ is
a generic trigonal curve of genus $7$.
But this is a different tale.
\end{example}

\refstepcounter{subsec}
\begin{example-p}
\label{gquot}
Other interesting complements
correspond to {\it Galois quotients\/}.
Let $G:X$ an (effective) action of
a finite group $G$ on a log pair $(X/Z,B)$,
with a boundary $B$ under (SM),
and $-(K+B)$ is nef.
For instance, $X=\mathbb P^n$ or a Fano variety,
Abelian variety, Calabi-Yau 3-fold,
or an identical contraction of
a non-singular point.
Then on the quotient $f:X\to Y=X/G$ we
have a unique boundary $B_Y$ such
that $K+B=f^*(K_Y+B_Y)$.
Hence, $(Y/Z,B_Y)$ will be a log pair
of the same type as $(X/Z,B)$.
The (minimal) complement index $n$ in this case
is an invariant of the action of the group.
In the global $1$-dimensional case $X=\mathbb P^1$
there exists an action of a finite group
$G\subset PGL(1)$ with $n\in N_1$.
All exceptional cases have the quotient description.
But in dimension $d\ge 1$ it may be that
some complement indices do not correspond to
quotients of $PGL(d)$ and not all
exceptional cases have the quotient description
for $X=\mathbb P^d$ or some other non-singular
Del Pezzos.
Even in the $1$ dimensional case
$Y=\mathbb P^1$ with $B=(n-1)P/n +(m-1)Q/m$
and $n>m\ge 2$ does not correspond to
a Galois quotient.
In higher dimensions we expect
more asymmetry as in a modern cosmology.
But among
symmetric minority we may meet
real treasures.

Quotients $X/G$, for finite groups of automorphisms of
K3 surfaces $X$, reflect geometry of
pairs $(X,G)$.
According to Nikulin, such pairs $(X,G)$ are bounded
when $G$ is non-trivial (cf. 
\cite[Theorem~18.5]{BPV}),
in particular, with a non-symplectic action.
On the other hand each Abelian surface $X$ has
an involution.
So, pairs $(X,\mathbb Z_2)$ are not bounded in that case.
However, we anticipate that the pairs $(X,G)$ 
with $K\equiv 0$ and non-singular $X$ are bounded
when $(X,G)$ is quite non-trivial (non-regular).
For instance, 
$X$ has irregularity $0$ and
$G\not\cong \mathbb Z_m$ with
proper divisors $m|12$. 
We may say that such pairs $(X,G)$ are 
{\it exceptional\/} for dimension $2$.
Indeed, if we have such a group $G$ then
it has a fixed point, according to 
the classification of algebraic surfaces. 
If the fixed point is not isolated it produces
a non-trivial boundary on $X/G$.
Otherwise it gives a non Du Val singularity on $X/G$
when the action is non-symplectic.
Then $X/G$ with a boundary belongs to a bounded family
by Corollary~\ref{int}.\ref{mcorol}
and Alekseev~\cite[ibid]{Al}.
Perhaps this holds sometimes for symplectic actions
as well as for $K3$ surfaces.
\end{example-p}

\refstepcounter{subsec}
\begin{remark}
\label{mod}
Moduli spaces mentioned after the Main Theorem
may have real parameters corresponding to
boundary coefficients.
If we would like to have usual algebraic moduli
we can forget boundaries or impose a condition
as (SM).
A little bit more generally, we may suppose that
$1$ is the only accumulation point for
possible boundary multiplicities.
\end{remark}

\refstepcounter{subsec}
\begin{remark}
\label{exc}
In addition to Conjecture~\ref{int}.\ref{econj}
we suggest that regular $RN_d=N_{d-1}$ is enough for
global non-exceptional and any local complements
in dimension $d$.
We verify it in this paper for $RN_2=N_1$.
An explanation is related to~\cite[Lemma 5.3]{Sh2}.
If we really have log canonical singularities
we can induce the problem from a lower dimension
(cf.~\cite[Proof of Theorem 5.6]{Sh2} and
the Induction Theorem~\ref{indcomp}.\ref{indth} below).
It means that in this case we need only
indices from $N_{d-1}$ for dimension $d$.
The same holds, if we can increase $B$ to $B'$ preserving
all requirements on $K+B$ for $K+B'$ and
$K+B'$ having a log canonical singularity.
A good choice is to simulate construction
of complements.
So, we take $B'=B^+=B+H/n$ where
$H\in |-nK-\lfloor (n+1)B\rfloor|$
for some $n$.
If $K+B^+$ is log canonical it gives
an $n$-complement.
This complement is bounded when $n\in N_d$ or
$N_{d-1}$.
Otherwise we have a log singular case in
which $H$ is called a {\it singular\/} element and
$K+B^+$ is a {\it singular\/} complement.
Keel and McKernan referred to such $H$ as
a tiger but it looks more like a must.
In this case we may consider
a weighted combination $B'=a B+(1-a)B^+$
and for an appropriate $0<a<1$, log divisor $K+B'$
will be log canonical but not Kawamata log terminal.
Therefore we have no reduction to dimension $d-1$
when we have no singular elements or such complements.
These cases correspond to exceptions.
In particular, $|-nK-\lfloor (n+1)B\rfloor|=\emptyset$ for
each $n\in RN_d$ in the exceptional case.
However, $|-nK-\lfloor (n+1)B\rfloor|\not=\emptyset$
for some $n\in RN_d$ in the regular cases.
(Cf. Corollaries~\ref{int}.\ref{nonvan1}-\ref{nonvan3}
below for $d=2$.)

In the local case singular elements are easy to construct
adding a pull back of a hyperplane section
of the base $Z$.
In the global case singular elements defines
an ideal sheaf of a type of Nadel's multipliers ideal sheaf.
If we have no singular elements and $B=0$
it implies an existence of a K\"ahler-Einstein metric
with a good convergence in singularities
(cf.~\cite{Na}).
It is then known that $X$ is stable in
the sense of Bogomolov (cf. \cite[Proposition~1.6]{Si}) 
and exceptions with $B=0$ should be bounded.
An algebraic counterpart of this idea is 
the Borisov-Alekseev's conjecture.
So, we may expect the same for exceptions with $B\not=0$
at least under (SM).
This will be proven for dimension $2$
in Section~\ref{gcomp}.
We discuss some exceptions in Section~\ref{ecomp}.

We think the same or something close holds
in any dimension 
(which will be discussed elsewhere), 
as well as for exceptional and non-exceptional cases
in the {\it formal sense\/}.
This means that they have an appropriate index:
some $n\in RN_d$ in the regular case, and
only some $n\in ER_d$ in the exceptional case.
In the local cases we replace $d$ by $d-1$.
For instance, any Abelian variety is
exceptional according to Definition~\ref{int}.\ref{dexc},
but they are non-exceptional from the formal point of view.
Their moduli are unbounded as we anticipate
in the non-exceptional case.
On the other hand, an elliptic fibre of type III or
any other exceptional type is formally exceptional.
In a certain sense their moduli are bounded.
However the log terminal singularities
of this type, i.~e., with the same graph
of a minimal resolution, will be unbounded,
but again bounded if we suppose
the $\varepsilon$-log terminal property.
For instance, the Du Val or surface canonical
singularities of the exceptional types are
bounded up to a certain degree 
(cf. Corollaries~\ref{csing}.\ref{exbound}-\ref{exboundf}).
This is a local version of the Borisov-Alekseev's conjecture.
\end{remark}

\refstepcounter{subsec}
\begin{remark}
Complements and their indices allow us
to classify contractions.
In particular, we divide them into
exceptional and regular.
This implies the same in some special situations.

For instance, any finite (and even reductive) group
representation, or more generally, action
on a Fano or on an algebraic variety with
a numerically semi-negative canonical divisor,
can be treated as exceptional or regular,
and they have such an invariant as the index
in accordance with their quotients in 
Example~\ref{int}.\ref{gquot}.
In particular, this works for subgroups
of $PGL$'s.
In the case of $G\subset PGL(1)$ all the exceptional
subgroups correspond to exceptional
boundary structures on $\mathbb P^1\cong \mathbb P^1/G$.
Crystallographic groups give
other possible examples.

We can apply the same ideas to a classification
of surface quotients, log terminal or canonical, or elliptic
singularities, or their higher dimensional analogs.  
The same
holds for elliptic fiberings or other fiberings with
complements.

Quite possibly, 
the most important applications of
complements are still to come; perhaps
they will be related to classifications
of contractions of 3-folds $X$ semi-negative
with respect to $K$.
They may help to choose an appropriate
model for Del Pezzo, elliptic fiberings and
for conic bundles.
A strength of these methods is that
they apply in the most general situation when
we have log canonical singularities and
contractions are extremal in an algebraic sense,
or even just contractions.
This also shows a weakness, because the distance
to very special applications, when we have such
restrictions as terminal singularities and/or
extremal properties contractions,
may be quite long and difficult.
In addition, exceptional cases are
still not classified completely and explicitly.

A primitive sample we give in \ref{ecomp}.\ref{emainth}.3.
\end{remark}

Another application to a log uniruledness
is given in \cite{KM}.
(Cf. Remark~\ref{clasf}.\ref{rull}.)

Now we explain the statement of the Main Theorem
and outline a plan of its proof.
We also derive some corollaries.
Let $S$ be a normal algebraic surface,
and $C+B$ be a boundary (or subboundary) in it.
We assume that $C=\lfloor C+B\rfloor$ is
the reduced part of the boundary, and
$B=\{C+B\}$ is its fractional part.
Let $f:S\to Z$ be a contraction.
Fix a point $P\in Z$.
We need to find $n\in RN_2$, such that
$K+C+B$ has an $n$-complement 
locally over a neighborhood of $P\in Z$, 
or prove that
other possibilities are bounded.
More precisely, the latter ones are only global:
$Z=P$, and $(S,C+B)$ are bounded.
In particular, underlying $S$ are bounded.
On the other hand, for $(S/Z,C+B)$ having an $n$-complement,
there exists an element $D\in |-n K-n C-\lfloor(n+1)B\rfloor|$
such that, for $B^+=\lfloor
(n+1)B\rfloor/n+D/n$, $K+C+B^+$ has only log canonical
singularities (cf. \cite[Definition~5.1]{Sh2}).  In the local
case the liner system means a local one$/P$.

As in Corollary~\ref{int}.\ref{mcorol}
we can prove more.

\refstepcounter{subsec}
\label{nonvan1}
\begin{corollary}
Let $(S/Z,C+B)$ be a {\rm quasi\/} log Del Pezzo, { \rm i.~e.,
the pair with nef$/Z$ $-(K+C+B)$ and under (EC)\/}, and 
(M) holds for $B$.
Then for almost all of them
there exists an index $n\in RN_2$ such that
$|-n K-n C-\lfloor (n+1)B\rfloor|\not=\emptyset$
or, equivalently, we have a non-vanishing
$R^0f_*\mathcal O(-n K-n C-\lfloor (n+1)B\rfloor)\not=0$
in $P$.
More precisely, the same hold for $n\in RN_1=\{1,2\}$,
whenever $K+C+\lfloor (n+1)B\rfloor/n$ is not
exceptionally log terminal, {\rm i.~e.,
there exists an infinite set of exceptional
divisors$/P$ with the log discrepancy $0$.}

Moreover, we need the condition (M) only in
the global case with only Kawamata log terminal
singularities, in particular, with $C=0$, and
the exceptions for $RN_2$ belong to it.
In this case we can choose a
required $n$ in $N_2$.
\end{corollary}

The last two statements will be proven in
Section~\ref{gcomp} and
a more precise picture is given 
in Example~\ref{csing}.\ref{sabclas}.
Note also that $n$ in the local case
in the corollary depends on $P$.
We prove that $N_2$ is finite but still we have no
explicit description of it.
By the Monotonicity~\ref{indcomp}.\ref{monl1}
$\lfloor(n+1)B\rfloor\ge n B$
under the condition (SM), or even (M) if $n\in RN_2$.

\refstepcounter{subsec}
\label{nonvan2}
\begin{corollary}
Suppose again conditions (EC) and (M). 
Then for almost all of log pairs $(S/Z,B+C)$
there exists an index $n\in RN_2$ such that
$|-n (K+C+B)|\not=\emptyset$
or, equivalently, we have a non-vanishing
$R^0f_*\mathcal O(-n(K+C+B))\not=0$ in $P$.
More precisely, the same hold for $n\in RN_1=\{1,2\}$,
whenever $K+C+B$ is not exceptionally log terminal.

The exceptions for $RN_2$ belong only to the global case
with Kawamata log terminal singularities,
in particular, with $C=0$.
In this case we can choose $n$ in $N_2$
but in (M) instead of $b_i\ge 6/7$
we should require that $b_i\ge m/(m+1)$
with maximal $m\in N_2$.

In particular, if the boundary $C+B$ is reduced, i.~e.,
$B=0$, then $|-n (K+C)|\not=\emptyset$
or $R^0f_*\mathcal O(-n(K+C))\not=0$ in $P$.
\end{corollary}

The last non-vanishing is non-trivial even when
$S$ is a Del Pezzo surface with quotient
singularities, because it states that
$|-nK|\not=\emptyset$ or $h^0(S,-nK)\not=0$
for bounded $n$.
Such Del Pezzo surfaces form an unbounded family
and moreover the indices of $K$ for them
are unbounded, 
as are their ranks of the Picard group.

\refstepcounter{subsec}
\label{nonvan3}
\begin{corollary}
Again under (EC) and (M): 
A log pairs $(S/Z,B+C)$ is exceptional when
$Z=\pt$, $(S,B+C)$ is Kawamata log terminal,
in particular, $C=0$, and, for each $n\in RN_2$,
$|-n K-\lfloor (n+1)B\rfloor|=\emptyset$ or
$|-n (K+B)|=\emptyset$, equivalently, we have a vanishing
$h^0(S,-n K-\lfloor (n+1)B\rfloor)=0$ or
$h^0(S,-n(K+B))=0$.

In particular, if the boundary $B$ is reduced, i.~e.,
$B=0$, then $|-n K|=\emptyset$
or $h^0(S,-nK)=0$.
\end{corollary}

Quite soon, in the proof of the global case
in the Inductive Theorem, we see that an inverse
holds at least in the weak Del Pezzo case
and in the formal sense.
It leads to the following question.

\refstepcounter{subsec}
\begin{question}
Under (WLF) and (SM), does any $(X/Z,B)$
with a formal regular complement have
a real non-exceptional regular complement in
the sense of Definition~\ref{int}.\ref{dexc}?
In general it is not true completely
($\reg{(X,B^+)}=1$), for instance,
locally for non log terminal singularities with
$2$-complements of type $\mathbb E2_1^0$
(see Section~\ref{clasf}); the singularity
is exceptionally log terminal.
It holds for log terminal singularities
of types $\mathbb A_m$ and $\mathbb D_m$.
But it is unknown even for surfaces $X=S/\pt$.
\end{question}

The inverse does not hold if we drop (WLF),
for instance, for Abelian varieties,
as it was mentioned in Example~\ref{int}.\ref{tcomp}.
However, it looks possible for $\mathbb P^1\times E$
where $E$ is an Abelian, in particular, an elliptic curve.
In this case $-K$ has a positive numerical dimension.
Moreover, each $1$-complement in such a case has
log singularities.
Alas,
these complements are not quite non-exceptional as well.
We explain it in Section~\ref{csing} in terms
of $\reg{(S,B^+)}$, which
specifies the question for $1$ and $2$-complements
with $\reg{(S,B^+)}=1$.

The above corollaries show that it is easier
to construct a required
$D\in |-n K-n C-\lfloor(n+1)B\rfloor|$
with small $n$ when $K+C+B$ has more log singularities.
An expansion of this fact is related
to the Inductive Theorem in Section~\ref{indcomp}, an analog
of arguments in \cite[Theorem 5.6]{Sh2}.
As in the last theorem, under conditions of
the Inductive Theorem we extend or lift $D$
and a complement from
its $1$-dimensional restriction or projection.
So, we say that we have an {\it inductive\/}
case or complement, and the latter has
regular indices and is non-exceptional in the global case.

Then we may try to change $(X/Z,B)$ in such a way
that a type of complement will be preserved
and new $(X/Z,B)$ will satisfy the Inductive Theorem.
For instance, we can increase $B$ which will
be used for local complements in Section~\ref{lcomp}.
Here, in the global case, arises a problem
with standard multiplicities.
Fortunately, all the other cases are global
and Kawamata log terminal where we use
a reduction to the Inductive Theorem or to
a Picard number $1$ case which,
in addition, is $1/7$-log terminal in points.
In the latter case we apply Alekseev's \cite{Al}.
This we discuss in Section~\ref{gcomp}.

\section{Inductive complements}
\label{indcomp}

\refstepcounter{subsec}
\label{noncex}
\begin{example}~Let $S=\mathbb P_E(F_3)$ be
a ruled  surface over
a non-singular curve $E$ of genus $1$,
which corresponds to a non-splitting vector bundle
$F_3/E$ of rank 2 \cite[p.~141]{BPV}
with the odd determinant.
We denote the ruling by $f\colon S\to E$.
It has a single section $\mathbb P_E(F_1)$
in its linear class by a Riemann-Roch and
a vanishing below.
We identify the section with $E$.

Note that $E|_E\sim \det F_2\sim O$
for a single point $O\in E$
(cf. \cite[Proposition 2.9 in Ch.V]{H}).
So, we have a natural structure of
an elliptic curve on $E$ with $O$ as a zero.
On the other hand $(K_S+E)|_E\sim K_E\sim 0$
according to the Adjunction Formula.
Thus $(K_S+2E-f^*O)|_E\sim 0$ \cite[Lemma 2.10 in Ch. V]{H}.
But $K_S+2E\equiv 0/E$.
So, as in Example~\ref{int}.\ref{contex}, 
$K+2E-f^*O\sim 0$ and
we may chose $-K_S=2E-f^*O$ as an
anti-canonical divisor.
In particular, $K_S^2=0$, and
$-K_S$ is nef because $F_3$ is not splitting.

The latter also implies that
the cone $\NEc{S}$ is two-dimensional and
has two extremal rays:
\begin{itemize}
\item the first one $R_1$ is generated by
a fibre $F$ of the ruling $f$,
\item the second one $R_2$ is generated by
$-K_S$.
\end{itemize}
Since $-K_S=2E-f^*O$, ray $R_2$ has no
sections.
We contend that {\it $R_2$ has
three unramified double sections
$C_i$, $1\le i\le 3$.
Each of them is $C_i\sim 2E-
f^*(O+\theta_i)=-K_S-f^*\theta_i$,
where $\theta_i$ is a non-trivial element
of the second order in $\Pic{E}$.
Hence $K_S$ has $2$-complements $C_i$:
$2(K_S+C_i)\sim 0$ but $K_S+C_i\not\sim 0$,
but none $1$-complement.
The sections are non-singular curves of genus $1$.
Moreover, $R_2$ is contractible with
$\cont{R_2}:S\to \mathbb P^1$ having
$4$-sections also non-singular and of genus $1$
as generic fibres .
The curves $C_i$ are the only multiple
(really, double) fibres  of $\cont{R_2}$.\/}

Since $-K_S+m E$ is ample for integer $m>0$,
we have vanishings
$h^i(S,m E)=h^i(K_S-K_S+m E)=0$ for $i>0$
according to Kodaira.
Hence by the Riemann-Roch,
$$
h^0(S,2E)=2E(2E-K_S)/2=3.
$$
Similarly, $h^0(S,E)=1$.
On the other hand a restriction
of $|2E|$ on $E$ is epimorphic and free,
because $h^1(S,E)=0$ and the restriction
has degree $2$ on $E$.
Thus $|2E|$ is free and defines
a finite morphism $g:S\to\mathbb P^2$
of degree $4$.
According to the restrictions, $E$ goes
to a line $L=g(E)$, and a generic fibre 
$F=f^*P$ of $f$ goes to a conic $Q=g(F)$,
which is tangent to $L$ in $g(P)$.
Otherwise, $|2E-F|\not=\emptyset$, and
we have a single double section $C\equiv -K_S$
passing through $Q\sim 2O-P$.
This defines a contraction onto $E$
which is impossible by Kodaira's formula 
\cite[Theorem~12.1 in Ch. V]{BPV},
since $-K_S$ is nef.
Thus $Q=g(F)$ is a conic and
$g$ embeds $F$ onto $Q\subset\mathbb P^2$.
Then, locally over $g(P)$, by the projection formula,
$g^*(L|_Q)=(g^*L)|_F=2E|_F=2P$.
Thus $Q$ is tangent to $L$ in $g(P)$.
Note, that the same image gives a fibre $F'=f^*P'$
for $P'\sim 2O-P$, because $F_3$ is invariant
for an automorphism (even involution)
$i:S\to S$ induced by
the involution $P\mapsto P'$ on $E$.
The latter holds by a uniqueness of $F_3$
with the determinant $O$.

Therefore we have a $1$-dimensional
linear system $|4E-f^*O|$,
a generic element of which is
a non-singular $4$-section of genus $1$.
This gives $g=\cont{R_2}$.
Since all non-multiple fibres  of $g$
are isomorphic, then
$g_*K_{S/E}\equiv\chi(\mathcal O_S)=0$.
Hence by Kodaira's formula we have
degenerations.
(Equivalently, the second symmetric product
of $F_3$ is not spanned on $E$;
cf. \cite[p.~3]{Md}.)
Double fibres  are the only possible degenerations,
because their components are in $R_2$ but
not sections of $f$.
They give double sections $C_i$ of $f$ in $R_2$.
On the other hand $2K_S\sim -4E+2f^*O\sim-g^*P$
for generic $P\in \mathbb P^1$.
Thus $K_S\sim_{\mathbb Q} -g^*P/2$.
In particular, $K_S$ has $2$-complement $g^*P/2$.
Again due to Kodaira we have $3$ double fibres:
$C_i$, $1\le 1\le 3$.
Note that $2 C_i\sim g^*P\sim -2K_S\sim 4E-2f^*O$.
Thus we get all properties of $C_i$,
except for $\theta_i\not= 0$.
This follows from a monodromy argument
since we have exactly $3$ such $\theta_i \in \Pic{E}$
with $2\theta_i=0$,
and $3$ double sections $C_i$ of $f$ in $R_2$.
(The corresponding double cover $C_i/E$
is given by $\theta_i$.)
\end{example}

\refstepcounter{subsec}
\label{elcomp}
\begin{corollary}
Let $S/E$ be an extremal ruling over
a non-singular curve of genus $1$.
It can be given as a projectivization
$S=\mathbb P_E(F)$ for a vector bundle
$F/E$ of rank $2$.
Then $S/E$ or $F/E$ has
\begin{itemize}
\item a splitting type if and only if
$K_S$ has a $1$-complement;
\item an exception in
Example~\ref{indcomp}.\ref{noncex}
if and only if $K_S$ has $2$-complement; and
\item an exception in Example~\ref{int}.\ref{contex}
if and only if $K_S$ does not have complements at all.
\end{itemize}
In addition, the cone $\NE{S}$ is closed
and generated by two curves or extremal rays:
\begin{itemize}
\item the first ray $R_1$ is generated by
a fibre $F$ of the ruling $S/E$,
\item the second one $R_2$ is given by
one of curves $G$ with $G^2\le 0$
given by a splitting in a splitting case:
$K_S+G+G'\sim 0$; other cases were discussed
in Examples~\ref{int}.\ref{contex} and
\ref{indcomp}.\ref{noncex}.
\end{itemize}
\end{corollary}

\begin{proof}
We need only to consider the splitting case
when $F=V\oplus V'$ is a direct sum of
two line bundles.
Then we have two disjoint sections
$G=\mathbb P_E(V)$ and $G'=\mathbb P_E(V')$.
We may suppose that $G^2\le 0$.
By the adjunction formula,
$K_S+G+G'\equiv 0$, moreover, $\sim 0$
by arguments in Example~\ref{int}.\ref{contex}.
Thus $K$ has a $1$-complement, and
$G$ generates $R_2$.
The latter holds for some $C$ with $C^2\le 0$.
We assume that curve $C$ is
an irreducible multi-section of $f$.
If $C^2<0$, then $C=G$.
Otherwise $C$ is disjoint
from $G$ and $G'$,
and $(K_S.C)=(K_S+G+G'.C)=0=(K_S+G+G'.G)$.
Thus $C=G$.
If $C^2=G^2=0$, $R_2$ is generated by $G$ as well.
\end{proof}

\refstepcounter{subsec}
\label{indth}
\theoremstyle{plain}
\newtheorem*{indth}{\thesection.\thesubsec.
Inductive Theorem}
\begin{indth}
Let $(S/Z,C+B)$ be a surface log contraction
such that
\begin{description}
\item{{\rm (NK)}}
$K+C+B$ is not Kawamata log terminal,
for instance, $C\not=0$, and
\item{{\rm (NEF)}} $-(K+C+B)$ is nef.
\end{description}
Then it has
locally$/Z$ a regular complement, {\rm i.~e.,
$K+C+B$ has $1-,2-,3-,4-$ or $6$-complement\/},
under (WLF) of Conjecture~\ref{int}.\ref{conjcs},
or assuming (M)
under any one of the following conditions:
\begin{description}
\item{{\rm (RPC)}}
$\NEc{S/Z}$ is rationally polyhedral with
contractible faces$/Z$, or
\item{{\rm (EEC)}} there exists an
{\rm effective\/} complement, {\rm i.~e.,
a boundary $B'\ge B$ such that $K+C+B'$
is log canonical and $\equiv 0/Z$,} or
\item{{\rm (EC)+(SM)}}, or
\item{{\rm (ASA)}} {\rm anti\/} log canonical divisor
$-(K+C+B)$ {\rm semi-ample}$/Z$, or
\item{{\rm (NTC)}} there exists a
{\rm numerically trivial contraction\/}
$\nu : X\to Y/Z$, {\rm i.~e., $\nu$ contracts
the curves $F\subset S/Z$ with $(K+C+B.F)=0$.}
\end{description}
\end{indth}

\begin{statement1}
We can drop (M) in the theorem, but
then it states just a boundedness
of $n$-complements.
More precisely,
$n\in\{1,2,3,4, 5,6, 7,8, 9,10,
11,12, 13,14$, $15,16, 17,18, 19,20,
21,22, 23,24,25, 26, 27,28, 29,30,
31,35, 36,40,
41$, $42, 43$,
$56, 57\}$
\end{statement1}

\begin{statement2}
If there exists an infinite number
of exceptional divisors with
the log discrepancy $0$ then
we have a $1$ or $2$-complement.
If we drop (M), we have
in addition $6$-complements.
\end{statement2}

There exist (formally) non-regular
complements in the Inductive Theorem
when (M) is not assumed.
Similarly, we have examples
with only $6$-complements
in \ref{indcomp}.\ref{indth}.2.

\refstepcounter{subsec}
\label{nonrcomp}
\begin{example}
Let $f:S\to \mathbb P^1$ be
an extremal ruling $\mathbb F_n$
with a negative section $C$.
Take a divisor
$B=V+H$ with a vertical part
$f^*(E)$, where $E\ge 0$, and
a horizontal part $H=\sum d_i D_i$
with $d_i\in \mathbb Z/m\cap (0,1)$ for
some natural number $m$.
The latter can be given, for instance, as
generic sections of $f$.

Suppose that the different $B_C=(B)_C$ has
only the standard multiplicities,
$V=f^*B_C$,
$\deg (K_C+B_C)=0$, and
$H$ disjoint from $C$.
Then $K+C+B$ as $K_C+B_C$ have
$n$-complements only for $n$ such
that $n$ is divided by the index of $K_C+B_C$
(cf. Monotonicity Lemmas~\ref{indcomp}.\ref{monl1}
and~\ref{indcomp}.\ref{monl3} below).

Then we do not have $n$ complements for $m=(n+1)$
when $K+C+B\equiv 0$
(cf. Monotonicity~\ref{indcomp}.\ref{monl2}.1
below).

Therefore we have a non-regular $n(n+1)$-complement
when $K_C+B_C$ has only
one regular $n$-complement, and
$n(n+1)\ge 7$.
Moreover, we can have no regular complements
in this case.
For instance, this holds, when $n=6$.

In general non-regular complements
in the Inductive Theorem
have a similar nature, as we see in its proof,
with the following modifications.
The ruling may not be extremal,
$S$ may be singular,
and $B_C$ may have non-standard multiplicities.
To find complements in such cases,
we use Lemmas~\ref{indcomp}.\ref{indcomp1},
\ref{indcomp}.\ref{compl2}-\ref{compl3}.

If we take $V$ with the horizontal multiplicities $1/3$,
and such $H$ that $(C,B_C)$ has just $2$-complements,
then $K+C+B$ will have $6$-complements as the minimal.
\end{example}

The following result clarifies
relations of conditions in
the Inductive Theorem.

\refstepcounter{subsec}
\label{relcond}
\begin{proposition}
Assuming that $K+C+B$ is log canonical
and nef$/Z$,
$$
{\rm (WLF)}\Longrightarrow 
{\rm (RPC)}\Longrightarrow
{\rm (NTC)}\Longleftrightarrow
{\rm (ASA)}\Longleftrightarrow
{\rm (EEC)}\Longleftarrow
{\rm (EC)}+{\rm (SM)}
$$
with the following exception.
For
${\rm (WLF)}\Longrightarrow{\rm (RPC)}$:
\begin{description}
\item{\rm (EX1)}
contraction $f:S\to Z$ is birational,
and up to a log terminal resolution,
$C=C+B$ is a curve with nodal singularities
of arithmetic genus $1$,
$/P$, $K+C\equiv 0/Z$, and $S$ has
only canonical singularities and outside $C$.
\end{description}
For
${\rm (EEC)}\Longrightarrow{\rm (NTC)}$:
\begin{description}
\item{\rm (EX2)} $Z=\pt$,
$K+C+B$ has a numerical dimension $1$,
$B'$ and $E$ are unique,
$(K+C+B')$ has a log
{\rm non-torsion\/} singularity of
{\rm genus\/} $1$ and of
{\rm numerical dimension\/} $1$, {\rm i.~e.,
on a crepant log resolution $C+B'=C$ form
a curve with only nodal singularities,
with the connected components of genus $1$,
and $E|_{C}\equiv 0$,
but not $\sim_{\mathbb Q}$.}
\end{description}
Nonetheless,
${\rm (WLF)}\Longrightarrow{\rm (NTC)}$
always, and
${\rm (WLF)}\Longrightarrow{\rm (RPC)}$
always in the analytic category or
in the category of algebraic spaces.
In (EX2) there exists a $1$-complement.
\end{proposition}

\refstepcounter{subsec}
\begin{remcor}
\label{rel}
In particular, for surfaces,
(WLF) always implies (ASA).
In other words, if $-(K+C+B)$ is log canonical,
nef and big$/Z$, then it is semi-ample.
It is well-known when $K+B+C$ is Kawamata
log terminal \cite[Remark~3-1-2]{KMM}
(cf. with arguments in the proof below).
For log canonical singularities,
it is an open question in dimension $3$ and higher.
In dimension $2$ we can replace
the last two conditions 
by a non-vanishing
\begin{description}
\item{(NV)}
$-(K+C+B)\sim_{\mathbb R} E/Z$ where 
$E$ is effective,
\end{description}
and $E^2>0$.
(For a definition of $\sim_{A}$ see
\cite[Defenition~2.5]{Sh3}.)
Note also that (EEC) implies (NV) but
not vise-versa (cf. Example~\ref{int}.\ref{contex}).
So, in general nef $-(K+C+B)/Z$ is not semi-ample if
it is not big and not $\equiv 0/Z$.
In the latter case $E=0$ and
according to the semi-ampleness conjecture
for $K+C+D/Z$ we anticipate semi-ampleness
\cite[Conjecture~2.6]{Sh3} 
(cf. \cite[Remark~6-1-15(2)]{KMM}
and Remark~\ref{int}.\ref{tcomp}).
This is the main difficulty in a construction
of complements: (EC).
In dimension $2$, at least
(NV) {\it holds when $-(K+C+B)/Z$ is nef\/}.
However as in Example~\ref{int}.\ref{contex}
$K+C+B+E$ may not be log canonical.
Does such a non-vanishing hold in higher dimensions?
In any case it implies
a log generalization of a Campana-Peternell's problem
in dimension $2$ \cite[11.4]{CP}, \cite{Md} and
cf.~\cite{Gr}: $-(K+C+B)$
is ample$/Z$ if and only if
$-(K+C+B)$ is positive on all curves $F\subset Y/Z$.
Indeed, then $E^2>0$ and $-(K+C+B)$ is ample due
to Nakai-Moishezon~\cite[Corollary 5.4]{BPV}, or
we can use the implication
${\rm (WLF)}\Longrightarrow{\rm (ASA)}$.
Therefore we have a complement in
a weak form (EEC).
Again by Example~\ref{int}.\ref{contex}
we cannot replace the above positivity be a weaker version:
nef and $(K+C+B.F)=0$ only for a finite set of curves.
As we will see in a proof of (NV),
the nef property of $-(K+C+B)$ implies
(EEC) in most cases and in a linear form.
For instance, the only exception is $S$
of Example~\ref{int}.\ref{contex} when
$S$ is not rational.
Do exceptions exist when $S$ is rational?
It appears that $-(K+C+B)$ satisfies (ASA)
in most of cases as well.
The exceptions again are unknown to the author.
For $3$-folds, similar questions are more difficult,
for example, the Campana-Peternell's problem.
In dimensional $2$, the most difficult cases
are related to non-rational or
non-rationally connected surfaces,
or more precisely,
to extremal fiberings over curves $E$ of genus $1$.
They are projectivizations of rank $2$
vector bundles$/E$.
So, similar cases are of prime interest
for $3$-folds:
projectivizations $X$ of rank $3$ vector bundles$/E$,
of rank $2$ over Abelian or K3 surfaces.
Is their cone $\NE{X}$ always closed
rationally polyhedral and
generated by curves as for Fanos?
What are the complements for $K$?
(Cf. Corollary~\ref{indcomp}.\ref{elcomp}.)

According to the same example,
we really need at least a contractibility
of extremal faces in (RPC) for (EEC) or (EC),
moreover for regular complements in the theorem.
As in (EX1), we sometimes have just
the rational polyhedral property but
not contractibility of any face.
If $-(K+C+B)$ is only nef, the cone
may not be locally polyhedral near
$(K+C+B)^\perp$ \cite[Example 4.6.4]{CKM}.
The same example gives an exception of (EX2).

In dimension $\ge 3$,
(ASA) is better than (NTC): (ASA) implies (NTC)
easily, but the converse is harder
as we see below.
In addition,
in a proof of the Inductive Theorem below,
we will see
that a semi-ampleness
of $-(K+B)$ on $S$ is enough.
Fortunately, it is good for an induction in
higher dimensions.
This will be developed elsewhere. 
It appears that the same should hold for (NTC);
then it would be best in this circumstance.
(NTC) also implies the Campana-Peternell's problem
(but not a solution of the latter).
\end{remcor}

\begin{proof} First we check that (WLF) implies
(NTC), and, except for (EX1), (RPC).
For the former, it is enough to prove (ASA).
For the latter, it is enough to prove that
$\NEc{S/Z}$ is rationally polyhedral
and, except for (EX1),
with the contractible faces
near a contracted face $\NEc{S/Z}\cap (K+C+B)^\perp$.
Indeed,
the cone satisfies both properties locally outside the face
by the LMMP.
A contraction in any face of
$\NEc{S/Z}\cap (K+C+B)^\perp$ preserves
this outside property.

To prove (ASA) in the $\mathbb Q$-factorial case
we may use the modern technique:  
to restrict a Cartier multiple
of $-(K+C+B)$ on the log singularities of $K+C+B$.
Then a
multiple of $-(K+C+B)$ will be free in such singularities,
because any nef Cartier divisor is free on a point, on a
rational curve, and on a curve of arithmetic genus $1$ with at
most nodal singularities assuming that in the last case the
divisor is canonical when it is numerically trivial.  
We meet only this
case after the restriction.
Then we may apply
traditional arguments to eliminate the base points outside the
log singularities.

For surfaces, however we may use more direct arguments,
which work in the positive characteristics. 
Since $-(K+C+B)$ is big we have a finite
set of curves $F/P$ with $(K+C+B.F)=0$.
We check that each of them generates an extremal ray
and (RPC) or (EX1) hold near them.
After a log crepant blow-up \cite[Theorem~3.1]{Sh3},
we assume that $S$ has only rational singularities
where as $K+C+B$ is log terminal.
In particular, $C$ has only nodal singularities.

In most cases we can easily check (RPC) and
derive (ASA) from this as shown below.
For instance, if $-(K+C+B)$ is ample it follows from the LMMP.
Otherwise we change $C+B$ such that it will hold
except for three cases: (EX1), or
\begin{description}
\item{\rm (i)} $C$ is a chain of rational curves$/P$,
$(K+C+B.C+B')=0$ in $C+B'$ is contractible to
a rational elliptic singularity, where
$C+B'/P$ is the connected component of $C$ in $B$; or
\item{\rm (ii)}
$Z=\pt$ and there exists
a ruling $g:S\to C$ with a section $C$, with $C^2<0$,
which is a non-singular curve of genus $1$, and
$S$ has only canonical singularities, and
outside $C$ and $B$,
$B$ has no components in fibres 
and does not intersect $C$.
\end{description}
In (EX1) we obviously have 
the rational polyhedral property but
contractions of some faces may not exist
in general.
In the last two cases we
check (RPC) directly near $(K+C+B)^\perp$.

According to \cite[Lemma~6.17]{Sh3},
we suppose that $-(K+C+B)=D=\sum d_i D_i$
is an effective divisor with irreducible curves $D_i$.
Moreover, we assume that $\Supp D$ contains
an ample divisor $H$.
So, if $K+C+B+\delta D$ is log canonical for some $\delta>0$,
then $K+C+B+\delta D-\gamma H$ gives
a required boundary $C+B\colon =C+B+\delta D-\gamma H$.
For instance, it works when $K+C+B$ is
Kawamata log terminal.
So, $C\not=0$ in the exceptional cases
which are considered below.
The curve $C$ is connected
by \cite[Lemma~5.7 and the proof of Theorem~6.9]{Sh2}
(see also  \cite[Theorem~17.4]{KC}).

More precisely,
$K+C+B+\delta D$ is not log canonical whenever
$\Supp D$ has a component $C_i$ which is in
the reduced part $C$ as well
(cf. \cite[1.3.3]{Sh2}).
Moreover, any component $C_i$
of $C$ is $/P$, in $D$, $(K+C+B.C_i)=0$ and $C_i^2<0$.
If $C_i$ is not$/P$ or $C_i/P$ with
$(K+C+B.C_i)<0$, we can decrease
the boundary multiplicity in $C_i$.
By an induction and connectedness of $C$
this gives a reduction to the Kawamata log terminal case.
If $C_i^2\ge 0$ we can do the same.
Since $H$ and $D$ intersect $C_i$, $C_i\subseteq\Supp{D}$.

In most cases,
$C$ is a chain of non-singular rational curves,
because $C$ is connected with nodal singularities
and
$$\deg K_C\le (C.K+C)\le (C.K+C+B)=-(C.D)\le 0.$$
The only exception arises when
$\deg K_C=0$, $C$ has the arithmetic genus 1,
$B=0$, and $S$ is non-singular in a neighborhood of $C$.

To check (RPC) locally,
we need to check (in the exceptional cases)
that $D=-(K+C+B)$, that
each nef $\mathbb R$-divisor $D$ near $-(K+C+B)$
is semi-ample, and that only curves $F/Z$
with $(K+C+B.F)=0$ are contracted by $D$.
Indeed, $D$ is big as $-(K+C+B)$.
So, we may assume that $D$ is effective as above
and $D\ge H$.
So, $D$ is numerically trivial only on a finite
set of curves $F$ and $F^2<0$.
We may also assume that $D$ is quite close to
$-(K+C+B)$, namely, $(K+C+B.F)=0$.
In other words we need to check that
any set of curves $F$ with $(K+C+B.F)=0$
is contractible in an algebraic category.
If the singularities after the contraction
are rational we can pull down $D$ to a
$\mathbb R$-Cartier divisor.
So, the pull down of $D$ is ample
according to Nakai-Moishezon \cite[Corollary~5.4]{BPV}.
Hence $D$ is semi-ample.
Of course, it applies directly when $D$ is
a $\mathbb Q$-divisor.
Otherwise we may present it as weighted combination
of such divisors 
(cf. \cite[Step~2 in the Proof of Theorem~2.7]{Sh3}),
because a small perturbation of coefficients
of $D$ preserves it positivity on all curves.
(Essential ones are components of $\Supp D$.)
The cone is rational finite polyhedral near
$K+C+B$ because the set of curves $F$ is finite.

If the curves $F$ with exceptional curves
on a minimal resolution form a tree of rational curves,
we could then contract any set of curves $F$ by
Grauert's and Artin's criteria
\cite[Theorems~2.1 and 3.2]{BPV}.
In addition, the singularities after contraction
are rational log canonical.

It works in case (i) when $C$ is
a chain of non-singular rational curves $C_i/P$.
Indeed, we can do contractions inductively,
because $F^2<0$.
First, for $F$ not in $C$.
Since the boundary multiplicity in $F$
is $<1$, $F$ is a non-singular rational curve.
And by a classification of
the log terminal singularities, $F$ with the curves
of a resolution form a tree of rational curves.
So, $F$ is contractible (by the LMMP as well).
Moreover, this preserves the log terminality
and rationality of singularities.
By \cite[Lemma~5.7]{Sh2} $C$ will form a chain
of rational curves again.
Finally, $C$ will be the whole curve where
$K+C+B\equiv 0/Z$.
Since all singularities are log terminal,
we have the required resolution of $C$.
Then we can contract $C$ inductively too.
This effectively gives case (i) because
$B=0$ near $C$.
Otherwise we can decrease $B$ and replace
$C+B$ by a Kawamata log terminal boundary.

Now we suppose that $C$ is a curve$/P$
with only nodal singularities of arithmetic
genus $1$.
If $K+C+B\equiv 0/P$, then we have
the exceptional case (EX1),
because $-(K+C+B)$ is big$/P$.
Since $S$ is non-singular in a neighborhood of $C$,
it easy to check that other curves$/P$ on
a minimal resolution form a chain of rational curves
intersecting simply $C$.
Moreover, they are $(-1)$ or $(-2)$-curves.
So, we can contract any set of such curves
with any proper subset of $C$ due to Artin.
Here the only problem is to contract the whole $C$.
Alas, this is not always true in the algebraic category,
e.~g., after a monoidal transform
in a generic point of $C$.
The cone in this case is rational polyhedral,
but contractions may not be defined for some faces
including the components of $C$.
Nonetheless $-(K+C+B)$ is semi-ample
because $K+C+B$ is semi-ample in this case.

Suppose now that $K+C+B\not\equiv 0/Z$.
Then $C$ is non-singular and it is case (ii).
Indeed we have an extremal contraction
$S\to Y/Z$ negative with respect to $K+B$.
It is not to a point or onto a curve because
$(K+C+B.C)=(K+C.C)=0$, and
respectively $C$ cannot be a section,
if $C$ is singular.
So, the contraction is birational.
By a classification of such contractions it
contracts a curve which does not intersect $C$.
After that we have again $K+C+B\not\equiv 0/Z$.
Induction on the Picard number$/Z$ gives a contradiction.
If $C$ is non-singular the only possible case
is when, after a finite number of birational
contractions, we have an extremal ruling
which induces a ruling $g$ as in (ii) with
section $C$.
In particular then $Z=\pt$.

$S$ has only canonical singularities and
outside $C$.
There are also no components of
$\Supp B$ in fibres  of $g$ and intersecting $C$.
In other words, $K+C+B$ is canonical
in points, and there are no components of
$B$ in fibres  of $g$.
Indeed, as we know there are no components of $B$
in fibres  of $g$ which intersect $C$.
After extremal contractions$/C$ it
holds for any component in fibres  of $g$.
So, after a blow-up in a non-canonical point
we get a contradiction.
This implies that it is really case (ii).

Now we verify (RPC) near $(K+C+B)^\perp$.
As above in case (i) we need to check
that any set of curves $F$ with $(K+C+B.F)=0$
is contractible.
Again, as in case (i),
any such $F$ besides $C$ is rational.
Hence such $F$ belongs to a fibre of $g$.
So, we have a finite number of them.
As above this implies
the polyhedral property for $\NEc{S}$.
In addition, a required contraction
corresponds to a face in
$\NEc{S}\cap (K+C+B)^{\bot}$,
and the latter is generated by curves in
fibres  of $g$ and, perhaps, $C$.

We only want to check the existence of
the contraction in this face.
If $C$ does not belong to such a face,
it follows from the relative statement$/C$.
Otherwise after making
contractions of the curves of the face
in fibres , we suppose that $F$ is generated by $C$.
Any birational contractions, with
disjoint contracted loci, commute.
So, it is enough to establish a contraction
of $C$, after contractions of the curves
in fibres  of $g$, which does not intersect $C$.
Equivalently, we can assume that
the fibres  of $g$ are irreducible.
Then $S$ is non-singular,
$g$ is extremal,
and $\rho (S)=2$,
where $\rho(S)$ denotes the Picard number.
Since $C^2<0$,
we contract $C$ by $h\colon S\to S'$ to a point
at least in the category of normal algebraic spaces.
However, we can pull down $K+C+B$ to a $\mathbb R$-Cartier
divisor $-D'=K_{S'}+h(C+B)$
because $(K+C+B)|_{C}$ is semi-ample
over a neighborhood of $h(C)$.
This follows from an existence of
$1$-complement locally
(cf. \cite[Corollary~5.10]{Sh2})
and even globally as we see later
in the big case of the Inductive Theorem.

Since $D'$ is nef and big, and $\rho(S')=1$,
then $D'$ is ample and $-(K+C+B)=h^*D'$ is semi-ample.
This completes the proof of
${\rm (WLF)}\Longrightarrow{\rm (RPC)}$.

(RPC) implies (NTC) by the definition.

(ASA) implies (NTC), because
any semi-ample divisor $D$ defines
a contraction which contracts
the curves $F$ with $(D.F)=0$.
A converse and other arguments in this proof
are related to the semi-ampleness of
log canonical divisors~\cite[Theorem~11.1.3]{KC}
(it assumes $\mathbb Q$-boundaries that can
be improved up to $\mathbb R$-boundaries as
in \cite[Theorem~2.7]{Sh3}).
So, let $\nu :X\to Y/Z$ be a numerical
contraction.
Then, by the semi-ampleness$/Z$,
$K+C+B=\nu^*D$ for a $\mathbb R$-divisor on $Y$
which is numerically positive on each curve
of $Y/Z$.
Then it is easy to see that $D$ is ample$/Z$ and
hence semi-ample$/Z$ except for the case
when $\nu$ is birational and $Z=\pt$.
However, this case is also known because for
complete surfaces $-(K+C+B)$ is ample when
$(K+C+B.F)<0$ on each curve $F\subset S$
which follows from (NV) as explained
in Remark-Corollary~\ref{indcomp}.\ref{rel}.

In the corollary we need only to verify (NV)
whenever $-(K+C+B)$ is nef$/Z$.
We do it now.
After a crepant resolution we
assume that $K+C+B$ is log terminal.
When $-(K+C+B)$ is numerically big, then (NV) follows
from \cite[Lemma~6.17]{Sh2}.
On the other hand if the numerical dimension
of $-(K+C+B)$ is $0/Z$ then by
\cite[Theorem~2.7]{Sh3}
$-(K+C+B)\sim_{\mathbb R} 0/Z$.
In other cases, $Z=\pt$ and
$D=-(K+C+B)$ has the numerical dimension $1$,
i.~e., $D$ is nef, $D\not\equiv 0$ and $D^2=0$.
We assume that (NTC) and (ASA) do not hold.
Otherwise $-(K+C+B)$ is semi-ample as we
already know and (NV) holds.
Now we reduce (NV) to the case when $D$ is
a $\mathbb Q$-divisor, or, equivalently,
$B$ is a $\mathbb Q$-divisor (cf. proof
in the big case of the Inductive Theorem).
If $C+B$ has a big prime component $F$,
then $-(K+C+B-\varepsilon F)$ satisfy (WLF)
and $-(K+C+B)$ is not assumed (NTC).
On the other hand, each component $F$ of $B$
with $F^2<0$ can be contracted by the LMMP
or due to Artin.
If this contraction is crepant, we preserve
the numerical dimension, the log terminally
and (NV).
Otherwise $(K+C+B.F)<0$ and we can decrease
the multiplicity of $B$ in $F$ which, as
in the above big case, gives (NV) by (WLF).
So, we assume that each $F$ in $B$ with
$F^2<0$ is contracted, i.~e., for each
prime component $F$ in $B$, $F^2=0$.
By the same reason such components
$F$ are disjoint.
Hence, if we slightly decrease
each irrational $b_i$
to a rational value, we get
a $\mathbb Q$-boundary $B'\le B$
and a divisor $-(K+C+B')$ which is nef
and with the same numerical dimension.
Assuming (NV), but not (NTC), we have
a unique effective $D'\sim_{\mathbb Q} -(K+C+B')$.
By the uniqueness $D$ has positive
multiplicities in components $F$ for
any such change.
Moreover they are bigger or equal to the change.
So, $D\sim_{\mathbb Q} -(K+C+B)$ and (NV) holds.

We assume now that $D$ is a $\mathbb Q$-divisor.
After a minimal crepant resolution we
suppose also that $S$ is non-singular.
Then, for positive Cartier multiples $m D$,
we have vanishing $h^2(S,m D)=h^0(S,K-m D)=0$,
at lease for $m\gg 0$,
because $D$ has positive numerical dimension.
Therefore by the Riemann-Roch:
\begin{align*}
h^0(S,m D)\ge& h^1(S,m D)+m D(m D -K)/2+
\chi(\mathcal O_S)\\
=&-m D K/2+\chi(\mathcal O_S)=
m D(C+B)+\chi(\mathcal O_S)\ge
\chi(\mathcal O_S).
\end{align*}
In particular, $h^0(S,m D)\ge 1$
when $\chi(\mathcal O_S)\ge 1$.
For instance, when $S$ is rational.
Assuming now that $S$ is non-rational,
we check then that, after crepant blow-downs,
$S$ is an extremal or minimal ruling over
a non-singular curve of genus $1$.
Indeed, we have a ruling $g:S\to E$ over
a non-singular curve of genus $1$ or higher,
since $D=K+C+B\not=0$, and $D\not=0/E$.
If $g$ is not extremal we have
a divisorial contraction in a fibre of $g$,
which is positive with respect to $D$.
Hence after a contraction we have (WLF).
Since $S$ is not rationally connected
it is only possible when $C\not=0$ by
\cite[Theorem~9]{Sh6}.
If $C\not=0$, a classification of
log canonical singularities and
\cite[Corollary~7]{Sh6}
implies that we have a connectedness by
rational and elliptic curves in $C$.
Hence $E$ has genus $1$, and
$C$ is a section of $g$.
If $S$ then, is not extremal we have
a crepant blow-down of a $(-1)$-curve
in a fibre of $g$.

In addition, $\chi(\mathcal O_S)=0$
and, for each $F$ in $C+B$,
$(D.F)=0$ by the above Riemann-Roch.
In particular, $C+B$ does not have
components in fibres .
The cone $\NE{S}=\NEc{S}$ is a closed angle
with two sides:
\begin{itemize}
\item $R_1$ is generated by a fibre of the ruling $g$,
\item the second one, $R_2$, is generated by
a multi-section $F$ with $F^2\le 0$.
\end{itemize}
Note that $(D.R_1)>0$ and $D=-(K+C+B)$ is nef.
Hence, if $F^2<0$, we can take $F=C$ as a section.
In this case, we find an effective divisor $D=B'-B$
increasing a component of $B\not=0$, or taking
another disjoint section $F$ otherwise.
Then $K+C+B+D=K+C+B'\sim_{\mathbb R} 0$ by
the semi-ampleness \cite[Theorem~2.7]{Sh3},
because $K+C+B'$ is log canonical
\cite[Theorem~6.9]{Sh2}.
Finally, if $F^2=0$, then $D=-(K+C+B)$, and
any component of $C+B$ generates $R_2$ as well.
Moreover, each curve in $R_2$ is non-singular
and of genus $1$.
They are disjoint in $R_2$.
There exist at most two of them when (NTC)
is not assumed.
As above $D=B'-B$, for some $B'>B$,
such that $K+C+B'\equiv 0$.
Using \cite{At},
it is possible to check that $K+C+B'\sim 0$
for an appropriate choice or
$K+C+B'$ is log canonical,
and log canonical,
except for the case in
Example~\ref{int}.\ref{contex}.
In the latter case, we have
a single curve $E$ in $R_2$,
$C+B=a E$ with $a\in [0,1]$, and $K\sim -2E$.
Hence $D=-(K+C+B)\sim (2-a)E>0$.
This completes the proof of (NV).

Now suppose (EEC): there exists an effective divisor
$E=B'-B=(K+C+B')-(K+C+B)\equiv -(K+C+B)/Z$ which is nef.
If $E$ is numerically big$/Z$,
then $E$ is semi-ample$/Z$.
This implies (NTC).
In other cases the numerical dimension of $E$
is $1$ or $0/Z$.
In the latter case $K+C+B\equiv -E=0/Z$
is semi-ample.
This implies (ASA) and (NTC).
Hence we assume that $E$
has the numerical dimension $1/Z$ and $Z=\pt$.
Then we can reduce this situation to
the case when $E$ is
an isolated reduced component in $C+B'$.
Otherwise $E$ is contractible by the LMMP.
Indeed, if $K+C+B'$ is purely log terminal
near a  connected component $F$ of $E$,
we may increase $C+B'$ in $F$:
for small $\varepsilon>0$,
a log canonical divisor
$K+C+B'+\varepsilon F\equiv \varepsilon F$
is log canonical and semi-ample.
This implies (NTC).
After a crepant blow-up,
we can suppose that
$K+C+B'$ is log terminal, and
a log singularity
is on each connected component $F$ of $E$.
Then essentially by Artin we can
contract the log terminal components in $F$.
So, $E$ is reduced in $C+B'$, and
each connected component $F$ of $E$ has
the numerical dimension $1$.
We check that $F$
is semi-ample when $F$ is not isolated in $C+B'$.
This implies (NTC).
Let $D$ be another component of $B'$
intersecting $F$ and with a positive multiplicity $b$ in $B$.
Then we can subtract a nef and big positive linear combination
$\delta D+\varepsilon F$ from $K+C+B'$
which gives the nef and big anti log divisor
$-(K+B'')=-(K+C+B')+\delta D +\varepsilon F$.
That is (WLF), which implies (RPC) and (NTC).

Hence we need to consider now the case when
$C+B'$ has as a reduced and isolated component
$F=\Supp{E}$.
We also assume that $K+C+B'$ is log terminal.
Then $F$ has only nodal singularities.

Since $E$ is nef and of numerical dimension $1$,
$E|_F\equiv 0$.

If $E|_F\not\sim_{\mathbb R} 0$,
we have (EX2), but not (NTC).
Nonetheless in these cases $K+C+B'\sim 0$ or,
equivalently, $K+C+B$ has $1$-complement.
In particular, $C+B'$ is reduced.
Indeed, if $E|_F\not\sim_{\mathbb R} 0$,
then $F$ is a curve
of arithmetic genus $1$.
To find the index of $K+C+B'$
we may replace $S$ by its terminal resolution,
when $S$ is non-singular.
If $S$ is rational we reduce the problem to
a minimal case, where $F=C=C+B'\sim -K$.
Hence $K+C+B'$ has index $1$.
Otherwise $S$ is not rationally connected.
Hence $F$ is a non-singular curve
with at most two components of genus $1$,
$S$ has a ruling $g:S\to G$ over
a non-singular curve $G$ of genus $1$.
If $g$ has a section $G$ in $F$,
then $C$ has another section $G'$,
whereas $C=C+B'=G+G'$ and $K+C+B'=K+C\sim 0$.
Indeed, after contractions of $(-1)$-curves,
which do not intersect a component of $F$,
where $E|_F\equiv 0$, but $\not\sim_{\mathbb R} 0$,
we preserve the last property on an extremal $g$.
It has section $G$ in $F$ such that
$E|_G\not\sim_{\mathbb R} 0$ or,
equivalently, $G|_G\not\sim_{\mathbb Q} 0$.
According to [At], this is a splitting case, i.~e.,
$G$ generates extremal ray $R_2$, and
it has only two curves $G$ and another section $G'$.
The latter holds because $R_2$ is not contractible.
So, we can assume that $F=C+B'$ is a double section
of $g$.
This is possible only when
the ruling $g$ is minimal.
Otherwise after a contraction of a $(-1)$-curve
in a fibre of $g$ intersecting $F$, we get
big $F\equiv -K$.
It is impossible by (ASA)
for $-K$, since $K$ has no log singularities
and hence $S$ is rationally connected.
So, $S/G$ is extremal with another extremal ray
$R_2$ generated by $F$, and $F^2=0$.
According to a classification of minimal rulings
over $G$, this is a non-splitting case,
because (NTC) is not assumed.
By Examples~\ref{int}.\ref{contex} and
\ref{indcomp}.\ref{noncex}  this is
impossible for the other rulings: either we have
no complements, or we have (NTC)
by Corollary~\ref{indcomp}.\ref{elcomp}.

Note that (NTC) does not hold in (EX2),
for example, due to uniqueness of $B'$ and $E$.
In these cases $K+C+B'\sim 0$,
but not only $\equiv 0$
by Corollary~\ref{indcomp}.\ref{elcomp}.
In addition, $B'$ and $E$ are unique in (EX2),
because otherwise a weighted linear combination
of different complements $B'$ gives
a complement $C+B'$ which has no log singularity
in a component of $F$.

In other cases $E|_F\sim_{\mathbb R} 0$.
Let $G$ be a connected component of $E$
in a reduced part of $C+B'$.
We prove then that a multiple of $G$ is
movable at least algebraically.
This implies (NTC), because
$K+B+C'+\varepsilon G'$ is semi-ample for
small $\varepsilon>0$ and
a divisor $G'$ disjoint from $\Supp{(C+B')}$,
including $E$, which is algebraically
equivalent to a multiple of $G$.
This implies $\equiv$ as well.
Note, that the LMMP implies, after
log terminal contractions for $K+C+B'$, that
$C+B'$ is reduced, or (NTC) holds.
Indeed, if $B'$ has a prime component $D$ with
multiplicity $0<b<1$, then $D$ is in $B$ and
disjoint from $F$.
After log terminal contractions
of such divisors with $D^2<0$,
we assume that $D^2\ge 0$ for others.
Since $D^2>0$ implies (WLF) and (RPC) for
$K+C+B'-\varepsilon D$, then $D^2>0$
implies (NTC) for $E$.
The same holds for $D^2=0$
by the semi-ampleness of
$K+C+B'+\varepsilon D$.
Hence $C+B'$ is reduced after contractions.
Now we can use a covering trick
\cite[Example~2.4.1]{Sh2}, because
$K+C+B'\sim_{\mathbb Q} 0$.
If $S$ is rational, the latter holds for any
numerically trivial divisor.
If $S$ is non-rational, then as in (EX2) above
we check that $K+C+B'\sim 0$ or
$\sim_{\mathbb Q} 0$ of index $2$ by
Corollary~\ref{indcomp}.\ref{elcomp}.
So, after an algebraic covering,
which is ramified only in $C+B'$ we
suppose that $K+C+B'$ has index $1$.
A new reduced boundary is an inverse image
of $C+B'$.
The same holds for $G$.
After a crepant log resolution
we assume that $S$ is non-singular
and $G$ is a Cartier divisor.
The curve $G$ has only nodal singularities
and genus $1$.
Contracting the $(-1)$-curves in $G$,
we assume that each non-singular rational component
of $G$ is a $(-m)$-curve with $m\ge 2$.
Then all such curves are $(-2)$-curves and
we can take reduced $G$.
The latter is obvious in other cases, too.
Indeed, after a normalization, we suppose that
$D$ is a $\mathbb Q$-divisor with multiplicities
$\le 1$ and with one component $D$
with multiplicity $1$.
Then $K+C+B'-G$ is log terminal near $G$,
$\equiv 0$ on $G$
and has multiplicity $0$ in $D$.
Since each such $D$ is not a $(-1)$-curve,
log divisor $K+C+B'-G$ is trivial near $G$
and $G$ is a reduced Cartier divisor.
According to our assumptions and reductions,
$G|G\sim_{\mathbb Q} 0$.
We prove that a multiple $mG$ is linearly movable.
First, we suppose that $S$ is rational.
Let $m G|G\sim 0$ for an integer $m>0$.
Then by a restriction sequence on $G$,
we have a non-vanishing
$h^1(S,(m-1)G)\not=0$ whenever $G$ is
linearly fixed and the restriction is
not epimorphic.
Hence by the Riemann-Roch,
\begin{align*}
h^0(S,(m-1)G)&\ge
h^1(S,(m-1)G)+(m-1)G((m-1)G- K)/2+\chi(\mathcal O_S)\\
&\ge 1+0+1=2
\end{align*}
and $(m-1)G$ is linearly movable.
In other cases, $S$ is not rationally connected.
Moreover, $F$ is a non-singular curve of genus $1$
or a pair of them, and there exists
a ruling $g: S\to G'$ to
a non-singular curve $G'$ of genus $1$.
If $G$ is a section of $g$, then we can
suppose that $S$ is extremal$/G'$ after contractions
of $(-1)$-curves, disjoint from $G$, in fibres $/G'$.
By \cite{At}, (EEC) and
Example~\ref{int}.\ref{contex},
it is a splitting case, i.~e., $g$ has
a section $G'$ in $C+B'$,
and $G'$ is disjoint from $G$.
Then as in Example~~\ref{int}.\ref{contex}
we verify that $G-G'\sim_{\mathbb Q} 0$,
because $g$ is extremal and
$(G-G')|_G=G|_G\sim_{\mathbb Q} 0$.
This means that a multiple of $G$
is linearly movable.
In other cases $C+B'=F=G$ is
a double section of $g$.
As above $S/G'$ is minimal since
$S$ is not rationally connected.
Then again $G$ is linearly movable
when $S$ has a splitting type.
If $G$ is a double section,
then $C+B'=G$ and we have
(NTC) by Corollary~\ref{indcomp}.\ref{elcomp}.
Really, we have no such a case
after the covering trick.

As in \cite[Proposition~5.5]{Sh2},
(ASA) implies (EEC).

Finally, (EC)+(SM) implies (EEC)
by a monotonicity result below.
\end{proof}

\refstepcounter{subsec}
\theoremstyle{plain}
\label{monl1}
\newtheorem*{monlemma}{\thesection.\thesubsec.
Monotonicity Lemma}
\begin{monlemma}
Let $r=(n-1)/n$ for a natural number $n\not=0$.
Then, for any natural number $m\not=0$,
$$\lfloor(m+1)r\rfloor /m \ge r.$$
Moreover, for $n\ge m+1$,
$$\lfloor(m+1)r\rfloor /m = 1.$$
\end{monlemma}

\begin{proof} Indeed, since $m r$ has
denominator $n$ and $r=1-1/n$, then
$$\lfloor(m+1)r\rfloor= \lfloor m r+r\rfloor
\ge m r.$$
This implies the inequality.

Since $\lfloor(m+1)r\rfloor /m$ has
denominator $m$, and it is always $\le 1$,
we get the equation.
\end{proof}

\begin{proof}[Proof of the Inductive Theorem: Big case]
Suppose that $-(K+C+B)$ is big
(as in \cite[Theorem~5.6]{Sh2}).
(Cf. \cite[Theorem~19.6]{KC}.)

The proof is based on the Kawamata-Viehweg vanishing
and the log singularities
connectedness~\cite[Lemma~5.7 and Theorem~6.9]{Sh2}.
Kawamata states that the vanishing works even
for $\mathbb R$-divisors.
In our situation it is easy to replace $B$
by a new $\mathbb Q$-divisor $\le B$
with the same regular complements.
We can make a small decrease of each
irrational $b_i$ in prime $D_i$ with $(K+B+C.D_i)<0$.
Finally, after this procedure $B$ will be rational.
Equivalently, each $b_i$ in $D_i$
with $(K+B+C.D_i)=0$ is rational.
Since $-(K+C+B)$ is big,
each such $D_i$ is contractible and
corresponding $b_i$ is rational.
\end{proof}

A higher dimensional version can be done similarly 
(cf. Proof of Theorem~\ref{csing}.\ref{singth}).
Other cases in the Inductive Theorem
are more subtle and
need more preparation.
(A higher dimensional version of that
will be presented elsewhere.) 
But first we derive some corollaries.

\refstepcounter{subsec}
\begin{corollary}
If $(S,B)$ is
a weak log Del Pezzo with $(K+B)^2\ge 4$,
then it has a regular complement.
\end{corollary}

\begin{proof}
We need to find $B'>0$ such that $(S,B+B')$ is
a weak log Del Pezzo, but
$K+B+B'$ is not Kawamata log terminal.
By the Riemann-Roch formula and
arguments in the proof of \cite[Lemma~1.3]{Sh1},
it follows from
the inequality:
\begin{eqnarray*}
(N(-K-B).N(-K-B)-K)/2 &=N(N+1)(K+B)^2/2+(-K-B.B)/2\\
&\ge 2N(N+1)> 2N(2N+1)/2.
\end{eqnarray*} 
\end{proof}

It implies non-vanishings as in
Corollaries~\ref{int}.\ref{nonvan1}-\ref{nonvan3}, 
in particular, the following non-vanishing.

\refstepcounter{subsec}
\begin{corollary} If $(S,0)$ is
a weak Del Pezzo with log canonical
singularities and $K^2\ge 4$,
then $|-12K|\not=\emptyset$ or $h^0(X,-12K)\not=0$.
\end{corollary}

\refstepcounter{subsec}
\begin{corollary}
A weak log Del Pezzo $(S,B)$ is exceptional 
only when $(K+B)^2< 4$.
\end{corollary}

\refstepcounter{subsec}
\begin{definition}
Let $D$ be a divisor of
a complete algebraic variety $X$.
We say that $D$ has
a {\it type of numerical dimension\/} $m$
if $m$ is the maximum
of the numerical dimension for
effective $\mathbb R$-divisors $D'$ with
$\Supp{D'}\subseteq\Supp{D}$.
Similarly, we define
a {\it type of linear dimension\/} where
we replace the numerical dimension of $D'$
by the Iitaka dimension of $D'$.
Note that both $0\le m\le \dim X$.
A {\it big\/} type is a type
with $m=\dim X$.
\end{definition}

Here are the basic properties.
\begin{itemize}
\item
The linear type $m$ is
the Iitaka dimension of $D$ when $D$
is an effective $\mathbb R$-divisor.
\item
Such $D$ is movable if and only if
$D$ has the linear type $m\ge 1$.
\item
The linear type is a birational
invariant for log {\it isomorphisms}, i.~e.,
we add exceptional divisors for extractions
and contract only divisors in $D$.
In particular, this holds for
the extractions and flips.
\item
For arbitrary extraction,
the numerical type is not
decreasing.
\item
For arbitrary log transform,
the linear type is not decreasing.
\end{itemize}

For surfaces we can prove more.

\refstepcounter{subsec}
\label{types}
\begin{p-definition}
Let $X=S$ be a complete surface,
and $D$ be a divisor.
Then the numerical type of $D$ is not
higher than that of the linear type.

Moreover, if $(S,C+B)$ is a log surface
such that
\begin{itemize}
\item
$K+C+B$ log terminal,
\item
$K+C+B\equiv 0$ , and
\item
$\Supp D$ is {\it divisorially\/} disjoint
from $LCS{(S,C+B)}$, i.e, they have no
divisorial components in common.
\end{itemize}
Then the numerical type of $D$ is
the same as the linear one.
In addition to the big type we have
a {\it fibre\/} type for a type of dimension $1$,
and that of an
{\it exceptional} type for a type of dimension $0$.

More precisely, $D$ is
a fibre {\it geometrically\/},
i.~e., there exists a fibre contraction
$g:S\to Y$ with an algebraic fibre $D'=g^*P$,
for $P$ in a non-singular curve $Y$, having
$\Supp{D'}\subseteq\Supp{D}$, if and only if
$D$ has a fibre type.
In addition, $D$ is supported in fibres
of $g$ and $g$ is defined uniquely by $D$.

Respectively, $D$ is
exceptional {\it geometrically\/},
i.~e., there exists a birational contraction of $D$
to a $0$-dimensional locus, if and only if
$D$ has an exceptional type.

The condition $K+C+B\equiv 0$ can be
replaced by (ASA).
\end{p-definition}

\begin{proof}
We use the semi-ampleness
\cite[Theorem~2.7]{Sh3}
for $K+C+B+\varepsilon D'
\sim_{\mathbb R}D'$ with small $\varepsilon>0$.
Note that $K+C+B\sim_{\mathbb R} 0$.

If $D$ has a fibre type we have
a fibering $g:S\to Y$ with
a required fibre $D'=g^*P$.
If $D$ has a horizontal irreducible component $D''$,
then $D'+\varepsilon D''$, will have a big type
for small $\varepsilon >0$.
{\it Horizontal\/} means not in fibres,
{\it vertical\/} in fibres.

If $D$ is of an exceptional type,
we can contract $D$ due to Artin or by the LMMP
(cf. with the proof of
Proposition~\ref{indcomp}.\ref{relcond}).

Finally,
if $K+C+B$ satisfies (ASA), we
have a numerical complement $K+C+B'\equiv 0$
with the required properties.
\end{proof}

\begin{proof}[Proof of the Inductive Theorem:
Strategy]
In this section we assume that $Z=\pt$.
The local cases, when $\dim Z\ge 1$,
we discuss in Section~\ref{lcomp}
where we will obtain better results.

Using a log terminal blow-up
\cite[Example~1.6 and Lemma~5.4]{Sh2},
we reduce the problem to the case
when $K+C+B$ is log terminal.
By our assumption it has a non-trivial reduced
component $C\not=0$.
We would like to induce complements from
lower dimensions
(in most cases, from $C$).

First, we prove the theorem in
\begin{description}
\item{ Case~I:}
$C$ is not a chain of rational curves.
\end{description}
After that we assume that $C$ is
a chain of rational curves.

By Proposition~\ref{indcomp}.\ref{relcond},
we can assume (NTC) or, equivalently, (ASA).
So, we have a numerical contraction
$\nu : S\to Y$ for $-(K+C+B)$,
where $Y$ is a non-singular curve
or a point, the latter for a {\it while\/}.
In the exceptional cases we have regular complements.

Let $\kappa^*$ be the numerical dimension
of $-(K+C+B)$, equivalently, $\dim Y$.

We construct complements in different cases
according to configuration of $C+B$ with
respect to $\nu$.
For $\kappa^*=1$, we distinguish two cases:
\begin{description}
\item{Case~II:}
$\nu$ has a multi-section in $C$.
\item{Case~III:}
$C$ is in a fibre of $\nu$.
\end{description}

For $\kappa^*=0$, most of cases
are inserted in above ones.

A pair $(S,K+C+B)$ with $\kappa^*=0$ is
considered as Case~II, when
$D=\Supp{B}$ has a big type.
As we see in a proof below, $\nu$ will
naturally arise when we need it,
and such $\nu$ will be a contraction
to a curve.

Such a pair is considered as Case~II, when
$D$ has a fibre type and
$C+D$ is of a big type.
Equivalently, $C$ has a multi-section for
$\nu=g$ with $g$ of
Proposition~\ref{indcomp}.\ref{types}.
Moreover, then $C$ has a double section of $g$.
In this case we will have a complement
by Lemma~\ref{indcomp}.\ref{indcomp2} .

Such a pair is considered as Case~III, when
$D$ has a fibre type, but now $C+B$ and $D$
sit only in fibres of $\nu=g$.
One component of $D$ gives
a (geometric) fibre of $\nu$.
Contraction $\nu=g :S\to Y$
plays here the same role as $\nu$
in Case~III with $\kappa^*=1$.

Finally,
\begin{description}
\item{Case IV:}
$D$ has an exceptional type.
\end{description}
Then we contract the boundary $B$ to points.

Case~I could be distributed in
the other cases.
We prefer to simplify the geometry of $C$,
which 
simplifies slightly the proofs of Cases II-IV.

We try to reduce each case to the Big one
for $K+C+B$, or for $K+C+B'$
with $B'=\lfloor (n+1)B\rfloor/n$, when
a complement index $n$ is suggested.
There are two obstacles here.
First, we cannot change $B$ or $B'$ in
such a way.
For instance, in Case~III with
curves of genus $1$ in fibres.
This leads to a separation into the cases.
Second, we cannot preserve complements
of some indices.
For instance, decreasing $B$ as,
Lemma~\ref{indcomp}.\ref{monl2} shows.
Here is a situation as in
Example~\ref{indcomp}.\ref{nonrcomp},
and it is a main difficulty in Case~II,
when we try to induce a (regular) complement
from $C$.
To resolve this difficulty,
we use Lemmas \ref{indcomp}.\ref{indcomp1}
and \ref{indcomp}.\ref{indcomp2} below.
In the cases when we cannot induce regular
complements, we find others by
Lemmas~\ref{indcomp}.\ref{compl1}-\ref{compl3}.
This occurs only when (M) is not assumed.

In Case~III, the most difficult
situation, 
when $\nu$ has curves of genus $1$
in fibres,
is quite concrete.
Then we use the Kodaira's classification
of degenerate fibres.
We also use {\it indirectly\/}
a Kodaira's formula \cite[p. 161]{BPV}
for canonical divisor $K$.

A possible alternative approach to Cases II-III
is to use directly an analog of Kodaira's formula:
\begin{equation}
\label{adj}
K+C+B\sim_{\mathbb R}
\nu^*(K_Y+\nu_*((K+C+B)_{S/Y})+B_Y)
\end{equation}
for a certain boundary $B_Y$,
which can be found locally$/Y$.
We assume that $K+C+B\equiv 0/Y$.
So, this is a fibering
of genus $1$ log curves.
In arbitrary dimension $n=\dim S$,
formula (\ref{adj}), with a boundary $D$
instead of $\nu_*((K+C+B)_{S/Y})+B_Y$,
was proposed, but not proved,
in the first draft of the paper.
Formula (\ref{adj}) also plays an important role
in a proof of an adjunction
in codimension $n$.
In this context,
a similar formula was proved,
in our surface case,
by Kawamata \cite{K2}, where
$B_Y$ corresponds to a divisorial part and
$f_*((K+C+B)_{S/Y})$ to a moduli part.
However, we have three difficulties in its application.
First, a relative log canonical
divisor $(K+C+B)_{S/Y}$ is nowhere defined,
and its properties are nowhere be found.
Second, the divisorial part
is given but not very explicitly.
Third, we have $\sim_{\mathbb R}$ or
$\sim_{\mathbb Q}$ for $\mathbb Q$-boundaries.
So, we need to control indices for
a complement on $(Y,B_Y)$ and
for $\sim_{\mathbb Q}$ in (\ref{adj}) for
an induced complement.

Finally, in Case~IV we decrease $C$
(cf. proof of
Corollary~\ref{int}.\ref{mcorol})
and use a covering trick.
In this case we have trivial regular complements.

A modification of a log model may change
{\it types\/}, i.e.,
possible indices, of complements.
Nonetheless for blow-down we
can always induce a complement of
the same index from above
by \cite[Lemma~5.4]{Sh2}.
An inverse does not hold in general.
For instance, there are no complements
after many blow-ups in generic points.
Lemmas~\ref{indcomp}.\ref{mlem1} and
\ref{gcomp}.\ref{mlem3} give certain sufficient
conditions when we can induce complements from below.
\end{proof}

\refstepcounter{subsec}
\begin{lemma}
\label{mlem1} 
(Cf.~\cite[Lemma 5.4]{Sh2}.)
In the notation of \cite[Definition~5.1]{Sh2},
let $f\colon X\to Y$ be
a birational contraction such that
$K_X+S+\lfloor (n+1)D\rfloor/n$ is
numerically non-negative on a
{\rm sufficiently general\/} curve$/Y$ in each
exceptional divisor of $f$.
Then
$$K_Y+f(S+D)\ n-complementary\ \ \Longrightarrow
\ K_X+S+D\ \ n-complementary.$$
\end{lemma}

By a {\it sufficiently general}  curve in
a variety $Z$ we mean a curve which belongs
to a covering family of curves in $Z$ 
(cf. \cite[Conjecture]{Sh4}). 
Note that for such curve $C$ and
any effective $\mathbb R$-Cartier divisor $D$,
$(D.C)\ge 0$.

\refstepcounter{subsec}
\begin{example-c}
\label{mlem2}
If $S+D=S$ is in a neighborhood of exceptional
locus for $f$, and $K_X+S+D$ is nef$/Y$,
then we can pull back the complements, i.e.,
for any integer $n>0$,
$$K_Y+f(S+D)\ \ n-complementary\ \ \Longrightarrow
\ K_X+S+D\ \ n-complementary.$$
\end{example-c}

\refstepcounter{subsec}
\label{negat} 
\theoremstyle{plain}
\newtheorem*{negatl}{\thesection.\thesubsec.
Negativity of a proper modification}
\begin{negatl}
(Cf.~\cite[Negativity~1.1]{Sh2}.)
Let $f:X\to S$ be
a proper modification morphism,
and $D$ be a $\mathbb R$-Cartier divisor.
Suppose that
\begin{description}
\item{\rm (i)}
$f$ contracts all components $E_i$ of $D$ with
negative multiplicities, {\rm i.e., such components
are exceptional for $f$\/};
\item{\rm (ii)}
$D$ is numerically non-positive on a
{\rm sufficiently general\/} curve$/S$ in each
exceptional divisor $E_i$ of (i).
\end{description}
Then $D$ is effective.
\end{negatl}

\begin{proof}
First, it is a local problem$/S$.

Second, according to Hironaka,
we may assume that $f$ is projective$/S$,
and $X$ is non-singular, in particular, $\mathbb Q$-factorial.
The pull back of $D$ will satisfy the assumptions.
The proper inverse image of a sufficiently general curve
will be again sufficiently general.

Third, there exists an effective Cartier divisor
$H$ in $S$ such that
\begin{description}
\item{\rm (iii)}
the support of $f^*H$ contains the components
of $D$ with negative multiplicities; and
\item{\rm (iv)} 
$f^{-1}H$ is positive on
the sufficiently general curves in (ii)
(it is meaningful since $f^{-1}H$ is Cartier).
\end{description}

For example, we may take, as $H$, a general hyperplane
through the direct image of an effective and
relatively very ample divisor$/S$, and through
the images of the components
of $D$ with negative multiplicities.

According to (iii-iv), Cartier divisor $E=f^*H-f^{-1}H$
is effective with positive multiplicities for
the components of $D$ having negative multiplicities,
and
\begin{description}
\item{\rm (v)} 
negative on the sufficiently
general curves in (ii).
\end{description}

Thus there exists a positive real $r$ such that
$D+rE\ge 0$ and has an exceptional component $E_i$ with 0
multiplicity unless $D\ge 0$.
However then $(D+rE.C)\ge 0$ for a general curve $C$
in $E_i$, which contradicts (ii) and (v).
\end{proof}

\begin{proof}[Proof of
Lemma~\ref{indcomp}.\ref{mlem1}]
We take a crepant pull back:
$$K_X+D^{+X}=f^*(K_Y+D^+).$$
It satisfies \cite[5.1.2-3]{Sh2} as
$K_Y+D^+$, and we need to check
\cite[5.1.1]{Sh2} only for the exceptional
divisors.
For them it follows from our assumption and
the Negativity~\ref{indcomp}.\ref{negat}.
\end{proof}

\refstepcounter{subsec}
\label{monl2}
\begin{monlemma}
Let $r$ be a real
and $n$ be a natural number.
Then
$$\lfloor(n+1)(r-\varepsilon)\rfloor /n =
\lfloor(n+1)r\rfloor /n,
$$
for any small $\varepsilon>0$, if and only if
$r\not\in\mathbb Z/(n+1)$.
\end{monlemma}

\begin{statement1}
Note, that for $r=k/(n+1)>0$,
we have
$\lfloor(n+1)r\rfloor /n=k/n>r$.
So, $r>0$ is not in $Z/(n+1)$ if
$\lfloor(n+1)r\rfloor /n=k/n\le r$.
\end{statement1}

\begin{proof}[Proof of the Inductive Theorem:
Case~I]
By (ASA) and \cite[Proposition~5.5]{Sh2},
we can assume that $K+C+B\equiv 0$.
Then by \cite[Theorem~6.9]{Sh2}, or
by the LMMP and \cite[Lemma~5.7]{Sh2}, $C$
has a single connected component,
except for the case when $C$
consists of two non-singular disjoint
irreducible components
$C_1$ and $C_2$ such that $S$ has
a ruling $g:S\to C_1\cong C_2/Y$
with sections $C_1$ and $C_2$
(cf. Theorem~\ref{clasf}.\ref{clasfth} below).
By \cite[Lemma~5.7]{Sh2},
an existence of the ruling follows
from the LMMP applied to $K+C_2+B$.
A terminal model cannot be a log Del Pezzo surface,
because $C_i$ will always be disjoint and
non-exceptional during the LMMP
(cf. Theorem~\ref{clasf}.\ref{dipole}).

In this case we can construct complements
to $K+C+B$ with given $C=C_1+C_2$.
(Cf. Lemma~\ref{indcomp}.\ref{indcomp1}
and~\ref{indcomp}.\ref{indcomp2} below.)
Making a blow-up of the generic point of $C_1$ and
contracting then the complement component of the fibre
we reduce to the case when $C_2$ is big.
Then, for small rational $\ep>0$,
$-(K+C+B-\ep C_2)=-(K+C+B)+\ep C_2$ is nef and big
which gives a required complement
by the Big case above.
More precisely, we have the same complement
as $K_{C_1}+(B)_{C_1}$ on $C_1$, where
$(B)_{C_1}$ denotes different
\cite[Adjunction~3.1]{Sh2}.
The blow-ups preserve complements by
\cite[Lemma~5.4]{Sh2}.
For contractions we may use
Example~\ref{indcomp}.\ref{mlem2}.
In this case it is easily verified directly
as well.

So, we suppose that $C$ is connected.
We also assume that each component
of $C$ is a non-singular rational curve.

Otherwise $C$ is
a (non-singular, if we want) irreducible curve
of genus $1$ with
$B=0$ and non-singular $S$ in a neighborhood of $C$
[Sh2, Properties 3.2 and Proposition 3.9].
In this case we suppose that $S$ is non-singular
everywhere after a minimal resolution.
We chose $B$ which is crepant for the resolution.
Then the LMMP and a classification of
contractions in a 2-dimensional minimal model program
gives a ruling $f\colon S\to C'$
with a surjection $C\to C'$, or $S=\mathbb P^2$
with a cubic $C+B=C$.
Since a blow-up in any point of $\mathbb P^2$ gives
a ruling, then $\mathbb P^2$ is
the only possible case,
when $S=\mathbb P^2$, and
$C=C+B\sim -K_{\mathbb P^2}$.
Hence $K+C+B=K+C\sim 0$ on $S=\mathbb P^2$, and
we have a $1$-complement.
Similarly, a $1$ or $2$-complement holds
in the case with the ruling $f$
and double covering $C\to C'$
(cf. Lemma~\ref{indcomp}.\ref{indcomp2}).
We want to check only that $B=0$
and $2(K+C+B)=2(K+C)\sim 0$
but not only $\equiv 0$.
After contractions of exceptional curves of the first kind
intersecting $C$ in fibres of $f$ we may suppose that
$f$ is extremal, i.e., with irreducible fibres.
This preserves all types of complements
by Example~\ref{indcomp}.\ref{mlem2}.
Note also that the same reduction holds for
the ruling $f\colon S\to C'=C$ having a section $C$.

Since $K+C+B\equiv 0$ and covering $C\to C'$ is double,
we have no boundary components in fibres of $f$ and $B=0$.
So, if $C$ is rational, then $f$ is
a rational ruling $\mathbb F_n$, and $K+C+B=K+C\sim 0$.
Otherwise $C'$ has genus $1$, and
$(K+C)$ is $1$- or $2$-complementary
by Corollary~\ref{indcomp}.\ref{elcomp}.
If $S$ is a surface of
Example~\ref{indcomp}.\ref{noncex},
then $C=C_i$ and we have a $2$-complement.
The surface $S$ of Example~\ref{int}.\ref{contex}
is impossible by (ASA).
If $f$ is a splitting case:
$-K_S\sim G+G'$.
Therefore $C\equiv G+G'$, as well as $C\sim 2G$ and $2G'$,
because $C\cap G=C\cap G'=G\cap G'=\emptyset$.
Thus
$2C\sim 2(G+G')\sim -2K_S$ which gives a $2$-complement
(but not $1$-complement).

Next we consider a case when $C=C'$ is
a non-singular curve of genus $1$
and a section of $f$.
Boundary $B=\sum b_i D_i\not=0$ has only horizontal
components $D_i$.
The curves $D_i$ are non-rational.
Hence $D_i^2\ge 0$ by the LMMP.
On the other hand, we have a trivial $n$-complement
for $(K+C+B)|_C=K_{C'}\sim 0$ for any natural number $n$.
Thus we have an $n$-complement when
we have $D_i^2>0$ with multiplicity
$b_i\not\in\mathbb Z/(n+1)$.
Indeed, we can then construct an $n$-complement
as in the Big case, for $(K+C+B-\varepsilon D_i)$
with small $\varepsilon>0$.
For $(K+B+C)$, we have the same complement
by the Monotonicity Lemma~\ref{indcomp}.\ref{monl2}.
Since $\mathbb Z/2\cap\mathbb Z/3=\emptyset$
in the unit interval $(0,1)$,
then $K+C+B$ is $1$- or $2$-complementary,
whenever some $D_i^2>0$.
Otherwise all $D_i^2=0$.
By Corollary~\ref{indcomp}.\ref{elcomp},
the curves $D_i$ and $C'$ are in $R_2$,
and are all disjoint non-singular curves
(of genus $1$).
Hence it is enough to find an $n$-complement
in a generic fibre$/Z$, which
we have for $n=1$ or $2$ by
\cite[Example~5.2.1]{Sh2}.

Finally, we suppose that $C$ is a (connected)
wheel of rational curve.
Then by the arguments of the case when $C$ is
a non-singular curve of genus $1$,
we see that $S$ is
a rational ruling $S\to C'$ with
a double covering $C\to C'$,
or $S=\mathbb P^2$ with a cubic $C+B=C$.
In both cases we have a $1$-complement.
\end{proof}

\refstepcounter{subsec}
\label{monl3}
\begin{monlemma}
Let $r<1$ be a rational number
with a positive integer denominator $n$,
and $m$ be a positive integer, such that
$n|m$.
Then
$$\lfloor(m+1)r\rfloor /m \le r.$$
Moreover, this gives $=$
if and only if $r\ge 0$.
\end{monlemma}

\begin{proof}
Let $r=k/n$.
Then
$$\lfloor(m+1)r \rfloor /m=
\lfloor(m+1)k/n\rfloor/m=
(k m/n+\lfloor k/n\rfloor)/m=
r+\lfloor r \rfloor/m\le r,$$
and $=r$ if and only if $r\ge 0$.
\end{proof}

\refstepcounter{subsec}
\begin{corollary} Let $m$ be a natural number,
and $D$ be a subboundary of index $m$ in codimension 1
and without reduced part,
i.e., $m D$ is integral with the multiplicities $<m$.
Then
$$m D\ge \lfloor (m+1)D\rfloor,$$
and $=$ if and only if $D$ is a boundary.
\end{corollary}

\refstepcounter{subsec}
\label{ncomp} 
\begin{lemma}
Let $C=\lfloor C+B\rfloor$ be
the reduced component
in a boundary $C+B$ on a surface $S$, and
$C'\subseteq  C$ be
a complete curve such that
\begin{description}
\item{\rm (i)} 
$K+C+B$ is (formally) log terminal
in a neighborhood of $C'$;
\item{\rm (ii)} 
$(K+C+B)|_{C'}$ has an $n$-complement; and
\item{\rm (iii)} 
$-(K+C+B)$ is nef on $C'$.
\end{description}
Then $(C_i.K+C+D)\le 0$ on each
component $C_i\subseteq C'$
with $D=\lfloor (n+1)B\rfloor/n$; and
$K+C+D$ is log canonical
in a neighborhood of $C'$.
\end{lemma}

In (i), {\it formally\/} means locally
in an analytic or etale topology.
This can be defined formally as well.

\begin{proof} First, we may suppose that
$C'$ is connected.

Second, by [Sh2, the proof of Theorem 5.6]
(cf. Proof of the Big case in the Inductive Theorem
and that of the Local case in Section~\ref{lcomp}),
the lemma holds when $C'$ is contractible
because then we have an $n$-complement
in a neighborhood of $C'$.

Third, it is enough for an analytic contraction,
because  under our assumptions in most cases,
it will be algebraic
due to Artin or by the LMMP.
It works and we get a rational
singularity after the contraction
when $C'$ is not isolated in $C+B$.
Otherwise $B=0$ and $C+B=C+D=C'$.
Then the lemma follows from (i) and (iii).

Finally, the contraction exists when
a certain numerical condition on the intersection form
on $C'$ is satisfied.
This will be negative after making
sufficiently many monoidal transforms
in generic points of $C'$.
Such a crepant pull back of $K+C+D$ preserves
the assumptions (i)-(iii) and the statements.
\end{proof}

\refstepcounter{subsec}
\label{indcomp1} 
\begin{lemma}
Let $(S,C+B)$ be a complete log
surface with a ruling $f\colon S\to Z$
such that
\begin{description}
\item{\rm (i)}
there exists
a section $C_1\hookrightarrow S$
of $f$ which is in the reduced part $C$;
\item{\rm (ii)} 
$(K+C+B)|_{C}$ has an $n$-complement
for some natural $n>0$;
\item{\rm (iii)} 
$C+D=C+\lfloor (n+1)B)\rfloor/n$
gives an $n$-complement near the generic
fibre of $f$, {\rm i.e., $K+C+D$
is numerically trivial on it};
\item{\rm (iv)}
$-(K+C+B)$ is nef; and
\item{\rm (v)} 
$K+C+B$ is (formally) log terminal in
a neighborhood of $C$, if we do not assume
that $K+C+B$ is log canonical everywhere but
just $C+B\ge 0$ outside $C$.
\end{description}
Then $K+C+B$ has an $n$-complement.
\end{lemma}

\begin{proof} The above numerical property (iv) and
\cite[Theorem 6.9]{Sh2} imply that
$K+C+B$ is log canonical everywhere
(cf.~\cite[the proof of Theorem 5.6]{Sh2}).

Making a crepant log blow-up we may
assume that $K+C+B$ is log terminal everywhere,
essentially by \cite[Lemma 5.4]{Sh2}.
Since $f$ is a ruling,
then $\NEc{S/Z}$ is rational polyhedral and generated
by curves in fibres of $f$
(cf.~(EX1) in \ref{indcomp}.\ref{relcond}).
Note that any contraction of
a curve $E\not\subseteq C$ in fibres of $f$
will preserve (i)-(v): (ii) by \cite[Lemma~5.3]{Sh2}
because the boundary coefficients
of $(K+C+B)|_{C}$ are not increasing.
This implies (v) as well.

We consider simultaneously
the boundary $C+D$ as in (iii).
After contractions of curve $E\not\subseteq C$
in fibres of $f$ with $(E.K+C+D)\ge 0$,
we may suppose that $-(K+C+D)$ is nef.
Indeed, it is true for the fibres
and on section $C_1$.
Since $K+C+D\equiv 0/Z$, applying the LMMP to $f$
we may suppose that $f$ is extremal.
Then $\NEc{S}$ is generated by a fibre and a section.
Note that $(C_i.K+C+D)\le 0$
for the curves $C_i\subseteq C$
by Lemma~\ref{indcomp}.\ref{ncomp}.

It is enough to construct an $n$-complement
after such contractions by
Lemma~\ref{indcomp}.\ref{mlem1}.

The boundary coefficients of $D$ belong to $\mathbb Z/n$.
In addition, by (iii)
\begin{description}
\item{\rm (vi)}
$K+C+D$ is numerically trivial$/Z$.
\end{description}
Therefore again by
Lemma~\ref{indcomp}.\ref{mlem1},
we may assume that the fibres of $f$
are irreducible or in $C$.
Since $C_1$ is a section,
we increase $C+D$ in fibres to $B^+$
in such a way that $(K+C+B^+)|_{C_1}$ is given by an
$n$-complement in (ii).
We contend that $K+B^+$ gives an
$n$-complement of $K+C+B$, too.

First, note that \cite[5.1.1]{Sh2}
holds by the construction.

Second, as above, $K+B^+\equiv 0$ because
it is true for the fibres and section $C_1$.

Third, $K+B^+$ is log canonical in a neighborhood of $C_1$
by the Inverse Adjunction \cite[Corollary~9.5]{Sh2}.
So, as above, $K+B^+$ is log canonical
everywhere, i.e., \cite[5.1.2]{Sh2}.

Finally, we need to check that $n(K+B^+)\sim 0$.
In particular it means that $n B^+$ is integral.
Since the log terminal singularities are rational
as well as any contractions of curves in fibres of $f$,
we can replace $(S/Z,C+B^+)$ by
any other crepant birational model.
For instance, we can suppose that $S$ is non-singular,
and all fibres of $f$ are irreducible.
Then $S$ is a non-singular minimal ruling$/Z$
with section $C_1$.
In that case $n(K+B^+)$ is integral and $\sim 0$
by the Contraction Theorem because these hold
for $n(K+B^+)|_{C_1}$.
The latter will be preserved
after any above crepant modification.
\end{proof}

\refstepcounter{subsec}
\label{indcomp2} 
\begin{lemma}
Lemma~\ref{indcomp}.\ref{indcomp1} holds even
if we drop (iii), but at the same time,
change (i) in it by
\begin{itemize}
\item 
there exists a curve $C_1$ in
$C$ with a covering $C_1\to Z$ of degree $d\ge 2$,
\end{itemize}
except for the case $C+B=C_i$ in
Example~\ref{indcomp}.\ref{noncex},
when $n$ is odd.
Moreover, then $d=2$ always.
In the exceptional case we have
a $2n$-complement for any natural number $n$.
\end{lemma}

\refstepcounter{subsec}
\label{invinv}
\begin{lemma}
Let $f\colon X\to Y$
be a conic bundle contraction with
a double section $C$.
If a divisor $D\equiv 0$
over generic points of codimension 2 in $Y$,
and, for any component $C_i$ of $C$,
$C_i\not\subseteq D$,
then the different $D_{C^{\nu}}$ is invariant under
the involution $I$ given by the double covering
$C^{\nu}\to Y$ on the normalization $C^{\nu}$.
\end{lemma}

\begin{statement1}
The same holds for $C$ assuming that
$K+C$ is log canonical in codimension 2
(cf. {\rm \cite[Theorem~12.3.4]{KC}\/}).
\end{statement1}

\begin{proof}[Proof-Commentary]
First, taking hyperplane sections, we
reduce the lemma to the case of a surface
ruling $X\to Y$ with a double curve $C^{\nu}$
over $Y$.

Second, we can drop $D$, because it is
pull backed from $Y$.

Finally, according to the numerical definition
of the different \cite{Sh2}, and because
$K+C\equiv 0/Y$ we may replace $X$ by any
crepant model $(X,D)$.
In particular, we may suppose that $X$ is
non-singular with an extremal ruling $f$.
According to M. Noether, we may assume that
$C$ is non-singular as well.
Then $D\equiv 0/Y$, because it is
supported in fibres of $f$.
So, we may drop $D$ again.

In this case the different is $0$ and invariant.

The same works for~\ref{indcomp}.\ref{invinv}.1.
We need the log canonical condition on $K+C+D$
only to define $(K+C+D)|_{C}$.
\end{proof}

\refstepcounter{subsec}
\label{invcomp}
\begin{lemma}
Let $C_1$ be a component of
a semi-normal curve $C$ with a finite Galois covering
$f\colon C_1\to C'$ of a {\rm main type\/} $\mathbb A$,
and let $B$ be a Weil $\mathbb R$-boundary supporting
in the normal part of $C$ and  Galois invariant on $C_1$.
Then $K+B$ has an $n$-complement which is Galois invariant
on $C_1$, if and only if it has an $n$-complement.
\end{lemma}

\begin{proof}[Proof-Commentary]
The type $\mathbb A$ means that we have
branchings at most over two points
$Q_1$ and $Q_2\in C'$ in each irreducible
component of $C'$.

According to \cite[Example~5.2.2]{Sh2},
$K+B$ has an $n$-complement if and only if
$$K+\lfloor B\rfloor +\lfloor (n+1)\{B\}\rfloor/n$$
is non-positive on all components of $C$.
This is a numerical condition
which can be preserved, if we first
replace $C$ by $C_1$, and even by any irreducible
component with the Galois covering $f:C=C_1\to C'$
given by the stabilizer of this component.
Intersections with other components we
include into the boundary with multiplicity $1$.
We also assume that $f$ is not an isomorphism.

If $C$ is singular, $B^+=0$ is invariant.
If $C$ is non-singular, by
the Monotonicity Lemma~\ref{indcomp}.\ref{monl3},
we suppose that $n B$ is integral.
If $K+B\equiv 0$ then we have a required
complement by the above criterion, and
it is invariant by our conditions.

Otherwise $\deg (K+B)<0$, and
$C$ is rational.
Then under our conditions on the branchings we
have them only over two points $Q_1$ and $Q_2\in C'$.
In addition, we have
unique points $P_1/Q_1$ and $P_2/Q_2$
respectively with maximal ramification
indices $\deg f-1$. 
Indeed, 
$$
-2=\deg K_C= (\deg f)(K_{C'}+(\frac{r_1-1}{r_1})Q_1+
(\frac{r_2-1}{r_2})Q_2)=-(\deg f)(\frac1{r_1}+
\frac1{r_2}),
$$
where $r_i|\deg f$ is the ramification multiplicity
in $P_i/Q_i$.
So, we can maximally extend $B$ in $P_1$ and $P_2$
preserving the following properties
\begin{itemize}
\item
$B$ is Galois invariant;
\item
$n B$ is integral; and
\item
$\deg (K+B)\le 0$.
\end{itemize}
Then $\deg (K+B)=0$, because $\deg K=-2$.
Again by the numerical condition we are done.
\end{proof}

\refstepcounter{subsec}
\begin{example}
For other types of
Galois action we may lose
Lemma~\ref{indcomp}.\ref{invcomp}.

Let us consider, for example, a type $\mathbb D$.
In this case we have a Galois covering
$f \colon C\to C'$ such that
\begin{itemize}
\item
the curves $C$ and $C'$ are
isomorphic to $\mathbb P^1$;
\item
$f$ is branching over
three points $Q_1,Q_2$ and $Q_3$;
\item
$f$ has $2$ branching
points $P_{1,1}$ and $P_{1,2}/Q_1$
with multiplicities $d$ where
$\deg f=2d\ge 4$; and
\item
$f$ has $d$ simple branching
points $P_{i,1},...,P_{i,d}/Q_i$ with
$i=2$ and $3$.
\end{itemize}
So, if $d=2m+1$ is odd and $n=d$, then for
$$B=\frac m d (P_{1,1}+P_{1,2})+
\frac 1 d (\sum P_{2,i}),$$
$\deg (K+B)=-2+2m/d+d/d=-1/d$,
whereas $B$ is invariant.
However, any $d$-complement will be
$B+(1/d)P$ which will not be
invariant for any choice of point $P$.
\end{example}

\begin{proof}[Proof of
Lemma~\ref{indcomp}.\ref{indcomp2}]
(iv) in~\ref{indcomp}.\ref{indcomp1} implies that $d=2$.
So, $D=0$ near the generic fibre of $f$,
and (iii) in~\ref{indcomp}.\ref{indcomp1}
is satisfied.
Hence by Lemma~\ref{indcomp}.\ref{indcomp1}
we may assume that $C_1$ is irreducible.

By (v) in~\ref{indcomp}.\ref{indcomp1} and
\cite[Lemma~3.6]{Sh2}, $C_1=C_1^{\nu}$ is
{\it non-singular\/}, except for the case
when $C_1$ is a Cartesian leaf, i.~e.,
an irreducible curve of arithmetic genus $1$
with a single nodal singularity.
In such a case $C+B=C+D=C_1$
and $K+C+B\equiv 0$ by the LMMP.
This gives an $n$-complement for any $n$
because $S$ is rational.

The double covering $C_1\to Z$
is given by an involution $I\colon C_1\to C_1$.
By (iii) in \ref{indcomp}.\ref{indcomp1},
$K+C+B\equiv 0/Z$.
Thus any contraction in fibres$/Z$ is crepant.
They preserve the different $(C-C_1+B)_{C_1}$.
Thus we may suppose that $S/Z$ is extremal.
Then $D=C-C_1+B=(K+C+B)-(K+C_1)\equiv
-(K+C_1)\equiv 0/Z$ and,
according to Lemma~\ref{indcomp}.\ref{invinv},
$(C-C_1+B)_{C_1}$ is invariant under $I$.
So, it has an invariant $n$-complement
by Lemma~\ref{indcomp}.\ref{invcomp},
when $C_1$ is rational.

Otherwise $C_1=C+B$ is a non-singular curve
of genus $1$ and we have
again an invariant $n$-complement
and for any $n$: $0$.

Then we may construct $B^+$ as in the proof
of Lemma~\ref{indcomp}.\ref{indcomp1},
and, after a reduction to extremal $f$,
check that $B^+$ gives an $n$-complement, except
for the case when $S$ is non-rational
and, by Corollary~\ref{indcomp}.\ref{elcomp},
$n$ is odd.
In the latter case, by the same corollary
we are in the situation of
Example~\ref{indcomp}.\ref{noncex}.
Indeed, in a splitting case,
$n(K+C+D)=n(K+C+B)=
n(K+C_1)=n(K_S+G+G')\sim 0$ for any $n$.
On the other hand, in the exceptional case,
$2n(K+C+D)=2n(K+C+B)=
2n(K+C_1)=2n(K_S+C_i)\sim 0$ for any $n$.
\end{proof}

\refstepcounter{subsec}
\label{invadj}
\begin{lemma}
Suppose that in
a surface neighborhood $S$ of a point
$P$, we have a boundary $B=C+\sum b_iD_i$
with distinct prime divisors $D_i$ such that
\begin{description}
\item{\rm (i)} 
$C$ is reduced and irreducible curve
through $P$;
\item{\rm (ii)} 
each $b_i=(m-1)/m$
for some integer $m>0$; and
\item{\rm (iii)} 
$K+B$ is log terminal.
\end{description}
Then $$(K+B)|_{C}=\frac{n-1}n P$$
with some integer $n>0$, and $m\mid n$.
\end{lemma}

\begin{statement1}
Formally at most one component $D_i$, with $b_i>0$,
is passing through $P$.
If we replace (iii) by
\begin{description}
\item{\rm (iv)}
$K+B$ is log canonical, but
not formally log terminal in $P$,
\end{description}
then $$(K+B)|_{C}=P\ ,$$
and formally
at most two components $D_i$, with $b_i>0$,
are passing through $P$.
Moreover, both have
multiplicities $b_i=1/2$, whenever
we have two of them.
\end{statement1}

\begin{proof}
By \cite[Theorem~5.6]{Sh2} and (i)-(iii),
$K+B$ has a $1$-complement $K+C+\sum D_i$
in a neighborhood of $P$.
Thus $K+C+\sum D_i$ is log canonical there.

So, by \cite[Corollary~3.10]{Sh2}
in a neighborhood of $P$,
a single divisor $D_i$ passes through $P$, and
$D_i|_{C}=1/l$ where $l$ is the index of $P$.
Therefore
$$(K+B)|_{C}=(\frac{l-1}l+\frac{m-1}{ml})P=
\frac{n-1}n P$$
with $n=l m$
(cf. \cite[the proof of Lemma~4.2]{Sh2}).

\ref{indcomp}.\ref{invadj}.1
follows from the above calculations
and \cite[Corollary~9.5]{Sh2}.
\end{proof}

\refstepcounter{subsec}
\label{adjunct}
\begin{corollary}
If a pair $(X,B)$ is log canonical and $B$
satisfies (M) or (SM), then for any
reduced divisor $C$ in the reduced part
$\lfloor B\rfloor$, the different
$(B-C)_C$ satisfies respectively (M) or (SM).
\end{corollary}

\begin{proof}
Hyperplane sections reduce the proof to
a surface case $X=S$.
It is enough to consider locally a
log terminal case.
For (M), note that $(l7-1)/l 7\ge 6/7$
for any natural number $l=n/7$.
\end{proof}

\refstepcounter{subsec}
\label{compl1}
\begin{lemma}
Let $b_i\in (0,1)$
and $n$ be a natural number.
Then pairs (i)-(ii) and (iii)-(iv) of
conditions below are equivalent:
\begin{description}
\item{\rm (i)}
$\sum b_i\le 1$,
\item{\rm (ii)}
$\sum \lfloor (n+1)b_i\rfloor/n>1$,
\item{\rm (iii)}
$b_i\in \mathbb Z/(n+1)$, and
\item{\rm(iv)}
$\sum b_i=1$.
\end{description}
\end{lemma}

\begin{proof}
The inequality (ii) holds if and only if,
for natural numbers $k_i$,
$b_i\ge k_i/(n+1)$ and $\sum k_i\ge n+1$.
Hence by (i) we have (iii)-(iv):
$1\ge \sum b_i\ge \sum k_i/(n+1)\ge 1$.
The converse follows from the same computation.
\end{proof}

\refstepcounter{subsec}
\label{compl2}
\begin{corollary}
In the notation of \cite[Example~5.2.2]{Sh2},
let $X=C$ be a chain of rational curves
with boundary $B$.
Then $(X,B)$ is
$1$-complementary, but not
$2$-complementary only when
\begin{itemize}
\item
all $b_i'\in \mathbb Z/3$ and $\deg B'=1$, or
\item
all $b_i''\in \mathbb Z/3$ and $\deg B''=1$.
\end{itemize}
In either of these cases we have
$4$- and $6$-complements.
\end{corollary}

\begin{proof} As in \cite[Example~5.2.2]{Sh2}.
\end{proof}

\refstepcounter{subsec}
\label{compl3}
\begin{lemma}
In the notation of \cite[Example~5.2.1]{Sh2},
let $X=C$ be an irreducible rational curve
with boundary $B$, and $n\in N_1$
be the minimal complementary index for $(X,B)$.
Then $(X,B)$ is
$(n+1)m$-complementary for a bounded $m$.
More precisely, if $\deg B<2$:
\begin{description}
\item{\rm for $n=1$:}
$m\in\{1,2,3,4,5,6,7,8,9,11\}$;
\item{\rm for $n=2$:}
$m\in\{1,2,3,4,5,6,7,8,10\}$;
\item{\rm for $n=3$:}
$m\in\{1,3,4,5,6\}$;
\item{\rm for $n=4$:}
$m\in\{2,3,4,5,6,8\}$; and
\item{\rm for $n=6$:}
$m\in\{3,4,5,6,8\}$.
\end{description}
\end{lemma}

\begin{proof}[Proof-Remark]
As in \cite[Example~5.2.2]{Sh2}.
But it is better to use a computer.
There exists a program by Anton~Shokurov
for $\deg B<2$, which can be easily
modified for $\deg B=2$.

On the other hand we can decrease $B$.
Then by Lemma~\ref{indcomp}.\ref{monl2},
we have the same complementary indices $(n+1)m$
as in the case with $\deg B<2$,
except for the case when $\deg B=2$
and $B$ has index $(n+1)m+1$.
Then we have an $(n+1)m+1$-complement.

From a theoretical point of view,
under the assumption that the number of
the elements in $\Supp{B}$ is bounded, we
can prove this result for any $n$,
if we check that, for any $B$, there
exists an $(n+1)m$-complement.
Of course, here is the difficult case
for $\deg B=2$ and $B$ having
irrational multiplicities, which
has been done in Example~\ref{int}.\ref{simul}.

Finally, we anticipate that the lemma
holds for arbitrary $n$ and
without the assumption that $n$
is the minimal complementary index.
For the letter, if $\deg B<2$,
a computer check shows that
\begin{description}
\item{\rm for $n=2$:}
$m\in\{1,2,3,\dots,15,16,18\}$;
\item{\rm for $n=3$:}
$m\in\{1,2,3,\dots,24,25,27\}.$
\end{description}
\end{proof}

\begin{proof}[Addition in
a proof of~\ref{indcomp}.\ref{indth}.1]
To obtain the indices in~\ref{indcomp}.\ref{indth}.1,
we unify $(n+1)m$ and $(n+1)m+1$
for $n$ and $m$ in
Lemma~\ref{indcomp}.\ref{compl3}.
\end{proof}

\begin{proof}[Addition in
a proof of~\ref{indcomp}.\ref{indth}.2]
To obtain the indices in~\ref{indcomp}.\ref{indth}.2,
we add a $6$-complement of
Lemma~\ref{indcomp}.\ref{compl2},
which follows from the proof below.
\end{proof}

\begin{proof}[Proof of the Inductive Theorem:
Case~II]
Here, we assume that $C$ has a multi-section $C_1$
of $\nu$.
If $ C_1$ is not a section we have
a regular complement from $C_1$
by Lemma~\ref{indcomp}.\ref{indcomp2} and
\cite[Examples~5.1.1-2]{Sh2}.
In the exceptional cases,
we take $n=1$ or $2$, which gives
regular complements again.

Hence we may assume that $C$ has a single
section $C_1$.
Note that $C$ is then connected by
\cite[Theorem 6.9]{Sh2}, and by our assumptions
$C$ is a chain of rational curves.
We also assume that $B$ has a big type
when $\kappa^*=0$.

Note also that $B\not=0$:
$B$ has horizontal components
for any contraction on a curve,
when $\kappa^*=0$.
Otherwise $C$ is a double section.

Let $B_C=(B)_C$ be the different
for adjunction $(K+C+B)|_{C}=K_{C}+B_C$.
Then according to our construction
the numerical dimension of $-(K+C+B)$ is
equal to that of $-(K_{C}+B_C)$, i.~e.
equal to $\kappa^*$, at least on $C_1$.

Suppose that $K_C+B_C$ is $n$-complementary
and $\kappa^*=1$.
Let $b_i$ be a multiplicity of $B$ in
a horizontal component $D_i$ such that
$b_i\not\in \mathbb Z/(n+1)$.
Then $-(K+C+B)+\varepsilon D_i=
-(K+C+B-\varepsilon D_i)$ is big
and with the same $n$-complements by
the Monotonicity Lemma~\ref{indcomp}.\ref{monl2}.
Therefore $K+C+B$ is $n$-complementary
unless all horizontal $b_i\in\mathbb Z/(n+1)\cap (0,1)$.
In this case we find a complement for
another index $m\in (n+1)\mathbb N$
by Lemmas~\ref{indcomp}.\ref{monl3} and
\ref{indcomp}.\ref{indcomp1}.
We have an $m$-complement for such $m$
by Corollary~\ref{indcomp}.\ref{compl2} and
Lemma~\ref{indcomp}.\ref{compl3}.
If $C$ is not irreducible we have
just regular complements.
Example~\ref{indcomp}.\ref{nonrcomp}
shows that we need
non-regular complements as well.

In addition,
for multiplicities under (M) or (SM),
we have only regular complements.
Indeed, then $b_i=1/2$,
because $K+C+B\equiv 0/Z$.
By \cite[Example~5.2.2]{Sh2},
we have for $K+C+B$
a (formally) regular complement,
except for the case with $n=1$,
when $K_C+B_C$ is $1$-complementary,
but not $n$-complementary for
all other regular indices $n$.
Note that $B_C$ also satisfies
(M) and (SM) by
Corollary~\ref{indcomp}.\ref{adjunct}.
Thus it is possible only when $C$
is irreducible (otherwise we have a $2$-complement)
and $B_C$ is reduced.
Additionally, by \cite[Example~5.2.1]{Sh2}
and in its notation,
$b_1\ge b_2\ge 1/2$ and all other $b_i=0$.
Hence we have a $2$-complement,
which concludes the case under (M) and (SM).

The same holds for an appropriate $\nu$
when $\kappa^*=0$.
Indeed, we have (RPC) by
Proposition~\ref{indcomp}.\ref{relcond},
because $(K+B+C)-\varepsilon D'$
satisfies (WLF) for some $D'$ with
$\Supp{D'}=\Supp{B}$ and
small $\varepsilon>0$.
Then we use arguments
in the proof of
Lemma~\ref{indcomp}.\ref{indcomp1}.
We mean that we contract exceptional curves
$E$ with $(K+C+B'.E)>0$, where
$B'=\lfloor (n+1)B\rfloor/n$.
This preserves the situation and
$n$-complements by the same reasons.
(RPC) is also preserved.
For a terminal model,
we have either a fibre contraction
such that $(K+C+B'.E)>0$
for its generic fibre, or
$-(K+C+B')$ is nef on the model.
In the former case, the contraction induces
a contraction $\nu$ as for $\kappa^*=1$.
Indeed, $C$ is a section of $\nu$.
Otherwise $C$ is in a fibre of $\nu$.
After the above contractions,
$(K+C+B'.C)>0$,
which contradicts 
Lemma~\ref{indcomp}.\ref{ncomp}.

In the other cases we assume
that $-(K+C+B')$
is nef after such contractions.
However, $K+C+B'$ may not be
log terminal, but just
log canonical and only near $C$.
So, if we would like to use
Lemma~\ref{indcomp}.\ref{indcomp1}
and the Big case, we need the following
preparation.
We contract all connected components
of the exceptional type in $B'$.
After that, by
Proposition~\ref{indcomp}.\ref{types},
we have a semi-ample divisor $D'$ with
$\Supp D'=\Supp B'$.
We contend that one can apply
Lemmas~\ref{indcomp}.\ref{indcomp1}-\ref{indcomp2}
and the Big case in that situation
with $C$ instead of $\lfloor C+B'\rfloor$.
In the lemmas we suppose
that $C$ has a multi-section for
a given contraction $f:S\to Z$ and
$K+C+B'\equiv 0/Z$.
To verify the lemmas and the Big case,
we replace $K+C+B'$ by
$K+C+B'-\varepsilon D'$.
The condition (ii) in \ref{indcomp}.\ref{indcomp1}
follows from that for $K+C+B$
and $K+C+B'$ by \cite[Lemma~5.3]{Sh2}.
Since $K+C+B'\equiv 0/Z$, we have
(iii) in \ref{indcomp}.\ref{indcomp1}
by Lemma~\ref{indcomp}.\ref{monl3}.
Other conditions (iv) and (v)
follow from the construction.
For instance, (v) holds because
otherwise we have a log canonical,
but not log terminal, point
$P\not\in \Supp{B'}$.
Then $P$ is not log terminal for
$K+C$ and $K+C+B$.
But that is
impossible by the LMMP for
$K+C+B$ versus $K+C+B'$.
Finally, $K+C+B'$ and
$K+C+B'-\varepsilon D'$ have
the same $n$-complements for
small $\varepsilon>0$ by
the Monotonicity~\ref{indcomp}.\ref{monl2}.

We continue the case with nef $-(K+C+B')$.
If $-(K+C+B')$ is big, we have
an $n$-complement as in the Big case,
as well for $K+C+B'$ and $K+C+B$ 
by our construction.
If $-(K+C+B')$ has the numerical dimension $1$,
we have a numerical contraction $f=\nu :S\to Z$
and $K+C+B'$ is $n$-complementary as
in the above case with $\kappa^*=1$, whenever
$f$ has a multi-section in $C$.
The same holds and by the same reasons
(cf. a proof in Case~III below)
if $B'$ has a horizontal component.
If $B'$ is in fibres of $f$,
then $B$ (an image of $B$)
is also in fibres of $f$,
because $K+C+B\equiv K+C+B'\equiv
K+C\equiv 0/Z$.
But $B$ has a big type.
Thus $B'$ has a horizontal component.

Finally, $K+C+B'\equiv 0$.
By the Monotonicity~\ref{indcomp}.\ref{monl2}
we will have an $n$-complement,
when $B'$ is of a big type.
As above,
if $B'$ has a fibre type
and in fibres,
this is only possible when $C$
is a multi-section for $g$ given by $B'$.
Then $K+C+B'$ and $K+C+B$ are
$n$-complementary by
Lemmas~\ref{indcomp}.\ref{indcomp1}
and~\ref{indcomp}.\ref{indcomp2}.
Since $B$ has a big type,
the case, when $B'=0$ and $K+C=K+C+B'\equiv 0$,
is impossible.
\end{proof}

\begin{proof}[Proof of the Inductive Theorem:
Case~III]
If $\kappa^*=1$ and we have
a horizontal element $D_i$ in $B$,
we can reduce this case to the Big one
for $K+C+B-\varepsilon D_i$ whenever
$D_i$ is chosen properly.
As in the proof of
Lemma~\ref{indcomp}.\ref{ncomp},
we can find
an $n$-complement, with $n\in RN_2$,
near $C$
(cf. Proof of the Local case in the Inductive Theorem
in Section~\ref{lcomp}).
In particular, by~\ref{indcomp}.\ref{monl2}.1,
we have horizontal
$D_i$ with the multiplicity $b_i$ of $B$
such that $b_i\not\in \mathbb Z/(n+1)$
(cf. Lemma~\ref{indcomp}.\ref{compl1}).
Then we can use the above reduction
by Lemma~\ref{indcomp}.\ref{monl2}.

If $\kappa^*=0$ we have $\nu$ given
by $B$ and $B$ have only vertical components.

Thus below we assume that
$B$ has only vertical components with
respect to $\nu:S\to Y$, and
$\nu$ is a fibering with curve fibres of
genus $1$.
$C$ is also vertical and,
according to Kodaira,
a modification of its fibre has
type $\rm I_b^*$, II, II$^*$, III, III$^*$, 
or IV, IV$^*$.
Near such fibres we have
respectively $n=2$-, $6$-, $4$- or
$3$-complements
(cf. Classification~\ref{lcomp}.\ref{class}).
In most cases it can be extended to
an $n$-complement on $S$.
In other cases we have $n(n+1)$-complements.
Under (M) or (SM), in the latter cases $n=1$,
and we have regular $n(n+1)=2$.

We can make it as in the proofs of
Lemma~\ref{indcomp}.\ref{indcomp1}
and Case~II.
We consider contractions of curves
$E\not\subseteq C$ with $(K+C+B'.E)\ge 0$ for
$B'=\lfloor (n+1)B\rfloor/n$.
In particular, $B'$ has multiplicities
in $\mathbb Z/n$.
Since $K+C+B=K+C+B'=K\equiv 0/Y$,
we can make $\nu$ extremal outside
a fibre $C$, i.~e., all other fibres 
are irreducible.
Then $K+C+B'$ is positive on fibres 
of some fibering $f:S\to Z$, or
$-(K+C+B')$ is nef.

Indeed, we can contract components of
$C$ preserving the numerical properties of
$K+C+B'$ because $K+C+B'\equiv 0/Y$ on
each component of $C$.
A terminal model will be extremal
and its cone has two extremal rays:
\begin{itemize}
\item the first one $R_1$ is generated by
a fibre $F$ of  $\nu$,
\item the second one $R_2$ is generated by
a multi-section section $E$.
\end{itemize}
If $(K+C+B'.E)> 0$, then $E$
induces a required fibering $f$.
We need to check that if $E$ is
contracted to a point, then we have
a required fibre contraction on $S$.
The former contraction induces
a birational contraction $f:S\to Z$,
for a birational inverse image of $E$.
After that $K+C+B'$ will be nef and big.
Moreover, it will also be positive 
on each curve $E\not\subseteq C$.
Again we can find extremal contractions
of such $E$ in $S$ subtracting $C$.
Finally, they give a required fibering,
because a terminal model cannot have
the Picard number $1$ by
Lemma~\ref{indcomp}.\ref{ncomp}.

In the case of such a fibre contraction (ruling), we
return to the original $S$.
By Lemma~\ref{indcomp}.\ref{ncomp},
$C$ will have a section of
induced $f:S\to Z$.
So, the horizontal
multiplicities of $B$ are in $\mathbb Z/(n+1)$
by Lemma~\ref{indcomp}.\ref{compl1}.
Then we use
Lemma~\ref{indcomp}.\ref{indcomp1}.
After contractions in a fibre of $C$ for $\nu$,
we have a fibre $C$ with the same $B_C$,
because $K+C+B\equiv 0/Y$.
On the other hand, by
Corollary~\ref{indcomp}.\ref{adjunct},
$B_C$ satisfies (M) and (SM)
whereas $\deg (K_C+B_C)=0$.
Thus by
Monotonicities~\ref{indcomp}.\ref{monl1}
and \ref{indcomp}.\ref{monl3},
$K_C+B_C$ is $n$-complementary if and only if
$B_C$ has index $n$.
Therefore $K_C+B_C$ is $m n$-complementary for any
natural number $m$.
Then, for $n:=(n+1)n$, we have (ii)-(iii) in
Lemma~\ref{indcomp}.\ref{indcomp1} on $S$.
Thus $K+C+B$ is $n(n+1)$-complementary.
Under (M) or (SM), $n=1$.

In other cases $-(K+C+B')$ is nef after
the above contractions.
We also assume that all fibres  of $\nu$,
except for fibre $C$, are irreducible.
By \cite[Theorem~6.9]{Sh2} after a complement
(cf. with the proof of Case~I),
$K+C+B'$ is log terminal, except for
the case when $C$ and $C'=\lfloor B'\rfloor$
are irreducible curves, and
$K+C+B'\equiv 0$.
Then $K+C+B'$ and $K+C+B$ are
$n$-complementary as in
the proof of Case~I.
So, we suppose that $K+C+B'$ is log terminal.

If $K+C+B'\equiv 0$, then
$B'$ has a fibre type, and
we have
an $n$-complement as in Case~II.
Indeed, we can use
Lemma~\ref{indcomp}.\ref{indcomp1}
when there exists a ruling $f:S\to Y$.
Otherwise, by the LMMP with
$K+C+B-\varepsilon D$ for an algebraic fibre 
$D=\nu^*(\nu(C))$, we
have a birational contraction of
a curve $E\not\subseteq C$ with
$(K+C+B'.E)=0$.
Since $E$ is a multi-section of $\nu$ and
$B'$ has a fibre type before the contraction,
it will have a big type afterwards.
Therefore, decreasing $B'$ after contraction,
we have an $n$-complement as in the big case
by Lemmas~\ref{indcomp}.\ref{monl2} and
\ref{indcomp}.\ref{mlem1}.

Finally, if $(K+C+B'.E)<0$ we
increase $B'$ adding vertical components,
but not in $C$, in such a way that
the new $B'$ will again have multiplicities
in $\mathbb Z/n$, and
$(K+C+B'.E)=0$.
Then $K+C+B'\equiv 0$, which
is the above case.
To verify $K+C+B'\equiv 0$, note
that $\nu$ is also a numerical
contraction for $K+C+B'$ with old $B'$,
and new $B':=B+\nu^* D$ for
an effective divisor $D$ on $Y$.

Before increasing $B'$,
we classify and choose
an appropriate model for fibre $C$.
We suppose that $C$ is minimal:
all $(-1)$-curves of $C$ are
contracted whenever $S$ is non-singular
near $C$.
It is possible only when $C$ is
reducible and the $(-1)$-curve
is not an edge in the chain of $C$.

Then fibre $C$ has one of
the following types
(see \cite{BPV} and
the Classification~\ref{lcomp}.\ref{class}):
\begin{description}
\item{($\rm I_b^*$)}
A minimal resolution of fibre $C$ has
a graph $\Dtd_{4+b}$, where $b+1\ge 1$ is
the number of irreducible components of $C$.
All curves of the resolved fibre are
$(-2)$-curves.
$B_C=(1/2)(P_1+P_2+P_3+P_4)$.
\item {(II)}
$B_C=(1/2)P_1+(2/3)P_2+(5/6)P_3$.
\item {(III)}
$B_C=(1/2)P_1+(3/4)(P_2+P_3)$.
\item {(IV)}
$B_C=(2/3)P_1+(2/3)P_2+(2/3)P_3$.
\end{description}
Curve $C$ is irreducible exactly
in the cases $(\rm I_0^*)$ and (II)-(IV).
Moreover, in cases II-IV, $C$ is $(-1)$-
or $(-2)$-curve, which splits
our classification into Kodaira's 
cases II-IV and II$^*$-IV$^*$,
respectively 
(cf. Classification~\ref{lcomp}.\ref{class}). 
Respectively, for I$_b^*$ and (II)-(IV), 
$n=2$ and $6,4,3$.
Type I$_b$ disappeared by our conditions on $C$.

When $C$ is a $(-1)$-curve,
we transform fibre $C$
of type (II)-(IV)
into a standard one $F_0$ of type II-IV 
in Kodaira's classification
\cite[p. 158]{BPV}.
Then $S$ is non-singular near $F_0$,
$F_0$ has multiplicity $1$ for $\nu$, and,
near $F_0$,
modified $C=(5/6)F_0, (3/4)F_0,$ and
$(2/3)F_0$ respectively.
Now the log singularity $C$ is hidden
in a point on $F_0$, and
new $C$ is not reduced.
By Lemma~\ref{indcomp}.\ref{mlem1},
we preserve the $n$-complements.
During the crepant modification all boundary
multiplicities have denominator $n$.

If $C$ is a $(-2)$-curve, then for types (II)-(IV) 
all curves of the resolved fibre are
$(-2)$-curves too.
Fibre $F_0=nC$ has, respectively,  
multiplicity $6,4$ and $3$.

In case $(\rm I_b^*)$, $F_0=2 C$.
This is an algebraic fibre. 

Now we choose $E$.
Since $K+C+B'$ is not nef,
there exists an extremal
contraction $f:S\to Z$ negative
with respect to $K+C+B'$.
It is not to a point.

First, suppose that $f$ gives a ruling
with a generic fibre $E$, i.~e.,
$E$ is a $0$-curve or $(K.E)=-2$.
Then $(C+B'.E)\le 2-1/n$, because
$E$ only crosses $C+B'$ in non-singular points,
and $(K+C+B'.E)\le -(1/n)$, when $<0$.
If $E$ is a section of $\nu$,
we can add to $K+C+B'$
a few copies of the generic fibre $F$
of $\nu$ with multiplicity $1/n$.

Otherwise $E$ is a double section of $\nu$,
except for cases II$^*$-IV$^*$, which we
consider later.
This follows from inequality
$(C+B'.E)\le 2-1/n$, because
$\nu$ has multiplicity $\le 2$ in $F_0$
and even $1$ in the cases II-IV.
On the other hand, 
$\mult{F_0}{C}\ge (n-1)/n\ge 2/3$
in the cases II-IV.
So, $E$ is a double section,
and, if $B'\not= 0$, then $(B'.E)\le 1/n$:
in case $\rm I_b^*$,
$(B'.E)=(C+B'.E)-(C.E)\le 1-1/n=1/2$,
and, in the cases II-IV,
$(B'.E)=(C+B'.E)-(((n-1)/n)F_0.E)\le
(2-1/n)-2(n-1)/n=1/n$.
Therefore such $B'=(1/n)F_1$ and $(F_1.E)=1$,
where $F_1\not= F_0$
is an irreducible fibre of $\nu$.
Then we increase $B'$ to $(2/n)F_1$.
If $B'=0$, we increase $B'$ by $(1/n)F$
in a generic fibre $F/Y$.

If $f$ is birational, $f$ contracts a $(-1)$-curve $E$
($(-1)$-curve on a minimal resolution).
Again in the cases I$_b^*$ and II-IV, 
$0>(C+K+B'.E)\ge -(1/n)$.
Moreover, $E$ is a section of $\nu$, and
we can increase $B'$ adding $(1/n)F$,
when $K+C+B'$ has index $n$ near $E$.
In case $\rm I_b^*$, $E$ passes through $P_i$.

In the other cases
$K+C+B'$ does not have index $n$
somewhere near section $E$.
It gives a singularity of $S$
on $E$ outside $C$.
On the above log model, $E$ has
another singularity $P_3$ of $S$.
(In the cases II-IV,
section $E$ crosses $F_0$ only
in a (single) non-singular point of $F_0$.)
Hence, by a classification of
surface log contractions,
there exists just the single
singularity outside $C$,
and we increase $B'$ only
in a corresponding fibre $F_1$.
After that
we need to check only that
$K+C+B'$ will have index $n$ near $F_1$.
After a crepant modification we
suppose that $S$ is minimal$/Y$
and non-singular near $F_1$.
On the other hand, $(K+C+B'.E)=0$
and $K+C+B'$ has index $n$ outside $F_1$.
Since $E$ is a section, $K+C+B$ has
index $n$ near $E$ or in a point $F_1\cap E$.
Using a classification of degenerations
due to Kodaira, we obtain that $K+C+B'$
has index $n$ everywhere in fibre $F_1$
and in $S$.

Finally, we consider the cases with II$^*$-IV$^*$.
When $f$ is a ruling, there exists a $(-1)$-curve
$E$ in a singular fibre of $f$.
In this case we choose $E$ through $P_3$.
If $f$ is birational again we have such a curve $E$.
In both cases $C$ passes through a singular point
$P_j\in C$ of $S$.
This time:
$0>(C+K+B'.E)\ge -(n-1)/n$.
More precisely,
if $E$ is passing through point
$P_j$ of type $\mathbb A_i$, then
$(C+K+B'.E)\ge -i/(i+1)$.
Since $E$ is not of a big type,
then on a minimal resolution it
can intersect only $(-m)$-curves of
the resolution with $m\ge i+1$ and $=$
only when $f$ is a ruling.
But then $(C+K+B'.E)\ge -i/(i+1)+(i-1)/(i+1)=-1/(i+1)$,
and $\ge -i/(i+1)(i+2)$ when $E$ is contractible. 
In particular, if $E$ is a section,
then $i+1=n$ and $-(n-1)/n(n+1)>-1/n$.
So, $E$ has at most two singularities of $S$
with the m.l.d. $<1$ (the m.l.d. is
the minimal log discrepancy),
and no intersections with $B'$ in other points
in such a case.
So, we can complete $B'$ as above when $E$ is a section.

Note also that in the ruling case, $E$ is
a section, except for the case when $E$ has only
one singularity $P_3$ of type $\mathbb A_3$, and
$E$ intersects the middle curve on a minimal resolution
of this point.
Then $n=4$,
$E$ is double section,
$B'$ intersects $E$ simply in 
a single branching point $Q\in E/Y$,
and $(K+C+B'.E)=1/2$ in $P_3$ and $1/4$ in $Q$.
Such intersection means the intersection of $E$
on a minimal crepant resolution with the boundary.
So, we complete $B'$ in the corresponding fibre.
Then $(K+C+B'.E)=1/2$ in $Q$ too for new $B'$. 
We need to verify that, for such $B'$, $K+B'$ 
will have index $n$.
This can be done on a Kodaira model $C'$.
If $E$ passes through a singular point on
Kodaira's model $C'$ after 
contractions to the central curves
as in our models $(I_b^*)$ and (II)-(IV), 
then for the fibre $C'$ of type I$_b^*$, 
$E$ passes through a simple Du Val singularity,
$B'=(1/2)C'$, and $K+(1/2)C'$ has index $4$.
Another possible type, assuming that $C'$ is
not reduced in $B'$, is only of type IV$^*$
with $B'=(3/4)C'$ and $E$.
Again $K+(1/2)C'$ has index $4$.
If $E$ passes transversally through
a non-singular point of $C'$,
then the multiplicity of $C'$ is $2$,
$B'=(1/2)C'$ and $K+(1/2)C'$ has index $4$.
If $E$ does not transversally pass $C'$,
the multiplicity of $C'$ is $1$, $(C'.E)=2$ 
near $C'$,
$B'=(1/4)C'$ and $K+(1/4)C'$ has index $4$.

In the other cases, $E$ is contractible,
and $E$ is an $l$-section with $2\le l\le 4$.
According to the classification of such
contractions, $E$ has a simple single
intersection with exceptional curve over
$P_j$, and it is an edge curve.
Moreover, $E$ has at most two singularities.
In the latter case $B'\not= 0$ only in two
corresponding fibres $C$ and $C'$, because
$-i/(i+1)(i+2)\ge 1/6$ and
$(K+C+B'.E)<0$.
In addition, $(K+C+B'.E)=1/(i+1)$ near $P_j$,
and $l=n/(i+1)$. 
We increase $B'$ to the numerically trivial case.
Near $C'$, $(K+C+B'.E)=i/(i+1)$ for new $B'$.
For such $B'$, we verify, that $K+B'$ has index $n$
near $C'$.
If $i=1$, $n=4$ and $l=2$, it was proved above
in the ruling case.
If $i=1$, $n=6$ and $l=3$, we can proceed 
similarly.
Since $E$ is $3$-section,
the type I$_b^*$ is only possible when
$B'=(1/3)C'$ with $E$ passing a singular point.
Then $K+B'$ has index $6$.
For type I$_b$, $B'=(1/6)C'$ when $C'$ has
the multiplicity $1$ and $B'=(1/2)$ when
$C'$ has the multiplicity $3$.
In both cases, $K+B'$ has index $6$.
The same holds for types II-IV
whereas the multiplicity of $C'$ is $1$.
In the case III$^*$, $B=(2/3)C'$ 
and $K+(2/3)C'$ has index $6$.
In the case IV$^*$, $B=(1/2)C'$
and $K+(1/3)C'$ has index $6$.

In the other case with two singularities, 
$i=2$, $n=6$ and $l=2$.
For type I$_b$, $B'=(2/3)C'$
when the multiplicity of $C'$ is $2$,
and $B'=(1/3)C'$ when the multiplicity is $1$.
The same hold for types II-IV
whereas the multiplicity of $C'$ is $1$.
For type I$_b^*$, $B'=(1/3)C'$
whereas $E$ does not pass through singularities.
Types II$^*$-IV$^*$ are impossible 
in this situation.
In all these cases $K+C+B'$ has index $6$.

Finally, $P_j$ is the only singularity of
$S$ on $E$.
Then it is enough to construct $B'$ near $E$.
So, it is possible to do this when
$E$ has the only (maximally branching) point 
with simple intersection with
a fibre somewhere over $Y$.
On the other hand, if $B'\not =0$ is in a fibre  $C'$
with a non-branching point, then,
near $C'$, $B'\ge (1/n)C'$ and $(B'.E)\ge l/n\ge 2/n$.
So, $n\ge 4$, and, for $n=4$, $(K+C+B'.E)=1/4$
in $P_j$, and $l=2$ is impossible.
So, $n=6$, $(K+C+B'.E)=1/3$ in $P_j$, 
and we can increase $C'$ up to $(1/3)C'$,
when $B'$ does not have other components.
Otherwise we have a third fibre $C''$
with a simple only intersection 
and a $6$-complement.

In all other cases we have two branchings of $E/Y$
in one fibre $C'$, $l=4$, and these
branchings are in a fibre $C'$ with $B'\ge (1/n)C'$.
(Cf. Lemma~\ref{indcomp}.\ref{invcomp}.)
Then $n=6$, $P_j=P_2$, $(K+C+B'.E)=2/3$
in $P_2$, and this case is impossible. 
\end{proof}

\begin{proof}[Proof of the Inductive Theorem:
Case~IV]
Here we suppose that $K+C+B\equiv 0$ and
$B$ has an exceptional type.
Then $B_C$ satisfies (SM), by
Lemma~\ref{indcomp}.\ref{adjunct}
after a contraction of $B$.
On the other hand,
$\deg (K_C+B_C)=(K+C+B.C)=0$.
Thus we have an $n$-complement on $C$ and
near $C$ (exactly) for such $n$ that
$B_C$ has index $n$.
Let us take such $n$.

Therefore $K+C+B=K+C+B'\equiv 0$
and has index $n$ near $C$,
where $B'=\lfloor (n+1)B\rfloor/n$.
This follows from \cite[Proposition~3.9]{Sh2},
when $B$ is contracted.
To establish that $K+C+B$ gives
an $n$-complement, we need to check
that $K+C+B$ has index $n$ everywhere in $S$.

After a contraction of $B$ we assume
that $C+B=C$.
Using the LMMP for $K$, we reduce
the situation to the case when
there exists an extremal fibre contraction
$f:S\to Z$.

If $F$ is to a curve $Z$,
we use Lemma~\ref{indcomp}.\ref{indcomp2}.
Perhaps we change $n$ by $2n$
when $n=1$ or $3$.

Otherwise $Z=\pt$,
$S$ has the Picard number $1$, and $C$ is ample.
Then we use a covering trick
with a cyclic $n$-covering $g:T\to S$.
On $T$, $K_T+D=g^*(K+C)$ has
index $1$ near ample $D=g^{-1}C$.
By \cite[Proposition~3.9]{Sh2},
$D$ is a non-singular curve (of genus $1$)
and $T$ is non-singular near $D$.
Hence $T$ is rational
and $K_T+D\sim 0$ (cf. proof of Case~I).
Therefore, $K+C+B$ has index $n$.
\end{proof}

\begin{proof}[Proof of \ref{indcomp}.\ref{indth}.2:
Global case]
From the proof of the Inductive Theorem.
\end{proof}

\section{Local complements}
\label{lcomp}

In the {\it local case\/},
{\rm when $Z\ge 1$\/},
we can drop most of the assumptions
in the Main and Inductive Theorems,
and in other results.

\refstepcounter{subsec}
\label{lmainth}
\begin{theorem}
Let $(S/Z,C+B)$ be a surface log contraction
such that
\begin{itemize}
\item
$\dim Z\ge 1$, and
\item
$-(K+C+B)$ is nef.
\end{itemize}
Then it has
locally$/Z$ a regular complement, {\rm i.~e.,
$K+C+B$ has $1-,2-,3-,4-$ or $6$-complement\/}.
\end{theorem}

\begin{proof}[Proof of the
Theorem~\ref{lcomp}.\ref{lmainth}: Special case]
Here (NK) of the Inductive Theorem
is assumed
in the following strict form:
\begin{description}
\item{(MLC)} ({\it maximal log canonical\/})
$K+C+B$ is not Kawamata log terminal
{\it in\/} a fibre$/P$, near which
we would like to find a complement, i.~e.,
$C$ has a component in the fibre, or
an exceptional divisor, with
the log discrepancy $0$ for $K+C+B$, has
the center in the fibre. 
\end{description}
But as in
Theorem~\ref{lcomp}.\ref{lmainth},
we drop other assumptions in
the Inductive Theorem except for (NEF).

After a log terminal blow-up, we
assume that $K+C+B$ is log terminal,
and by our assumption $C\not=0$.
Moreover, $C$ has a component in
a fibre of $f$.

According to the Big case,
we suppose that $Z$ is a non-singular curve,
and $K+C+B\equiv 0/Z$.
Hence $f$ is a fibering of log curves
of genus $1$.
Note that $C$ is connected near the fibre
by the LMMP and \cite[Lemma~5.7]{Sh2}.

Let $B_C=(B)_C$ be the different
for adjunction $(K+C+B)|_{C}=K_{C}+B_C$.
Then $K_C+B_C\equiv 0/P$.
Take a regular $n$ such that
$K_C+B_C$ is $n$-complementary near the fibre. 
We contend that $K+C+B$ is
$n$-complementary near the fibre too.

Divisor $-(K+C+B')$ is nef on
a generic fibre of $f$ where
$B'=\lfloor (n+1)B\rfloor$.
Indeed, by the LMMP and
Proposition~\ref{indcomp}.\ref{ncomp},
the same holds for the fibre after
a contraction of non-reduced components
of $B$ in the fibre$/P$.

We decrease $C+B$ in
a horizontal component whenever one exists.
Then, by Lemma~\ref{indcomp}.\ref{monl2} and
\ref{indcomp}.\ref{monl2}.1,
for an appropriate choice of
the horizontal component,
we have the same $n$-complements
and use the Big case.

So, we can assume that $C+B$ has
no horizontal components.
Thus $f$ is a fibering of curve of genus 1.

Then, as in Case~IV in
the proof of the Inductive Theorem,
$B$ has an exceptional type.
As there we can verify that
$B_C$ satisfies (SM),
$K_C+B_C\equiv 0$, and
$K_C+B_C$ has index $n$.
Moreover, $K+C+B$ has (local) index $n$
near the fibre$/P$, and $B'=B$.

Therefore, to check that $K+C+B$
has a (trivial) $n$-complement
near the fibre, we need to check only
that $K+C+B$ has a {\it global} index $n$, i.~e.,
\begin{equation}
\label{indn}
n(K+C+B)\sim 0/P.
\end{equation}

The log terminal contractions, i.~e.,
contractions of the components in $B$,
preserve the (formal) log terminal property
of $K+C+B$ and (\ref{indn})
according to \cite[2.9.1]{Sh1} \cite[3-2-5]{KMM}.

Thus we assume after contractions that
$C+B=C$ is the fibre$/P$ and
$K+C$ is (formally) log terminal.
Such a model of the fibre is
its {\it weak log canonical model\/}.
It is not unique.
For instance, we can blow up
a nodal point of $C$.
However, it is unique if
we impose the following
{\it minimal\/} property:
\begin{itemize}
\item
all $(-1)$-curves of $C$ are
contracted whenever $S$ is non-singular
near $C$.
\end{itemize}
Such a model will be called
{\it log minimal}.
Its uniqueness follows
from the MMP and
a classification below
(see Classification~\ref{lcomp}.\ref{class}).
Note also that it can be non log terminal,
but it is always formally log terminal.

We check (\ref{indn}) for
the minimal $n$ such that $K_C+B_C$ has index
$n$, i.~e., for the index of $K_C+B_C$.

This essentially
follows from Kodaira's classification of elliptic fibres 
\cite[Section 7 of Ch. V]{BPV} and his formula for
a canonical divisor of an elliptic fibering
\cite[Theorem~12.1 in Ch. V]{BPV} \cite{Shf}.
See also the Classification~\ref{lcomp}.\ref{class}
below.
The latter gives a (non-standard) classification of
the degenerations
for a fibering with the generic curve of genus 1.
\end{proof}

First, we add types $_m\rm I_b,b\ge 0,$ to
types of log models in the proof of Case~III
in the Inductive Theorem.
In the following cases,
$S$ is non-singular near $C$.
\begin{description}
\item{($_m\rm I_0$)}
$C$ is a non-singular
curve of genus $1$, and $f^*P=m C$.
\item{($_m\rm I_1$)}
$C$ is an irreducible rational curve
of genus $1$ with one node,
and $f^*P=m C$.
\item{($_m\rm I_b$)}
$C$ is a wheel of $b\ge 2$
irreducible non-singular rational curves $C_i$,
and $f^*P=m C$.
Each $C_i$ is a $(-2)$-curve.
\end{description}

\refstepcounter{subsec}
\label{class}
\theoremstyle{plain}
\newtheorem*{class}{\thesection.\thesubsec.
Classification of degenerations
in genus $1$ (Kodaira)}
\begin{class}
Any degeneration of non-singular curves
of genus $1$ has up to a birational transform
a log minimal model of one of
the following types:
$_m\rm I_b$, $_m\rm I_b^*$, 
{\rm II, II$^*$, III, III$^*$}, or 
{\rm IV, IV$^*$}.
Each of these models has
a unique birational transform into
a Kodaira model with the same label.

In addition, for the log model of type
$_m\rm I_b$, $_m\rm I_b^*$, 
{\rm II, II$^*$, III, III$^*$}, or 
{\rm IV, IV$^*$},
$K+C$ has, respectively, the index
$1,2,6,4$ or $3$.
\end{class}

Note that $K+C$ is log terminal,
except for the type $_m\rm I_1$,
when $C$ is a Cartesian leaf.

\begin{proof}
Adding a multiplicity of the fibre
we can suppose (MLC) (cf. below Proof of
Theorem~\ref{lcomp}.\ref{lmainth}: General case).

Then, according to the proof the Special case in
Theorem~\ref{lcomp}.\ref{lmainth},
we have a log minimal model $C/P$.

According to a classification of
formally log terminal singularities,
$C$ is a connected curve with only
nodal singularities.
On the other hand $(C,B_C)$ has
a log genus $1$.
Thus $C$ has an arithmetic genus $\le 1$.

If the genus is $0$, then
$C$ is a chain of rational curves,
and the possible types were given
in the proof of Case~III
in the Inductive Theorem.

If $C$ is not irreducible then,
by (SM) and since $K_C+B_C\equiv 0$,
$B_C$ is the same as in the type $\rm I^*_b$.
So, each $P_i$ is a simple double singularity of $S$,
and a minimal resolution gives
a graph $\Dtd_{4+b}$.
This fibre has the type $\rm I_b^*$
due to Kodaira.
In that case, $n=2$, $f^*P=2 C$, and
$K\sim 0/P$.
Thus $2(K+C)\sim 0/P$.

The same holds when $C$ is irreducible,
and $B_C$ is the same as in the type $\rm I^*_b$.

In the other cases with the genus $0$,
$C$ is an irreducible non-singular rational curve.
Then $\deg B_C=2$ and under (SM)
we have only $B_C$ as in
types $\rm I^*_b$, (II), (III), or (IV) as above.
We need to consider only the types (II)-(IV).
In all of them we have three singularities
$P_i$.
If $C$ is $(-1)$-curve, 
each of them is simple, i.~e., has
a resolution with one irreducible curve.
Otherwise $C$ is $(-2)$-curve and
the singularities are Du Val.
This gives respectively
Kodaira's types II-IV and II$^*$-IV$^*$.
For instance, the curves, in a minimal resolution
of points $P_i$ in type~IV,
are $(-3)$-curves.
Now we can easily transform fibres of
the types~II, III and IV into the same due Kodaira,
because $C$ is a $(-1)$-curve.
In type~III, $C$
will be transformed into three $(-2)$-curves
with a simple intersection in a single point.
The latter is a blow-down of the old $C$.
(Cf. \cite[p. 158]{BPV}. 
After the transform,
$C=((n-1)/n)F_0$,
where $F_0$ is the modified fibre
(cf. Proof of the Inductive Theorem:
Case~III).

For types I$_b^*$ and II$^*$-IV$^*$,
the transform is a minimal resolution.

Then, for Kodaira's types~II-IV, 
and II$^*$-IV$^*$,
$K\sim 0/P$ and $F_0\sim 0/P$.
(More generally, $D\sim 0/P$
for any integral $D$ such that $D\equiv 0/P$,
and for the types $\rm I_b^*$, II-IV
and II$^*$-IV$^*$.
This follows from the fact that
$F_0$ does not have
non-trivial unramified coverings.)
Hence, on Kodaira's model,
$n(K+C)=n(K+((n-1)/n)F_0)=nK+(n-1)F_0\sim 0/P$.
An alternative approach will be discussed
at the end of the proof.

In the other cases $B_C=0$, and
the genus is $1$.
Since $K+C$ is formally log terminal,
$S$ is non-singular near $C$ \cite[3.9.2]{Sh2},
$C$ is a curve with only nodal singularities
and of arithmetic genus $1$.
In particular, $n=1$.

If $C$ is irreducible, then for some natural number $m$,
$f^*P=m C$, and $K\sim (m-1)C/P$ by
Kodaira's formula \cite[ibid]{BPV}
(cf. formula (\ref{adj}), with $B_Y=P$,
in Section~\ref{indcomp}).
So, $K+C\sim 0/P$.

Similarly, we can do the next case,
when $C$ is reducible.
Since $K_C+B_C\equiv 0$,
the irreducible components $C_i$ of $C$
form a wheel as in $_m\rm I_b$.
By the minimal property, $K/P$ is nef.
On the other hand, $K\equiv 0$ in
a generic fibre. 
Hence $K\equiv 0/P$, each $C_i$ is
a $(-2)$-curve, and $f^*P=m C$.
So, in this case, a log minimal model $C/P$
coincides with a Kodaira model of
the type $_m\rm I_b$,
$K\sim (m-1)C/P$,
and $K+C\sim 0/P$.

Finally, we may also use a covering trick
\cite[Section~2]{Sh2} to reduce
a proof of (\ref{indn}) to a case with $n=1$
or to the type $_m\rm I_b$.
The latter is a {\it crucial\/} fact:
$K+C\sim 0/P$ for the type $_m\rm I_b$.
It can be induced from dimension $1$.
\end{proof}

\begin{proof}[Proof of
Corollary~\ref{int}.\ref{ind2}:
Local case]
This is the Special case because
any local trivial complement satisfies (MLC).
On the other hand,
any regular $n$ divides $I_1=12$.
\end{proof}

\begin{proof}[Proof of \ref{indcomp}.\ref{indth}.2:
Local case]
From \cite[Example~5.2.2]{Sh2},
because the complements are induced
from the $1$-dimensional case.
\end{proof}

\begin{proof}[Proof of the
Theorem~\ref{lcomp}.\ref{lmainth}: General case]
According to the Big case
in the Inductive Theorem,
we suppose that $Z$ is a non-singular curve,
and $K+C+B\equiv 0/Z$.

By \cite[Lemma~5.3]{Sh2} we can
increase $B$.
We do it in such a way that
$K+C+B+p f^*P$ is {\it maximally log canonical\/}
for some real $p\ge 0$:
\begin{itemize}
\item
$K+C+B+p f^*P$ is log canonical, but
\item
$K+C+B+p' f^*P$ is not so for any $p'>p$.
\end{itemize}
Such a $p$ exists, and (MLC)
is equivalent to these conditions.
\end{proof}

\begin{proof}[Proof of
Corollary~\ref{int}.\ref{mcorol},
Main and Inductive
Theorems: Local case]
Follows from
Theorem~\ref{lcomp}.\ref{lmainth}.
\end{proof}

\section{Global complements}
\label{gcomp}

\refstepcounter{subsec}
\begin{theorem}
\label{gmainth}
Let $(S,C+B)$ be
a complete algebraic log surface
such that
\begin{description}
\item
(M) of the Main Theorem holds, and
\item{{\rm (NEF)}} $-(K+C+B)$ is nef.
\end{description}
Then its complements
are bounded under
any one of the following conditions:
\begin{description}
\item{{\rm (WLF)}}
of Conjecture~\ref{int}.\ref{conjcs};
\item{{\rm (RPC)}}
of the Inductive Theorem;
\item{{\rm (EEC)}} of the Inductive Theorem;
\item{{\rm (EC)+(SM)}}
of Conjecture~\ref{int}.\ref{conjcs};
\item{{\rm (ASA)}}
of the Inductive Theorem; or
\item{{\rm (NTC)}}
of the Inductive Theorem.
\end{description}

More precisely, for almost all
such $(S,C+B)$,
we can take a regular index in $RN_2$.
The non-regular complements are
exceptional in the sense of
Definition~\ref{int}.\ref{dexc}.
\end{theorem}

\refstepcounter{subsec}
\label{unb}
\begin{lemma}
There exists $c>0$ such that
all $b_i\le 1-c$, for
any log surface $(S,B)$,
under the assumptions of
Theorem~\ref{gcomp}.\ref{gmainth},
and such that $\rho(S)=1$,
$S$ is $1/7$-log terminal, and
$(S,B)$ does not have
the regular complements.
\end{lemma}

\begin{proof}
If $B=0$, any $c>0$ fits.
Otherwise,
$S$ is a log Del Pezzo.
Such an $S$ is bounded by \cite[Theorem 6.9]{Al}.
By (M) the same holds for $(S,\Supp{B})$.

So, we may assume that $S$ is fixed, as are
the irreducible components of $\Supp B=\cup D_i$.
Consider a domain
$$\mathcal D=
\{D=\sum d_i D_i\mid K+D \text{ is log canonical and }
-(K+B) \text{ is nef}\}.$$
It is a {\it closed\/} polyhedron by
\cite[Property~1.3.2]{Sh2}
and  by a polyhedral property of $\NEc{S}$.
Take $c=1-d$ where $d=\max \{d_i\}$ for
$D\in\mathcal D\cap\{D\ge B\}$.
\end{proof}

\begin{proof}[Proof: Strategy]
We are looking for the exceptions.
Thus we assume that $(S,B)$
does not have $1-,2-,3-,4-$ and $6$-complements.
We prove that $(S,C+B)$ belongs to a bounded family.
Equivalently, $(S,\Supp{(C+B)})$ is bounded.
Moreover, we verify that complements are
bounded and exceptional as well.

We will suppose (ASA) or (NTC)
by Proposition~\ref{indcomp}.\ref{relcond}.
For the exceptions in the proposition, we
have regular complements.

According to~\cite[Theorem~2.3]{Sh3} and
\cite[Lemma~5.4]{Sh2},
we can suppose that $(S,C+B)$ is log terminal.
In particular, $S$ has only rational singularities.
So, $S$ is projective.
Moreover, then $(S,C+B)$ is Kawamata
log terminal by the Inductive Theorem.
In particular, $C=0$.
The change preserves (ASA).

In addition, we suppose that
$K+B$ is $1/7$-log terminal in
the closed points of $S$.
Otherwise we make a crepant blow-up
of the exceptional curves $E$ with
a log discrepancy $\le 1/7$.
This preserves all our assumptions.
We have a finite set of such $E$
by \cite[Corollary~1.7]{Sh3}.

In other words, now $K+B$ is
$1/7$-log terminal in the closed points.

If $K+B$ is $1/7$-log terminal everywhere,
or, equivalently, if $B$ does not have
an irreducible component $D_i$
with $\mult{D_i}{B}\ge 6/7$ and satisfies (SM),
then $S=(S,B)$ is bounded
by (M) and \cite[Theorem~6.9]{Al},
except for the case, when $B=0$ and
$S$ has only canonical singularities.
In the former case we have
a bounded complement.
If $(S,B)$ satisfies (WLF), we
construct a complement as
in \cite[Proposition~5.5]{Sh2}.
Similarly, we proceed in the other cases
by (ASA).
Since $(S,B)$ is bounded in a strict sense, i.~e.,
in an algebraic moduli sense,
a freeness of $-(K+B)$ is bounded.
In the case when $B=0$ and $S$
has only canonical singularity,
$(S,0)$ has a regular complement
according to (ASA) and to a classification
of surfaces.
In such a case,
we can even suppose that $S$ is non-singular.

Now we can assume that $B$ has
an irreducible component $D_i$
with $\mult{D_i}{B}\ge 6/7$.
Then we reduce all required boundednesses
to a case with the minimal
Picard number $\rho=\rho(S)=1$.
We find a birational contraction
$g: S\to S_{\rm min}$ such that
$S_{\rm min}$ has all the above properties
and $\rho(S_{\rm min})=1$.
Moreover, $g$ does not contract
the irreducible components $D_i$
with $\mult{D_i}{B}\ge 6/7$, and
\begin{description}
\item{\rm (BPR)}
there exists a boundary $B'\ge B$,
with $\Supp{B'}$ in divisors $D_i$
having $\mult{D_i}{B}\ge 6/7$, such
that $g$ contracts only curves $E$
with log discrepancies $\le 1$ for
$K_{\rm min}+B'_{\rm min}$, and
$-(K_{\rm min}+B'_{\rm min})$ is nef,
where $K_{\rm min}=K_{S_{\rm min}}$ and
$B'_{\rm min}=g(B')$.
\end{description}
In particular, $B_{\rm min}=g(B)\not=0$.
This reduction will be called
a {\it minimization\/}.
It uses the Inductive Theorem and
the Main Lemma below.

By the LMMP $-(K_{\rm min}+B_{\rm min})$ is nef.
Hence $B_{\rm min}$ and $-K_{\rm min}$ are ample,
because $\rho(S_{\rm min})=1$ and $B_{\rm min}\not= 0$.
So, $S_{\rm min}$ is a log Del Pezzo surface.
Since $K_{\rm min}+B_{\rm min}$ is
$1/7$-log terminal in the closed points,
then, by \cite[Monotonicity~1.3.3]{Sh2},
$S_{\rm min}$ does the same.
Therefore, due to Alekseev,
we have a bounded family
of such Del Pezzo surfaces
\cite[Theorem 6.9]{Al}.
For $\Supp{B_{\rm min}}$ we have only a bounded
family of possibilities,
because all $b_i\ge 1/2$ and $\rho(S_{\rm min})=1$.

The condition (BPR) above
guarantees a
{\it boundedness for the partial resolution\/} $g$.
First, by the Inductive Theorem
$K_{\rm min}+B'_{\min}$ is Kawamata log terminal and
$B'_{\rm min}$ is reduced, because
$K_{\rm min}+B'_{\rm min}$ does not have
the regular complements as
$K_{\rm min}+B_{\rm min}$.
Hence the multiplicities
of $B'_{\rm min}$ and $B$
are {\it universally\/} bounded
by Lemma~\ref{gcomp}.\ref{unb} and (M):
all $b_i\le 1-c$ for some $c>0$.
Thus $(S,\Supp{B})$ is bounded,
because it resolves only
exceptional (for $S_{\rm min}$)
divisors $E$ with log discrepancies $\le 1$
for $K_{\rm min}+B'_{\rm min}$ by (BPR)
(cf. \cite[Second Main Theorem and
Corollary~6.22]{Sh3}).

This will be done more explicitly 
in Theorem~\ref{ecomp}.\ref{emainth} below.
For another approach see 
the following remark.
\end{proof}

\refstepcounter{subsec}
\begin{remark}
In the strategy above, $(S,\Supp{B})$
is bounded by \cite[Theorem~6.9]{Al}
and Lemma~\ref{gcomp}.\ref{unb}.
Indeed, $K+B$ is $\varepsilon$-log terminal
for any $c>\varepsilon>0$.
However, we prefer a more effective and explicit
property (BPR)
(cf. Proof of
Theorem~\ref{ecomp}.\ref{emainth}
in Section~\ref{ecomp}).
\end{remark}

In the same style
as Lemma~\ref{indcomp}.\ref{mlem1},
we can prove its following improvement.

\refstepcounter{subsec}
\label{mlem3}  
\theoremstyle{plain}
\newtheorem*{mainl}{\thesection.\thesubsec.
Main Lemma}
\begin{mainl}
In the notation of \cite[Definition~5.1]{Sh2},
let $f\colon X\to Y$ be
a birational contraction such that
\begin{description}
\item{\rm (i)}
$K_X+S+D$ is numerically non-negative on a
{\rm sufficiently general\/} curve$/Y$ in each
exceptional divisor of $f$; and
\item{\rm (ii)} 
for each multiplicity $d_i=\mult{D_i}{D}$ of a prime divisor
$D_i$ in $D$,
$\lfloor (n+1)d_i\rfloor/n\ge d_i$,
whenever a non-exceptional on $Y$ divisor $D_i$ 
intersects an exceptional divisor of $f$.
\end{description}
Then
$$K_Y+f(S+D)\ n-complementary\ \ \Longrightarrow
\ K_X+S+D\ \ n-complementary.$$

In addition, we can assume that $D$ is
just a subboundary.
\end{mainl}

\refstepcounter{subsec}
\label{mlem4}
\begin{example-c} 
By the Monotonicity Lemma \ref{indcomp}.\ref{monl1},
(ii) in \ref{gcomp}.\ref{mlem3} holds
whenever all coefficients are standard,
i.~e., they satisfy (SM). 

Respectively, (i) in \ref{gcomp}.\ref{mlem3}
holds when
$K+S+D$ is nef$/Y$.

Then by the Main Lemma
we can pull back the complements, i.e.,
for any integer $n>0$,
$$K_Y+f(S+D)\ \ n-complementary\ \ \Longrightarrow
\ K_X+S+D\ \ n-complementary.$$
\end{example-c}

\begin{proof}[Proof of the Main Lemma]
We take a crepant pull back:
$$K_X+D^{+X}=f^*(K_Y+D^+).$$
It satisfies \cite[5.1.2-3]{Sh2} as
$K_Y+D^+$, and we need to check
\cite[5.1.1]{Sh2} only for the exceptional
divisors.
For them it follows from our assumption and
the Negativity~\ref{indcomp}.\ref{negat}.
Indeed, on the exceptional prime divisors $D_i$,
$D^{+X}\ge D$ and has multiplicities in $\mathbb Z/n$.
Hence $D^{+X}\ge S+\lfloor (n+1)D\rfloor/n$
according to 
the Monotonicity~\ref{indcomp}.\ref{monl3} above  
and \cite[Lemma~5.3]{Sh2}.
Indeed, for any multiplicity $d_i^+<1$ in
$D^{+X}$, we have
$d_i^+\ge\lfloor (n+1)d_i^+\rfloor/n\ge
\lfloor (n+1)d_i\rfloor/n$.
\end{proof}

\begin{proof}[Proof of
Theorem~\ref{gcomp}.\ref{gmainth}:
Minimization]
Let $D$ denote a boundary
with the coefficients
$$d_i=\cases
1&\text{ if }b_i\ge 6/7,\\
b_i
&\text{ otherwise.}\endcases$$
Hence by
the Monotonicity Lemma~\ref{indcomp}.\ref{monl1},
for any $n\in RN_2$,
we have
\begin{itemize}
\item
$\lfloor(n+1)B\rfloor/n=
\lfloor D\rfloor+\lfloor(n+1)\{D\}\rfloor/n\ge D\ge B$.
\end{itemize}
Therefore, $K+D$ is log canonical.
Indeed, locally by \cite[Corollary~5.9]{Sh2},
there exists an $n$-complement $(S,B^+)$
with $n\in RN_2$.
In addition, $B^+\ge \lfloor(n+1)B\rfloor/n\ge D$.
Hence by \cite[Monotonicity~1.3.3]{Sh2},
$K+D$ is log canonical.

Since $B$ has a multiplicity $b_i\ge 6/7$,
then $D$ has a non-trivial reduced part
and $K+D$ is not Kawamata log terminal.

By the Inductive Theorem,
$-(K+D)$ does not satisfy (ASA), because
$K+D$ as $K+B$ does not have
the regular complements.

Moreover, we contend that when $\rho=\rho(S)>1$,
then for any $\mathbb R$-divisor $F$, such that
$D\ge F\ge B$ and $-(K+F)$ is semi-ample,
there exist
an exceptional curve $E$ and
a divisor $B'$ such that
\begin{itemize}
\item
$(K+D.E)>0$ and $\mult{E}{B}\le 5/6$,
\item
$D\ge B'\ge F$, $(K+B'.E)=0$, and
$-(K+B')$ is semi-ample.
\end{itemize}
Then we contract $E$ to a point:
$h:S\to Z$, and replace $(S,B)$ by
$(Z,h(B))$.
On the first $S$ we take $F=B$.
Then we take $F=h(B')$.
The contraction preserves the properties.
In particular, $(Z,h(B))$ does not have
the regular complements by
the Main Lemma.
We contract only curves
with $(K+B.E)\le 0$, and,
by the Local case and (M),
with $b_i\le 5/6$.
Indeed, near $E$ we have
a regular complement $(S,B^+)$, $B^+\ge D$,
$\mult{E}{B^+}=\mult{E}{D}=1$, and
$(K+D.E)\le (K+B^+.E)=0$.
Hence, $K+B$ will always be 
$1/7$-log terminal, and
we do not contract the curves with $b_i\ge 6/7$.
Contracted $E$, or any other exceptional
divisor of $S$ with a log discrepancy $\le 1$
for $K+B'$,
will have the same log discrepancy
for $K_Z+h(B')$.
By \cite[1.3.3]{Sh2}, these discrepancies do not
increase for $K+F'$ with any $F'\ge h(B')$.
Thus all contracted $E$ will have
log discrepancies $\le 1$ for $K+B'$.
Finally, an induction on $\rho$
gives required $S_{\rm min}=S$ with $\rho=1$.

We find $E$ case by case with respect
to the numerical dimension $\kappa^*$
of $-(K+B)$.

First, (WLF) when $\kappa^*=2$.
Then we have (RPC).
In particular, $-(K+B)$ is not nef,
because it is not semi-ample.
Then there exists an exceptional curve $E$
with $(E.K+D)>0$.
Since $\rho>1$, otherwise we
have a fibre extremal contraction $S\to Z$,
which is positive with respect to $K+D$.
The latter is impossible by (M),
because $-(K+B)$ is nef.
Therefore,
we need to find $B'$ and $E$
with above properties.
Take a closed polyhedron
$$
\mathcal D=\{B'\mid
D\ge B'\ge F,\text{ and}
-(K+B')\text{ is nef}\}.
$$
It is polyhedral by (RPC).
Take a maximal $B'$ in $\mathcal D$.
Then $-(K+B')$ satisfies (ASA).
It is Kawamata log terminal by
the Inductive Theorem, because $B'\ge B$.
So, we cannot increase $B'$ only 
because $(K+B'.E)=0$ for
some extremal curve and is positive when
we increase $B'$.
By (M) it is possible only
for birational contractions.
Properties (WLF) and (RPC) will be
preserved.

Second, $\kappa^*=1$, and we
have a numerical contraction $\nu :S\to Y$
for $K+B$.
By (M) the horizontal multiplicities
satisfy (SM), and $D=B$ in
the horizontal components.
Thus $-(K+D)\equiv 0/Y$.
Otherwise we have a vertical exceptional curve $E$
with $(K+D.E)>0$.
As above we contract $E$.
This time, we can take $B'=B$,
because $K+B\equiv 0/Y$.
After such contractions,
$-(K+D)\equiv 0/Y$ and is not nef,
because it is not semi-ample.
Note that $Y$ is rational, because
$K+B$ is negative on the horizontal curves.
As above we have no extremal fibre contractions,
positive with respect to $K+D$.
Thus we have an exceptional
(horizontal) curve $E$
with $(E.K+D)>0$.
After that contraction we will have (WLF)
and do as above.

Third, $\kappa^*=0$ or
$K+B\equiv 0$.
In this case we take $B'=B$
and need only to contract
some $E$ with $(K+D.E)>0$.
If $B$ has a big type,
we again have (WLF) and (RPC).
If $B'$ has a fibre type,
then, by Proposition~\ref{indcomp}.\ref{types},
we have a fibration $S\to Y$ of
genus $1$ curves, whereas $B'$ and $D$ have only
vertical components.
As above, after contractions,
we suppose that $K+D\equiv 0/Y$.
Since $B\not=0$ and forms a fibre,
we have a horizontal extremal curve $E$
with the required properties.
After its contraction, we have (WLF).
Finally, $B$ has an exceptional type.
Then decreasing $B$ in the non-standard
multiplicities, we can find $E$, which
is outside $\lfloor D\rfloor$, but
intersects $\lfloor D\rfloor$.
Thus $(K+D.E)>0$.
If after a contraction of such $E$, we change
a type of $B$, we return to
a corresponding type: big or fibre.
\end{proof}

\begin{proof}[Proof of
Theorem~\ref{gcomp}.\ref{gmainth}:
Bounded complements]
Here, we check that complements
are bounded.
Since $(S,\Supp{B})$ is bounded,
it is enough to establish that
complements, for all
$$
B'\in \mathcal D=\{
\Supp{B'}=\Supp{B},
-(K+B')\text{ is nef and log canonical}\},
$$
are bounded.
Note that each $K+B'$ is semi-ample
by Proposition~\ref{indcomp}.\ref{types},
because $K+B$ is semi-ample,
Kawamata log terminal, and
(NV) of Remark~\ref{indcomp}.\ref{rel}
holds for $K+B'$.
Thus, for each $\mathbb Q$-boundary $B'$
we have an $n$-complement such that
$n B'$ is integral.
Therefore we have $n$-complements
near $B'$ by
the Monotonicity~\ref{indcomp}.\ref{monl2}.

Hence we have bounded complements,
according to Example~\ref{int}.\ref{simul}.
Indeed, we can restrict our problem on 
any ample non-singular curve; as
we see later in Section~\ref{ecomp},
the cases with non-standard coefficients
are reduced to a case with $\rho(S)=1$. 

A more explicit approach will be
given in Theorem~\ref{ecomp}.\ref{emainth}.
\end{proof}

\begin{proof}[Proof of
Theorem~\ref{gcomp}.\ref{gmainth}:
Exceptional complements]
As in the strategy we assume (ASA),
a log terminal property for $K+B$
and an absence of the regular complements.
Then we have a (non-regular) complement
$(S,B^+)$.
Here we check that $K+B^+$ is
Kawamata log terminal for any (such) complement.

After a crepant blow-up we
suppose that $K+B^+$ has a reduced
component, and derive a contradiction.
Let $D$ denote a boundary
with the coefficients
$$d_i=\cases
1&\text{ if }b_i^+=1,\\
6/7&\text{ if }1>b_i^+\ge 6/7,\\
b_i
&\text{ otherwise.}\endcases$$
Then $B^+\ge D$ and $D$ satisfies (SM).
By the Monotonicity~\ref{indcomp}.\ref{monl1},
in the non-reduced components of $D$
and $B^+$,
$\lfloor(n+1)B\rfloor/n
\le \lfloor(n+1)D\rfloor/n$, for any $n\in RN_2$.
Hence,
by \cite[Lemmas~5.3-4]{Sh2},
$(S,D)$ does not have
the regular complements.
Hence $-(K+D)$ does not satisfy
(ASA) by the Inductive Theorem.

We contend then that $\rho>1$
and we have an exceptional curve
$E$ with $(K+D.E)>0$, and
automatically $\mult{E}{D}<1$.
Indeed, if $\rho=1$, then
$K+B^+\equiv 0$, log canonical,
and $B^+\not=0$, $K$ are ample.
Hence $-(K+D)$ is nef,
because $B^+\ge D$.
This is impossible by the Inductive theorem.

Therefore, $\rho>1$.
If we have an exceptional curve
$E$ with $(K+D.E)>0$,
we contract this curve.
Again we will have no regular complements
by Example~\ref{gcomp}.\ref{mlem4}.
Such a contraction will be to a rational
singularity, because $E$ will not be
in $\lfloor D\rfloor$.

We prove case by case that such an $E$ exists,
except for a case when $K+B\equiv 0$.

Indeed, if we have (WLF) or $\kappa^*=2$,
then we have (RPC), and the latter can
be preserved after a crepant log resolution above.
Take a weighted linear combination of $B$ and $B^+$.
In addition, $-(K+C)$ is not nef when
(ASA) fails.
Thus we have an extremal contraction
$S\to Z$ which is positive with respect to $K+D$,
and $\dim Z\ge 1$.
If $Z$ is a curve, $K+B^+\equiv 0/Z$,
but $K+D$ is numerically positive$/Z$,
which is impossible for $B^+\ge D$ as
in the above case.
Hence we have $E$.
An induction on $\rho$ and contractions
of such $E$ give a contradiction in this case.

Suppose now that $\kappa^*=1$ and we
have a numerical contraction $\nu :S\to Y$
for $K+B$.
By (M) the horizontal multiplicities
satisfy (SM), and $B^+=D=B$ in
the horizontal components.
Thus $-(K+D)$ is nef
on the horizontal curves.
Moreover, it is nef.
Otherwise we have a vertical exceptional curve $E$
with $(K+D.E)>0$.
As above we contract $E$.
Finally, $-(K+D)$ is nef,
and $-(K+D)\equiv 0$ by the Inductive Theorem.
It is possible only when $B^+=D$.
But then we have (ASA) which does not hold
in our case too.

So, we get $\kappa^*=0$ or
$K+B\equiv 0$.

Now let $D$ denote a support of
the non-standard multiplicities in $B$.
If $D$ has a big type, then we get
(WLF) for $K+B-\varepsilon D'$ for
some effective $\mathbb R$-divisor $D'$
with $\Supp{D'}\le D$ by
Proposition~\ref{indcomp}.\ref{types}.
Again we do not have regular complements
by the Monotonicity~\ref{indcomp}.\ref{monl2}:
the non-standard multiplicities $>6/7$
under (M).
We do the same when $D$ has a fibre type.
Finally, $D$ has an exceptional type
and we contract $D$ to points.
(Cf. Proof of Corollary~\ref{int}.\ref{mcorol}:
Numerically trivial case in Section~\ref{int}.)

So, $K+B\equiv 0$, satisfies (SM), but
does not have the regular complements
by Lemma~\ref{gcomp}.\ref{mlem3}
and the Monotonicity~\ref{indcomp}.\ref{monl1}.
Then we have only the trivial complements.
In that case $B^+=B$,
which contradicts 
a Kawamata log terminal property of $K+B$.
This gives a contradiction to our assumption
on an existence of a non-regular and non-exceptional
complement $(S,B^+)$.
\end{proof}

\begin{proof}[Proof of
the Main Theorem:
Global case]
Follows from
Theorem~\ref{gcomp}.\ref{gmainth}.
\end{proof}

Now we slightly improve
Proposition~\ref{indcomp}.\ref{relcond}.

\refstepcounter{subsec}
\label{drel}
\begin{proposition}
Assuming that $K+C+B$ is log canonical
and nef$/Z$,
$$
{\rm (EC)}\Longrightarrow
{\rm (NTC)}\Longleftrightarrow
{\rm (ASA)}
$$
with the exception (EX2) of
Proposition~\ref{indcomp}.\ref{relcond}.
Nonetheless, in (EX2) there exists
a $1$-complement.

Moreover, we can replace (EC) by its weaker form:
\begin{description}
\item{\rm (EC)'}
there exists a boundary $B'$ such that
$K+B'$ is log canonical and $\equiv 0/Z$,
\end{description}
{\rm i.~e., (EC) for $S$\/}.
\end{proposition}

\begin{proof}
Let $(S,B')$ be as in (EC)'.
If we replace $C+B$ by a weighted linear
combination of $C+B$ and $B'$
we can suppose that $K+C+B$ and $K+B'$
have the same log singularities:
\begin{itemize}
\item
the exceptional and non-exceptional divisors with
the log discrepancy $0$, and
\item
the exceptional and non-exceptional divisors with
the log discrepancies $<1$.
\end{itemize}
Note that (EX2) will mean that $B'$ is unique and
$B'\ge B+C$, i.~e., (EEC) holds.

After a log terminal resolution we
suppose that $K+B'$ is log terminal.
By the above properties, a support $D$ of
curves, where $C+B>B'$, is disjoint
divisorially from $\LCS{(S,B')}$.
If $D$ has an exceptional type,
we can contract it, when $K+C+B\equiv 0$ on $D$.
Then $B'\ge B$ and we have (EEC),
which implies the proposition
by Proposition~\ref{indcomp}.\ref{relcond}.
If $K+C+B$ is negative somewhere on $D$,
then (WLF) and (RPC) hold for $K+C+B-\varepsilon D'$
with some $\varepsilon>0$, and $D'$
having $\Supp{D'}\le D$.

On the other hand, if $D$ has a big type,
(WLF) and (RPC) hold for $K+C+B-\varepsilon D'$
with some $\varepsilon>0$, and nef and big $D'$
having $\Supp{D'}\le D$.
Here we may have one exception (EX1), when
$K+C+B$ satisfies (NTC).

In addition, the proposition holds when
$K+C+B\equiv 0$ by its semi-ampleness.

Finally, $D$ has a fibre type and $-(K+C+B)$
has a numerical dimension $1$ and $Z=\pt$
(a global case).
If a fibering given by $D$ does not agree
with $K+C+B$, i.e., $(K+C+B.F)<0$ on
a generic fibre. 
Then $K+C+B-\varepsilon F'$ satisfies (WLF)
for a divisor $F'$ with $\Supp{F'}\le D$, which
defines the fibering.
Otherwise, the fibering gives
a numerical contraction for $K+C+B$.
\end{proof}

\begin{proof}[Proof of
Corollary~\ref{int}.\ref{mcorol}:
Global case]
Follows from
Theorem~\ref{gcomp}.\ref{gmainth}
and Proposition~\ref{gcomp}.\ref{drel}.
\end{proof}

\begin{proof}[Proof of
Corollary~\ref{int}.\ref{nonvan1}]
Again in the global case, it follows from
Theorem~\ref{gcomp}.\ref{gmainth}
and Proposition~\ref{gcomp}.\ref{drel}.

In the local case, use Theorem~\ref{lcomp}.\ref{lmainth}.
\end{proof}

\refstepcounter{subsec}
\label{mrem}
\begin{remcor}
We can improve most of the above results as well.

\ref{gcomp}.\ref{mrem}.1.
In the Main and Inductive Theorems we can replace
(WLF) by (EC)' of
Proposition~\ref{gcomp}.\ref{drel}.

We anticipate that the Main Theorem and
Corollary~\ref{int}.\ref{mcorol}
hold without (M)
as does the Inductive Theorem.
Of course, then exceptional complements
may be unbounded
(cf. Example~\ref{indcomp}.\ref{nonrcomp}).

\ref{gcomp}.\ref{mrem}.2.
By the Monotonicity~\ref{indcomp}.\ref{monl1},
in Corollaries~\ref{int}.\ref{ind2}-\ref{ind3},
we can replace (SM) by
\begin{description}
\item{\rm (M)'}
the multiplicities $b_i$ of $B$
are {\rm standard\/}, i.~e., $b_i=(m-1)/m$ 
for a natural number $m$, or
$b_i\ge I/(I+1)$,
where $I$ is maximal among
the indexes under (SM): $I|I_2$.
\end{description}
(Cf. Classification~\ref{csing}.\ref{singth}.1 below.)
\end{remcor}

\section{Exceptional complements}
\label{ecomp}

In this section we start a classification
of the exceptional complements.
By the Main and Inductive Theorems, they arise
only in the global case $(S,B)$ when
$K+B$ is Kawamata log terminal.
By Remark~\ref{gcomp}.\ref{mrem}.1,
we can assume just (EC)' and (M) as
additional conditions.
In a classification we describe
such $(S,B)$, which will also be called
{\it exceptional\/}, and
their {\it minimal\/} complements.
Here we do this completely in a few cases.
An importance of a complete classification
of the exceptions
will be illustrated in Section~\ref{csing}.
We will continue the classification
elsewhere.

Since the exceptional complements
are bounded, the following invariant
\begin{align*}
\delta (S,B)=\#\{E\mid E &\text{ is
an exceptional or non-exceptional divisor}\\
&\text{ with the log discrepancy }
a(E)\le 1/7\text{ for }K+B\}
\end{align*}
is also bounded.
It is independent on the crepant modifications.

\refstepcounter{subsec}
\label{emainth}
\begin{theorem}
$\delta\le 2$.
\end{theorem}

\begin{statement1}
If $\delta =0$, then $(K,B)$
is $1/7$-log terminal, and
$B$ has only multiplicities
in $\{0,1/2,2/3,3/4,4/5,5/6\}$.
However, the m.l.d. 
of $K+B$ is only $>1/7$.
\end{statement1}

A minimum of such m.l.d.'s exists,
but it is yet unknown explicitly.

In the other cases, $\delta\ge 1$,
the m.l.d. of $K+B$ is $\ge 1/7$,
and after a crepant resolution we assume
that $K+B$ is $1/7$-log terminal in the closed points.
To classify the original $(S,B)$ we need
to find crepant birational contractions of
the $1/7$-log terminal pairs $(S,B)$.
To classify the latter pairs we consider
their minimizations $g:S\to S_{\rm min}$
as in the strategy of the proof
of Theorem~\ref{gcomp}.\ref{gmainth}.
In this section some results on $(S_{\rm min},B_{\rm min})$
and their classification are given.
According to the strategy,
Theorem~\ref{ecomp}.\ref{emainth}
is enough to prove for $(S_{\rm min},B_{\rm min})$.
So, we assume in this section
\begin{itemize}
\item
$\rho (S)=1$,
\item
$K+B$ is $1/7$-log terminal in the closed points,
\item
$B$ has a multiplicity $b_i\ge 6/7$,
\item
$-(K+B)$ is nef, but
\item
$K+D$ is ample for $D=\lfloor(n+1)B\rfloor/n$
with any $n\in RN_2$.
\end{itemize}

To find all such $1/7$-log terminal pairs $(S,B)$
with $\rho(S)>1$ we need to find $K+B'\equiv 0$
with $B'\ge B$ and $\rho(S)=1$.
The former pairs are crepant partial
resolutions of $(S,B')$.
See the strategy and (BPR) in Section~\ref{gcomp}.

Let $C=\lfloor D\rfloor$ denote a support of
the curves $C_i$ with $\mult{C_i}{B}\ge 6/7$,
and $D$ be the same as in
the Minimization of Section~\ref{gcomp}.
Let $F$ be the rest of $B$ or, equivalently, be 
the fractional part of $D$:
$F=\{D\}=\sum b_i D_i$ for $D_i$ with
$b_i=\mult{D_i}{B}\le 5/6$.
By the Inductive Theorem,
$K+D$ is ample for such $D$.

\begin{statement2} For $\delta=1$,
a curve $C$ is irreducible,
has only nodal singularities
and at most one node.
The arithmetic genus of $C$ is $\le 1$.
Divisor $F$ does not pass the node.
\end{statement2}

Abe found a classification in
the {\it elliptic\/} case when
$C$ has arithmetic genus $1$ \cite{Ab}. 

\begin{statement3} For $\delta=2$, $C=C_1+C_2$,
where $C_1$ and $C_2$ are irreducible
curves with only normal crossings
in non-singular points of $S$.
Divisor $F$ does not pass $C_1\cap C_2$
and $b_1+b_2<13/7$, where
$b_i=\mult{C_i}{B}$.
Constant $c$ below is as in 
Lemma~\ref{gcomp}.\ref{unb}

For $C$, we have only the
following three configurations:
\begin{description}
\item {$(\text{I}_2)$}
$C=C_1+C_2$ and $C_i$'s form a wheel; and
\item {$(\text{A}_2)$}
$C=C_1+C_2$, and $C_i$'s form a chain.
\end{description}

Moreover, the curves $C_i$ are non-singular rational 
$m_1\ge m_2\ge 0$-curves, except for 
Case~{\rm ($\rm A_2^6$)} below.

In the case $(\text{A}_2)$, the only possible cases are:
\begin{description}
\item{\rm ($\rm A_2^1$)}
$S=\mathbb P^2$ whereas $C_1$ and $C_2$ are
straight lines, 
$F=\sum d_i D_i$ and $1<\sum d_i\deg D_i\le 3-b_1-b_2$,
assuming that $K+B$ is log terminal. 
\item{\rm ($\rm A_2^2$)}
$S$ is a quadratic cone, whereas
$C_1$ is its section and $C_2$ is its generator,
$2b_1+b_2\le 8/3$.
$F=(2/3)D_1$, where $D_1$ is
another section not passing the vertex. 
$c=1/21$.
\item{\rm ($\rm A_2^3$)}
$S$ is a normal rational cubic cone
whereas $C_1$ is its section and $C_2$
is its generator, $3b_1+b_2\le 7/2$. 
$F=(1/2)D_1$ where
$D_1$ is also a section.
Both sections $C_1$ and $D_1$ 
do not pass the vertex
and $\# C_1\cap D_1\ge 2$. 
$c=1/14$. 
\item{\rm ($\rm A_2^4$)}
$S$ has $B=(6/7)(C_1+C_2)+(1/2)D_1$, $m_1=1$ and $m_2=0$, whereas $S$
has only two singularities
$P_1\in C_1$ and $P_2\in C_2$,
and $D_1$ is a non-singular rational
$1$-curve with a single simple intersection with $C_2$,
a single simple intersection with $C_1$ and
with another single intersection with $C_1$ in $P_1$.
Singularity $P_i$ is
Du Val of type $\mathbb A_i$.
\item{\rm ($\rm A_2^5$)}
$S$ has $B=(6/7)(C_1+C_2)+(1/2)D_1$, $m_1=1$ and $m_2=0$, whereas $S$
has only two singularities
$P_1\in C_1$ and $P_2\in C_2$,
and $D_1$ is a non-singular rational
$1$-curve with a single simple intersection with $C_2$,
a single simple intersection with $C_1$ and
with another single intersection with $C_1$ in $P_1$.
Singularity $P_1$ is simple with $(-3)$-curve
in a minimal resolution,
singularity $P_2$ is Du Val of type $\mathbb A_3$.
$c=1/7$, and $(S,B)$ is $14$-complementary
and the complement is trivial.
\item{\rm ($\rm A_2^6$)}
$S$ has $B=(6/7)(C_1+C_2)$, whereas $S$
has only two singularities
$P_1,P_2\in C_2$.
Curve $C_1$ has arithmetic genus $1$
and has only nodal singularities;
at most $1$.
Curve $C_2$ is a rational non-singular $(-1)$-curve. 
Singularities $P_i$ are Du Val of type $\mathbb A_i$.
$c=1/7$, and $(S,B)$ is $7$-complementary
and the complement is trivial.

In the case $(\text{I}_2)$ the only possible cases are:
\item{\rm ($\rm I_2^1$)}
$S$ is a quadratic cone whereas
$C_1$ and $C_2$ are its two distinct sections,
$b_1+b_2\le 7/4$.
$F=(1/2)L$,
where $L$ is a generator of cone $S$.
$c=3/28$.
\item{\rm ($\rm I_2^2$)}
$S$ has $B=(6/7)(C_1+C_2)$, $m_1=1,m_2=2$, whereas $S$
has only two singularities $P_i\in C_i$. 
Singularities $P_i$ are Du Val of type $\mathbb A_i$.
$c=1/7$, and $(S,B)$ is $7$-complementary
and the complement is trivial.
\end{description}
\end{statement3}

\refstepcounter{subsec}
\label{logter}
\begin{proposition}
Under the assumptions of this section,
$K+D$ is formally log terminal, except for the case
when $P$ is non-singular and near $P$, $D=C+(1/2)C'$
with non-singular irreducible curves $C$ and $C'$
having a simple tangency; $\mult{C}{B}< 13/14$.
\end{proposition}

Note that the latter log singularity appears
only on cones: Cases (A$_2^{1-3}$) when $D_1$
is tangent to $C_1$.

\begin{proof} By the proof of
the Minimization in Section~\ref{gcomp},
$K+D$ is log canonical,
and ($1/7$-)log terminal outside $C$.
So, we need to check a log terminal property
formally (locally in an analytic topology)
in the points $P\in C$.

First, suppose that $C\not= D$ in
a neighborhood (even Zariski) of $P$.
Then $K+C$ is purely log terminal
and $C$ is non-singular by
\cite[Lemma~3.6]{Sh2}.
But $S$ may have a singularity of index $m$
in $P$ \cite[Proposition~3.9]{Sh2}.
If we have formally
two distinct prime divisors
(two branches) $D_1$ and $D_2$
through $P$ in $\Supp (D-C)$,
then, by~\ref{indcomp}.\ref{invadj}.1,
$K+D$ will be log canonical in $C$ only when
$b_1=b_2=1/2$ and,
in a neighborhood of $P$, $K+D=K+C+(1/2)(D_1+D_2)$.
(We would like to remind the reader that all non-reduced
$b_i=(n-1)/n$ with $n=1,2,3,4,5$ or $6$.)
So, if $P$ is singular,
then by a classification of surface
log canonical singularities,
the curves $E$ of a minimal resolution
form a chain, whereas a birational transform
of $C$ intersects simply one end of the chain,
and that of $D_i$ intersects simply another end.
The intersection points are outside
the intersections of the curves $E$.
Thus the log discrepancy $a=(E,K+B)$,
in any $E$ for $K+B$, is the same
as $a(E,K+b C+D_1)$ for $K+b C+D_1$,
where $b=\mult{C}{B}\ge 6/7$ near $P$.
Therefore, by Corollary~\ref{indcomp}.\ref{adjunct},
$a=a(E,K+B)=1-\mult{P}{(b C)_{D_1}}
\le 1-\mult{P}{((6/7)C)_{D_1}}\le 1-6/7=1/7$
for an exceptional divisor $E$.
This is impossible by our assumptions.
Hence $P$ is non-singular.
Then the monoidal transform in $P$ gives $E$
with the same property.

Hence we may have formally at most one irreducible
component (or a single branch)
$C'$ of $\Supp (D-C)$ through $P$.
By the form of it and a classification 
of log canonical singularities \cite{K1}, $K+D$ will
be log terminal in $P$, except for the case
when $C'$ has multiplicity $1/2$ in $D$, and $B$
intersects only an end curve in a minimal resolution
of $P$, or $D=C$ near $P$.
Again, as in the above case,
we have a contradiction with assumptions, except
for the case when $P$ is non-singular and $C'$
has a simple tangency with $C$ in $P$.
Such singularity will be $1/7$-log terminal
for $K+B=b C+(1/2)C'$, only when
$b<13/14$.

Finally, we suppose that $C=D$ near $P$,
then by \cite[Theorem~9.6~(6)]{K1} $P$ has type $\mathbb D_m$ with
$m\ge 3$ when $C$ is formally irreducible in $P$.
Then it has an exceptional divisor $E$ on a minimal resolution
with $a(E,K+B)\le 1-b\le 1/7$, which is again impossible.
Indeed, $a(E,K+C)=0$, $a(E,K)\le 1$ and
$a(E,K+B)=a(E,K+b C)=a(E,K+C)+(1-b)\mult{E}{C}=
(1-b)a(E,K)\le 1/7$.

Otherwise by \cite[Theorem~9.6~(7)]{K1}, $C$ has two branches
in $P$ and again, by the $1/7$ log terminal property in $P$,
$P$ will be non-singular.
Hence $K+B=K+C$ is formally log terminal here.
\end{proof}

\refstepcounter{subsec}
\label{node}
\begin{corollary}
$C$ has only nodal singularities, and
only in non-singular points of $S$.
$C$ is connected.
Moreover, each irreducible component of $C$ intersects
all other such components.
\end{corollary}

\begin{proof}
The first statement follows from
Proposition~\ref{ecomp}.\ref{logter}.
Since $\rho=1$, each curve on $S$ is ample,
which proves the rest.
\end{proof}

Let $g$ be the arithmetic genus of $C$.
\refstepcounter{subsec}
\label{genus}
\begin{proposition}
$g\le 1$.
\end{proposition}

In its proof and in a proof of
Theorem~\ref{ecomp}.\ref{emainth} below,
we use the following construction.

We reconstruct $S$ into a non-singular
minimal rational model $S'$.
Make a minimal resolution
$S^{\rm min}\to S$ of $S$, and then
contract $(-1)$-curves: $S^{\rm min}\to S'$,
where $S'$ is minimal.
Then $S'=\mathbb P^2$ or $S'=\mathbb F_m$,
because $S$ and $S'$ are rational.

By Corollary~\ref{ecomp}.\ref{node},
the resolution $S^{\rm min}\to S$
preserves $C$ up to an isomorphism.
A birational transform of a curve $C_i$ or another one
on $S'$ are denoted again by $C_i$ or as the other one
respectively.
According to the LMMP,
$(S',B')$ is log canonical and
$-(K_{S'}+B')$ is nef, because
the same holds for $K+B$ and its crepant
blow-up $K^{\rm min}+B^{\rm min}$,
where $K^{\rm min}=K^{S^{\rm min}}$
and $B'$ is the image of $B^{\rm min}$.
An image of $B^{\rm min}$ is not less
than a birational image of $B$.
So, $-(K_{S'}+B)$ is nef for $S'=\mathbb P^2$.

\begin{statement1}
Moreover, on a minimal rational model $S'$,
$g(C)\le 1$ and that $g(C')=0$ for
each (proper) $C'\subset C$, except for
Case~($\rm A_2^6$) in \ref{ecomp}.\ref{emainth}.3.
\end{statement1}

\begin{proof}
Suppose that $g\ge 2$.

$S^{\rm min}\to S'$ cannot contract all $C$,
because it has to be a tree of
non-singular rational curves.
However, we may increase $g$ after contraction
of some components of $C$ and
other curves on $S^{\rm min}$.

If $S'=\mathbb P^2$,
then $-(K_{S'}+B)$ and
$-(K_{S'}+(6/7)C)$ are nef, and
$\deg C\le 3$,
which is only possible for $g\le 1$.

Therefore $S'=\mathbb F_m$ with $m\ge 2$.
Indeed, original $S'\not\cong \mathbb F_0$,
because $\rho=1$.
So, if final $S'=\mathbb F_0$,
we had before a contraction of $(-1)$-curve.
Then we can reconstruct $S$ into $S'=\mathbb F_1$,
and so into $\mathbb P^2$.

By Corollary~\ref{ecomp}.\ref{node},
we have at most one fibre $F$
of $\mathbb F_m$ in $C$.
If a unique {\it negative\/} section $\Sigma$
is not in $C$,
then $\sigma=\mult{\Sigma}{B'}\le
2 -2\times (6/7)=2/7<1/3$,
because $C$ is not a section of $\mathbb F_m$ over
the generic point of $\Sigma$.
(Otherwise $C$ will be rational with
only double singularities, and $g=0$.)
By the nef property of $-(K+B)$,
we have inequality
$0\ge (K_{S'}+B'.\Sigma)\ge
(K_{S'}+(2/7)\Sigma.\Sigma)$.
This implies that $\Sigma$ is $(-2)$-curve,
and $m=2$.
If we had before a contraction of a $(-1)$-curve we
can reconstruct $S'$ to the above $S=\mathbb P^2$.
Hence $(S'=S^{\rm min},B'=B^{\rm min})$ is
the minimal resolution of $S$, and
$C\cap \Sigma=\emptyset$
(cf. Lemma~\ref{ecomp}.\ref{ineq} below).
Since $\rho(S)=1$, $S$ is a quadratic cone
(or a quadric of rank $3$ in $\mathbb P^3$)
with a double section $C=D$ ($\sim -K$)
not passing through its vertex.
So, by the Adjunction $g=1$.

Finally, $\Sigma$ is a component of $C$, then we
have another component $\Sigma '$ in $C$,
which is also a section$/\Sigma$.
If $C=\Sigma+\Sigma'$ has $g\ge 2$,
then $(\Sigma.\Sigma')\ge 3$.
This is impossible, because
$0\ge (K_{S'}+B'.\Sigma)\ge
(K_{S'}+\Sigma+(6/7)\Sigma'.\Sigma)=
\deg (K_{\Sigma}+(6/7)(\Sigma'|_{\Sigma}))
\ge -2+3(6/7)=4/7$.
One last case $C=\Sigma+\Sigma'+F$, whereas
$C$ has $g\ge 2$, only when $(\Sigma'+F.\Sigma)\ge 3$,
because $\Sigma\cap\Sigma'\cap F=\emptyset$
by the log canonical property: $3\times (6/7)>2$.
Then we may act as above replacing $\Sigma'+F$ by
$\Sigma'$.

Now we prove \ref{ecomp}.\ref{genus}.1.
If $C$ has a component, say $C_1$,
of the arithmetic genus $g\ge 1$.
Then according to the above,
$g=1$, $C=C_1+C_2$, where $C_2$
is non-singular rational ($m_2$)-curve,
and $C_1$ intersects $C_2$ in one point.
Moreover, if $m_2\ge 0$, then $C_2$ is not
exceptional on $S'$ and $C$ in $S'$ has genus $\ge 2$,
which is impossible as we know.
On the other hand, $(K_X.C_2)<0$
since $\rho (S)=1$.
Therefore $m_2=-1$.
However, this is only possible for Case~(A$_2^6$)
by the following lemma. 
\end{proof}

\refstepcounter{subsec}
\label{move}
\begin{lemma}
Each non-singular irreducible rational
(proper) component $C_i\subset C$
is movable on a minimal resolution of $S$, i.e.,
$C_i$ is an $m$-curve with $m\ge 0$, except for 
Case~{\rm ($\rm A_2^6$)} in \ref{ecomp}.\ref{emainth}.3.
\end{lemma}

\refstepcounter{subsec}
\label{ineq} 
\begin{lemma}
Let $P$ be a log singularity
$(S,B)$ such that
\begin{description}
\item{\rm (i)}
$B\ge b C>0$, where $C$ is an
irreducible curve through $P$;
\item{\rm (ii)} 
$B$ is a boundary, and
\item{\rm (iii)}
$P$ is a singularity of $S$.
\end{description}
Let $E$ be a curve on a minimal resolution of $P$
intersecting the proper inverse image of $C$,
and let $d=1-a(E,K+B)$ be the multiplicity of
the boundary on the resolution for the crepant
pull back.

Then $d\ge ((m-1)/m)b$ where $m=-E^2$.

Moreover, $d\ge (1/2)b$ always, and $=$ holds
only when $P$ is simple Du Val, $B=b C$ near $P$,
and $K+C$ is log terminal in $P$.

Otherwise $d\ge (2/3)b$.

In addition, $d\ge ((m-1)/m)b$ whenever
$P$ is log terminal for $K+C$ of index $m$.
In this case we may also include in $C$
components of $B$
with standard multiplicities
{\rm (cf. Lemma \ref{indcomp}.\ref{invadj})}.
\end{lemma}

\begin{proof} Let $C$ be the inverse image of
$C$ in $S^{\rm min}$.
Since the resolution is minimal all
multiplicities of $B^{\rm min}$ are non-negative.
So, we may consider a contraction of only $E$, and
suppose that $B=b C$.
Then we may find $d$ from the following
equation $(K+b C+d E.E)=0$.
If $E$ is singular or non-rational, $d\ge 1\ge (1/2)b$
and even
$\ge (2/3)b$ because $b\le 1$ by (ii).

Otherwise $E$ is $(-m)$-curve with $m\ge 2$.
Hence $d=(m-2)/m+(C.E)(1/m)b\ge ((m-2)/m+ (1/m))b=
((m-1)/m)b\ge (1/2)b$, because $0\le b\le 1$.
Moreover $=$ only for $m=2$, $B=b C$ near $P$,
and $(C.E)=1$ in $P$ when $b>0$.
The next calculation shows that $P$ is
a simple Du Val singularity when $b>0$
and $d= (1/2)b$.

If $m\ge 3$, then $d\ge (2/3)b$.
The same holds if we replace $E$ by a pair of
intersected $(-2)$-curves.

Finally, $d(b)=1-a(E,K+B)$ is
a linear function of $b$.
So, it is enough to check the last inequality
for $b=0$ and $1$.
For $b=0$, $a\le 1$ and $d\ge 0$ by (iii).
This gives the required inequality.
For $b=1$, $a(E,K+B)\le a(E,K+C)=1/m$
by \cite[3.9.1]{Sh2}.
Hence $d\ge (m-1)/m\ge ((m-1)/m)b$.

We may also include in $C$
components with standard coefficients
by Lemma~\ref{indcomp}.\ref{invadj}.
\end{proof}

\begin{proof}[Proof of Lemma~\ref{ecomp}.\ref{move}]
Since $-K$ is ample on $S$ we should only eliminate
the case when $C_i$ is a $(-1)$-curve.

Since $C^2_i>0$ on $S$ and by
Proposition~\ref{ecomp}.\ref{logter},
$C_i$ has at least two singularities $P_1$ and $P_2$.
They are distinct from the intersection
points $P=(C\setminus C_i)\cap C_i$.
Such an intersection point $P$
exists because $C_i\not= C$ and
by Corollary~\ref{ecomp}.\ref{node}.

So, we may calculate $(C_i.K+B)$ on a minimal
resolution $S^{\rm min}/S$.
We denote again by $C'$ and $C_i$ respectively
a proper inverse image of $C'=C\setminus C_i$
and $C_i$.
Over $P_1$ and $P_2$ we have, respectively, single
(non-singular rational)
curves $E_1$ and $E_2$
intersecting $C_i$ in $S^{\rm min}$.
Let $b$, $d$, $b_1\le b_2$ be the multiplicities of
$B^{\rm min}$ in $C'$ (in any component through $P$),
$C_i,E_1,E_2$ respectively.
Then by our assumptions $b,d\ge 6/7$.
On the other hand $b_1\ge (1/2)d$ because $P_1$
is singular, and $b_1=d/2$ is attained only
for a simple Du Val singularity in $P_1$.
Since $C^2_i>0$, then $P_2$ is not such
a singularity and $b_2\ge (2/3)d$ by
Lemma~\ref{ecomp}.\ref{ineq}.
Otherwise we have a third singularity $P_3$
of $S$ in $C_i$, and
$b_3=\mult{E_3}{B^{\rm min}}\ge (1/2)d$,
where $E_3$ intersects $C_i$ in $S^{\rm min}$.

For three points and more:
$(K+B.C_i)\ge -1+b-d+b_1+b_2+b_3\ge
-1+b+d/2\ge -1+6/7+3/7=2/7>0$.
Therefore we have only two singularities  
and $0\ge (K+B.C_i)=
(K^{\rm min}+B^{\rm min}.C_i)\ge
(K^{\rm min}+b C'+d C_i+b_1 E_1+b_2 E_2.C_i)\ge
-1+b-d+b_1+b_2\ge
-1+b-d+(1/2)d+(2/3)d=-1+b+(1/6)d\ge
-1+6/7+1/7=0$.
Hence we have the equality,
whereas $C_i=C_2$ and $C'=C_1$,
$B=(6/7)(C_1+C_2)$, $P_1$ is a simple Du Val singularity,
and $P_2$ is a Du Val singularity of type $\mathbb A_2$.
Otherwise it is a simple singularity 
with $(-3)$-exceptional curve, because $b_2=(2/3)d$
and $C_2^2>0$.
However, this is impossible for $m_2=-1$.

So, we can resolve 
singularities $P_1$ and $P_2$ by $(-2)$-curves
$E_1$ and $E_2,E_3$ respectively,
where $E_2$ intersects $C_2$ on the resolution.
If we contract successively $C_2,E_2$ and $E_3$,
we transform $E_1$ into a ($1$)-curve
which is tangent to the transform of $C_1$
with order $3$.
Hence $S'=\mathbb P^2$ and $C_1$ is
a cubic in it.
Finally, $K+B\equiv 0$ and we have
a trivial $7$-complement
according to the computation. 
\end{proof}

\begin{proof}[Proof of
Theorem~\ref{ecomp}.\ref{emainth}]
$C$ has at most
$\delta=3$ irreducible components $C_i$ by
Proposition~\ref{ecomp}.\ref{genus}.
Otherwise $g\ge 2$ and $2$ is attained
when $C$ has four non-singular rational components
with one intersection point for each pair of
components by
Corollary~\ref{ecomp}.\ref{node}.

Thus we prove~\ref{ecomp}.\ref{emainth}.1-2,
and can assume that $C$ has at least $\delta\ge 2$
components $C_i,1\le i\le \delta$.

The $1/7$-log terminal property of $K+B$ implies
that $F$ does not pass $C_1\cap C_2$
and $b_1+b_2<13/7$.

By Proposition~\ref{ecomp}.\ref{logter},
$F$ does not pass the nodes of $C$:
$2\times (6/7)+1/2>2$.
In particular, $F$ does not pass
$C_i\cap C_j$ for $i\not= j$.

Except for Case~(A$_2^6$),
each component $C_i$ is non-singular and
rational.
By~\ref{ecomp}.\ref{genus}.1 this holds on $S'$.
Then, for $S$, it is implied by
Lemma~\ref{ecomp}.\ref{move}.

So, excluding Case~(A$_2^6$) in what follows, 
we have only the
following three configurations of 
non-singular rational curves $C_i$:
\begin{description}
\item {$(\text{I}_3)$}
$C=C_1+C_2+C_3$ and $C_i$'s form a wheel;
\item {$(\text{I}_2)$}
$C=C_1+C_2$ and $C_i$'s form a wheel; and
\item {$(\text{A}_2)$}
$C=C_1+C_2$, and $C_i$'s form a chain.
\end{description}
This follows from
Corollary~\ref{ecomp}.\ref{node} and
Proposition~\ref{ecomp}.\ref{genus}.

To eliminate some of these cases we
prove that each $C_i$, with two nodes
$P_1$ and $P_2$ in $C$, is an $m_i$-curve
with $m_i\ge 1$.
We know or can suppose that $C_i$ is non-singular rational,
or an $m_i$-curve.
Moreover, by Lemma~\ref{ecomp}.\ref{move},
$m_i\ge 0$.
We suppose that $m_i=0$, and
derive a contradiction.

Since $C^2_i>0$ on $S$, then
$C_i$ has at least one singularity $P$ of $S$,
$P\not=P_1$ and $P_2$.

As above we may calculate $(K+B.C_i)$ on a minimal
resolution $S^{\rm min}/S$.
We denote by $C'$, $C''$ branches of
$C\setminus C_i$ in $P_1$ and $P_2$
respectively.
We identify them with their proper inverse
images on $S^{\rm min}$.
Over $P$ in $S^{\rm min}$,
we have, respectively, one
(non-singular rational) curve $E$ intersecting $C_i$.
Let $b'$, $b''$, $d$, $b$ be the multiplicities
of $B^{\rm min}$ in $C',C'',C_i,E$.
Then by our assumptions $b',b'',d\ge 6/7$.
On the other hand, $b\ge (1/2)d$
by Lemma~\ref{ecomp}.\ref{ineq}.

Therefore we get a contradiction: $0\ge (K+B.C_i)=
(K^{\rm min}+B^{\rm min}.C_i)\ge
(K^{\rm min}+b' C'+b'' C''+d C_i+b E.C_i)\ge
-2+b'+b''+b\ge
-2+b'+b''+(1/2)d\ge
-2+6/7+6/7+3/7=1/7$.

Now we are ready to verify that
$g(C)=0$ and we
have Case~($\rm A_2$), except for
two Cases ($\rm I_2^1$) and (I$_2^2$) 
in \ref{ecomp}.\ref{emainth}.3
with configuration (I$_2$).
We want to  eliminate Case~($\text{I}_3$) and
the other cases in ($\text{I}_2$).
According to that which was proved above
and to the construction, $S'$
again has the same curves $C_i$
as components of $C$:
$m_i'$-curves with $m_i'\ge m_i\ge 1$.

First, we consider $S'=\mathbb P^2$.

Since $4\times (6/7)>3$, they are all
$1$-curves in Case~($\text{I}_3$).
Hence there are no contractions of
$(-1)$-curves onto $C\subset S'$
for $S^{\rm min}\to S'$.
In particular, we preserve curve $E$
over any singularity $P\in C$ of $S$,
or a curve $E$ with a standard multiplicity
$0<\mult{E}{B}<1$.
Either has multiplicity $b\ge 3/7$.
So, we get 
$0\ge \deg(K_{S'}+b_1 C_1+b_2 C_2+b_3 C_3+b E)\ge
-3+b_1+b_2+b_3+b\ge -3+3\times (6/7)+3/7=0$.
Hence $b=3/7$ and at most one of curves $C_i$,
say $C_1$, has a singularity.
This is impossible when $S\not=S'$ because
then $S$ is a rational cone by
the condition $\rho (S)=1$, and
because then $S^{\rm min}=\mathbb F_1$,
when $S'=\mathbb P^2$.
So, $S=S'=\mathbb P^2$ and $F=0$.
But then $(S,B)$ is $1$-complementary
with $B^+=C$, which contradicts 
our assumptions.

We can do the same in Case~($\text{I}_2$),
when the ($1$)-curve $C_i$ does not have
singularities of $S$,
because then $S=S'=\mathbb P^2$ with $F=0$ and
$1$-complement $B^+=C$.
The case when $C_1$ and $C_2$ are
both ($1$)-curves on $S'$ is only possible when
$g(C)=0$.
So, we have a ($1$)-curve, say $C_1$, on $S$ and on $S'$ with
a single simple Du Val singularity $P_1$ of $S$.
It is only possible in Case~(I$_2^2$).
More precisely, $m_1=m_1'=1,m_2=2,m_2'=4,
B'=(6/7)(C_1+C_2)+(3/7)E_1$ and 
the line $E_1$ is tangent to the conic $C_2$
in a point $P\not\in C_1\cap C_2$ on $S'$.
The inverse transform $S'-\to S$ can be done
as follows.
Surface $S^{\rm min}$ is obtained by
successive monoidal transforms:
first, in $P$ which gives the ($-1$)-curve $E_2$,
then in $E_1\cap E_2$,
which gives the ($-1$)-curve $E_3$,
then in $E_1\cap E_3$,
which gives the ($-1$)-curve $E_4$.
Curves $E_1,E_2$ and $E_3$ are ($-2$)-curves
on $S^{\rm min}$ and
$B^{\rm min}=(6/7)(C_1+C_2)+(3/7)E_1+(2/7)E_2+(4/7)E_3$.
To obtain $S$ we contract $E_1$ to $P_1$
and $E_2,E_3$ to $P_2$.

Now we suppose that $S'=\mathbb F_{m}$ with $m\ge 2$
but never $=\mathbb P^2$.
Since $C_i$ are $m_i$-curves with $m_i\ge 1$,
they are sections of $\mathbb F_{m}/\Sigma$,
$C_i\not=\Sigma$, and
only in Case~($\text{I}_2^1$).

Indeed, as in the proof of
Proposition~\ref{ecomp}.\ref{genus}
$\sigma\le 2/7$,
and $\Sigma$ is an exceptional $(-2)$-curve 
in $S$.
Hence $S^{\rm min}=S'=\mathbb F_2$ and 
$S$ is a quadric of rank 3,
having just one singularity.
Moreover, since $\sigma\le 2/7$, this 
is only possible in Case~($\rm I_2^1$)
with two conic sections (not through the singularity)
$C_1$ and $C_2$.
However $\deg B\le 4$ with respect to
$C_1\sim C_2$.
Therefore, $F=(1/2)L$,
where $L$ is a generator of the quadric,
$c=3/28$, and
$b_1+b_2\le 7/4$.
Since $\mathbb F_0\not=S^{\rm min}$, we reduce
$S'=\mathbb F_0$ to one of the above cases.

In particular, we have proved that $\delta\le 2$.
In our assumptions $\delta=2$.
We know also that $C_1$ and $C_2$ are
$m_1$ and $m_2$-curves, and say $m_1\ge m_2\ge 0$.
Excluding Case~($\rm I_2$) in what follows,
we suppose ($\rm A_2$): $\# C_1\cap C_2=1$.

First, we consider cases with $S'=\mathbb F_m$,
but never $=\mathbb P^2$,
in particular, $m\ge 2$.

As above $S'=\mathbb F_m$
is possible only when $C_1$ is a section and $C_2$
is a fibre of $\mathbb F_m/\Sigma$.
If both $C_1$ and $C_2$ are sections of $S'=\mathbb F_m$
and $m\ge 2$,
then $m=2,S^{\rm min}=S'=\mathbb F_2$, and 
$\Sigma\not=C_1$ and $C_2$.
So, $(C_1.C_2)\ge 2$ which is
impossible under (A$_2$).
By the same reason, $C_1$ is a section
but not a multi-section.

Since $C_2\subset S'$ is a $0$-curve, 
$m_2=m_2'=0$, and
there were no contractions on $C_2$.
So, by Proposition~\ref{ecomp}.\ref{logter},
\begin{description}
\item
$K+D$ is log terminal near $C_2$ on $S$.
\end{description}
Otherwise $F=\{D\}\ge (1/2)D_1$ where curve $D_1$ is
tangent to $C_2$.
Since $(K+D.C_2)>0$, we have on $C_2$ 
a singular point of $S$ or
one more (non-tangency)
intersection point with $F$.
This gives one more curve $D_2$ on $S'$ with
$d_2=\mult{D_2}{B'}\ge (1/2)b_2\ge 3/7$
by Lemma~\ref{ecomp}.\ref{ineq}.
But this is impossible:
$0\ge (K+B.C_2)\ge -2+b_1+2(1/2)+d_2\ge
-2+6/7+1+3/7=2/7$.

Since $K+D$ is log terminal near $C_2$ on $S$
and  $(K+D.C_2)>0$,
we have at least in total
two singularities of $S$ on $C_2$ or
intersection points $P_1$ and $P_2$ with $F$.
Moreover, one of them is not
\begin{description}
\item{\rm ($\Atd{1}{}$)}
a simple Du Val singularity of $S$,
near the singularity $F=0$; nor
\item{\rm ($\Atd{1}{*}$)}
a simple (in a non-singular
point of $S$) intersection point with
a component $D_i$ of $F$ with
$\mult{D_i}{F}=1/2$, near the point
$F=(1/2)D_i$.
\end{description}
Otherwise we have in total three
singularities or intersection points
with $F$, which is impossible,
because then $0\ge (K_{S'}+B'.C_2)\ge 
-2+6/7+3\times (3/7)=1/7>0$
by Lemma~\ref{ecomp}.\ref{ineq}.

We assume that $P_2$ is not ($\Atd{1}{}$) nor
($\Atd{1}{*}$).
So, if $P_2$ is a non-singular point of $S$,
then $F$ has a component $D_2$
passing through $P_2$ with
$\mult{D_2}{F}=(i_2-1)/i_2$ and $i_2\ge 3$.
Moreover, near $P_2$, $D_2$ has
a simple intersection with $C_2$
and $F=((i_2-1)/i_2)D_2$.
In addition, $P_1$ has type ($\Atd{1}{}$) or 
($\Atd{1}{*}$).
Otherwise some $S'=\mathbb P^2$.
Indeed, 
$0\ge (K+B.C_2)\ge -2+6/7+2\times (4/7)=0$.
Then $K+B,K^{\rm min}+B^{\rm min}$ and $K_{S'}+B'\equiv 0$,
$b_1=b_2=6/7$, and the modifications are
crepant.
By Lemma~\ref{ecomp}.\ref{ineq},
$K+D$ has the index $3$ in each $P_i$.
Moreover, $B'=(6/7)(C_1+C_2)+(4/7)(E_1+E_2)$
for divisors $E_1, E_2$ on $S^{\rm min}$ and $S'$
over $P_1$ and $P_2$.
Hence each $P_i$ is
a singularity and $\rho(S^{\rm min})\ge 3$.
So, we may suppose that $S'=\mathbb F_m$
and $m\ge 3$.
Then $\Sigma=E_i$ for some $E_i$.
This gives a contradiction  
$4/7=\sigma\ge 1/3+(1/3)b_2=1/3+2/7$.

Assuming that  
\begin{description}
\item
$P_1$ has type ($\Atd{1}{}$) 
and this is the only singularity of $S$ on $C_2$,
\end{description}
we verify then that $(S,B)$ has type (A$_2^2$).
Indeed, a fractional component $F$ of $D$, with
multiplicity $(l-1)/l, l\ge 3$, intersects $C_2$.
Since $\rho(S)=1$, $S$ is a quadric cone, 
$S^{\rm min}=S'=\mathbb F_2$, and
$0\ge (K+B.C_2)\ge -2+b_1+(1/2)b_2+(l-1)/l\ge
-2+6/7+3/7+(l-1)/l$.
Hence $(l-1)/l\le 5/7$, and $l=3$
which gives Case~(A$_2^2$).
In particular, $F=(2/3)D_1$ for
a section $D_1$ not passing $P_1$.

The next case when
\begin{description}
\item
$P_1$ has type ($\Atd{1}{}$) 
and $P_2$ is a singularity of $S$,
\end{description}
is reduced to $\mathbb P^2$. 
Indeed, $\rho(S^{\rm min})\ge 3$,
and we can suppose that $m\ge 3$. 
Then $\Sigma\not=E_1$, where
$E_1$ is the exceptional curve$/P_1$
on $S^{\rm min}$ or $S'$, because
$E_1$ will be a section of $S'=\mathbb F_m$ with
$E_1^2\ge -2$ and even $\ge m$ on $S'$.
On the other hand, 
$0\ge (K+B.C_2)\ge -2+b_1+\sigma+(1/2)b_2\ge
-2+6/7+\sigma+3/7$ and $\sigma\le 5/7$. 
But, since $(K+B'.\Sigma)\le 0$, 
$5/7\ge \sigma\ge
(m-2)/m+b_2/m\ge
(m-2)/m+(1/m)(6/7)$, which gives
$m\le  4$.
As above after a modification we assume
that $m\le 3$.
So, $m=3$.
Curves $C_1$ and $E_1$ do not
intersect $\Sigma$ simultaneously.
Otherwise we have a contradiction:
$5/7\ge \sigma\ge
1/3+(b_1+b_2+d_1)/3\ge
1/3+5/7$, where $d_1=\mult{E_1}{B'}=(1/2)b_2$. 
Therefore $(E_1.C_1)\ge 3$ and 
the intersection points $E_1\cap C_1$ 
are outside of $\Sigma$. 
We need
to make at least $3$ blow-ups in $E_1$
to disjoint $E_1$ and $C_1$ on $S^{\rm min}$.
Hence we can get $S'=\mathbb P^2$. 

In the next case 
\begin{description}
\item
$P_1$ has type ($\Atd{1}{*}$).
\end{description}
Then $P_2$ is a singularity of $S$,
because $C_2^2>0$ and $m_2=m_2'=0$.
We suppose that $P_2$ has index $l\ge 3$
for $K+D$.
Then we get type (A$_2^3$).
We have no other singularities on $C_2$,
equivalently, $B=b_2 C_2$ and $S$
is non-singular near each other point. 
So, $\Sigma$ is a curve in a resolution
of $P_2$ intersecting $C_2$.
It is an ($-m$)-curve on $S'=\mathbb F_m$
with $m\ge 2$.
Singularity $P_1$ gives
a fractional component: $F\ge (1/2)D_1$,
where $D_1$ is a section of $S'=\mathbb F_m$.
As above, $0\ge(K+B.C_2)\ge -2+b_1+\sigma+1/2\ge
-2+6/7+\sigma+1/2$, and $\sigma\le 9/14$.
By Lemma~\ref{ecomp}.\ref{ineq},
$9/14\ge \sigma\ge ((l-1)/l)b_2\ge ((l-1)/l)6/7$.
So, $l\le 4$.
Moreover, for $l=4$ we have equations:
$b_2=6/7$ and $\sigma=9/14$.
As in the last part of the proof
of Lemma~\ref{ecomp}.\ref{ineq},
this is possible only when $F=0$ near $P_2$
and $P_2$ is a Du Val singularity of type $\mathbb A_3$.
But then $m=2$ and we can reconstruct $S'$
into $\mathbb P^2$.

Therefore $l=3$.
Moreover, according to the same reasons
this is not a Du Val singularity of type $\mathbb A_2$.
Thus $P_2$ is a simple singularity which
can be resolved by ($-3$)-curve.
So, $m=l=3$ and $S$ is a cubic cone.
This is Case~(A$_2^3$).
Moreover, $F=(1/2)D_1$, and both sections
$C_1$ and $D_1$ do not pass the vertex.
Since $K+B$ is log terminal $\# C_1\cap D_1\ge 2$.
Since $K+B$ has a non-positive degree,
$3b_1+b_2\le 7/2$ and $c=1/14$. 

Finally, we consider the case when some
$S'=\mathbb P^2$.
We suppose that $C_1$ and $C_2$
are respectively $m_1$- and $m_2$-curves
on $S$ with $m_1\ge m_2\ge 0$, and
$m_1'$- and $m_2'$-curves on $S'=\mathbb P^2$
with $m_1', m_2'\ge 1$, and 
$m_1'\ge m_2'$ whenever $m_1=m_2$.

Note that $K+D$ is log terminal in this case.
By Proposition~\ref{ecomp}.\ref{logter},
if $K+D$ is not log terminal, then
$F\ge (1/2)D_1$ where curve $D_1$ is
tangent to $C$, as it is
for $C$ and $D_1$ in $S'=\mathbb P^2$.
If $D_1$ is a line in $\mathbb P^2$,
then $3\times (6/7)+1/2>3$ and $K+B'$ is ample.
Therefore $m_1'=m_2'=1$ and
$D_1$ is a conic in $\mathbb P^2$.
It was checked above that $D_1$ on $S$ is
not tangent to the $0$-curves $C_i$.
In other words, $D_1$ is tangent to $C_i$
with $m_i=m_i'=1$.
In addition, we have on $C_i$ a singular point of $S$ or
one more (non-tangency) 
intersection point with $F$.
This gives one more curve $D_2$ with
$d_2=\mult{D_2}{B'}\ge (1/2)b_i\ge 3/7$
by Lemma~\ref{ecomp}.\ref{ineq}.
But this is impossible:
$0\ge \deg(K_{S'}+B')\ge 
\deg(K_{S'}+b_1 C_1+b_2 C_2+(1/2)D_1+d_2 D_2)\ge
-3+b_1+b_2+2(1/2)+d_2\ge -3+2\times (6/7)+1+3/7=1/7$.

If $m_1=m_2=0$, we have contraction
$S^{\rm min}\to S'=\mathbb F_0=
\mathbb P^1\times\mathbb P^1$
given by the linear system $|C_1+C_2|$
on $S^{\rm min}$.
Such cases we consider later.
In other cases $m_1\ge 1$.
We verify that the latter is possible only
for types (A$_2^1$), (A$_2^4$) or (A$_2^5$).
So, first we check that $m_1=1$.

Otherwise, $m_1'=4\ge m_1\ge 2$ and $C_1$ is
a conic on $S'=\mathbb P^2$; 
$m_2'=1$ and $C_2$ is
a line on $S'=\mathbb P^2$,
because $4\times (6/7)>3$.
Since $K+B'$ is log terminal, $C_1$ and $C_2$
have two intersection points $R_1$ and $R_2$.
Moreover, $m_2=0$, because $C_1+C_2$ has
configuration (A$_2$) on $S$ by our assumptions.

As we know, in total $C_2$ has two
singularities of $S$ or intersections with
$F$.
On the other hand,
we have contractions of curves onto $C_2$
for $S^{\rm min}\to S'$ only over one of
points $R_i$.
Thus there exists point $P_1\not= R_1$ and $R_2$ 
on $C_2$ which is singular on $S$ or belongs to $F$.
Moreover,
\begin{description}
\item
$P_1$ on $S$ has type ($\Atd{1}{}$) or 
($\Atd{1}{*})$.
\end{description}
Otherwise, by Lemma~\ref{ecomp}.\ref{ineq},
$P_1$ gives on $S'=\mathbb P^2$ 
a curve $D_1\not= C_1$ and $C_2$ with
$d_1=\mult{E_1}{B'}\ge (2/3)b_2\ge 4/7$.
Then $K+B'$ is ample, because
$3\times (6/7)+4/7>3$.
It is impossible.

By the same reasons, $d_1=3/7$,
$P_1$ has type ($\Atd{1}{}$),
$D_1$ is ($-2$)-curve on $S^{\rm min}$,
$b_1=b_2=6/7$, $F=0$,
$K+B,K^{\rm min}+B^{\rm min}$ and $K_{S'}+B'\equiv 0$,
and the modifications are
crepant.
If $\# D_1\cap C_1=2$,
then we should make at least two blow-ups
in this intersection to disjoint $D_1$
and $C_1$ on $S^{\rm min}$.
According to (A$_2$) for $C_1+C_2$
on $S$ and $S^{\rm min}$, we need to
make at least one blow-up in $R_1$ or $R_2$.
So, $m_1\le 1$.
This contradicts our assumptions.
Thus $D_1$ is tangent with $C_1$.
But again to disjoint $D_1$ and $C_1$
we need two blow-ups.
That leads to the same contradiction.   
As we see later, the latter case is possible for
$m_1=1$ in type (I$_2^2$), or for $\rho=2$
and $m_1=m_2=0$ as will be discussed below.

So, $m_1=1$, and we have contraction
$S^{\rm min}\to S'=\mathbb P^2$
given by the linear system $|C_1|$
on $S^{\rm min}$; $m_1'=m_2'=1$. 
So, $C_1$ and $C_2$ are
lines on minimal $S'=\mathbb P^2$.

If $m_2=1$, we have no contractions onto $C$ 
for $S^{\rm min}\to S'=\mathbb P^2$
and no singularities of $S$ on $C$.
Then $S=\mathbb P^2$, and this is the first exception
(A$_2^1$).
Indeed, if, say, $C_1$ has a singularity,
it gives rise to a curve $E_1\subset S^{\rm min}$
and $S'=\mathbb P^2$,
which intersects $C_2$ on both $S'$ and $S^{\rm min}$.
It is impossible.

Therefore, $m_2=0$.
Moreover, by the same arguments $C_2$
has at most one singularity $P_2$ of $S$, and
$C_1$ has a singularity $P_1$ of $S$
(except for (A$_2^1$)).
Really, $S$  has a singularity $P_2\in C_2$,
because $m_2=0$; 
$P_2$ is a single singularity of $S$ on $C_2$.
It corresponds to 
point $P_2\in C_2\subset S'=\mathbb P^2$
into which we contract and only once a $(-1)$-curve.

So, $F\ge (1/2)D_1$ for some curve 
$D_1\not=C_1$ and $C_2$ on $S$.
Moreover, $D_1$ is a curve on $S'=\mathbb P^2$.
On the other hand,
$P_1$ gives another curve $E_1$ on $S'=\mathbb P^2$ with 
$d_1=\mult{E_1}{B'}\ge ((m-1)/m)b_1$,
where $m\ge 2$ is the index of $K+D$ in $P_1$.
Moreover, by the construction $E_1$ is
a line on $S'=\mathbb P^2$.
Curve $D_1$ is also
a line on $S'=\mathbb P^2$.
Otherwise, $0\ge \deg(K_{S'}+B')\ge 
\deg(K_{S'}+b_1 C_1+b_2 C_2+(1/2)D_1+d_1 E_1)\ge
-3+b_1+b_2+2(1/2)+d_1\ge -3+2\times (6/7)+1+3/7=1/7$.
By the same reasons,
other components of $F$ are contracted on $S'=\mathbb P^2$.
In other words, $\Supp{B'}=C_1+C_2+D_1+E_1$.
Equivalently, $D_1$ is the only component
of $F$ which is passing a non-singular point of $S$ on $C_1$.
Moreover, lines $C_1,C_2,D_1$ and $E_1$ are in general
position in $S'=\mathbb P^2$, because 
the intersection point $P_3$ of
$E_1$ and $D_1$ does not belong to $C$.

Curve $D_1$ crosses $C_1$ in a single non-singular
point of $S$.
We contend that $D_1$ passes $P_1$ on $S$ too.
Indeed, we can increase $B$ to $B''$ in $D_1$ 
in such a way that $K+B''\equiv 0$.
Then $K^{\rm min}+(B'')^{\rm min}$ and $K_{S'}+B'''\equiv 0$
and the modifications are crepant, where
$B'''=b_1 C_1+b_2 C_2+d''D_1+d_1'E_1$
with $d_1'=\mult{E_1}{(B'')^{\rm min}}$ and
with $d''=\mult{D_1}{B''}$, 
and $(B'')^{\rm min}$ corresponds to 
the crepant resolution $S^{\rm min}\to S$.
By the above arguments or the Inductive Theorem,
$d_1'$ and $d''<1$.
However $\delta=d_1'+d''-1\ge 0$,
if we assume that $D_1$ and $E_1$ are disjoint on
$S^{\rm min}$.
Moreover, we have contractions of curves
for $S^{\rm min}\to S'$ onto $D_1$ only over 
the intersection point $P_3$,
because their multiplicities in $(B'')^{\rm min}$ are
non-negative.
Curves on $S^{\rm min}/P_3$ form a chain
with $D_1$ and $E_1$.
So, $D_1$ is a rational $n$-curve with $n\le 0$
on $S$ with a single singularity $P_3$ 
of type $\mathbb A_l$.
Moreover, $n=0$, because $D_1^2>0$ on $S$.
Also $l\ge 2$, because $S$ is not a cone:
it has too many singularities.
In addition, by Lemma~\ref{ecomp}.\ref{ineq}
and our description of the modification near $D_1$,
$\delta\ge d_2=\mult{E_2}{(B'')^{\rm min}}\ge (2/3)d''\ge
(2/3)(1/2)\ge 1/3$, where $E_3/P_3$ in $S^{\rm min}$
intersects $D_1$.
This gives a contradiction:
$0=\deg(K_{S'}+B''')=
\deg(K_{S'}+b_1 C_1+b_2 C_2+d''D_1+d_1'E_1)\ge
-3+b_1+b_2+d''+d_1'\ge -3+2\times (6/7)+1+1/3=1/21$.
Note another contradiction, that the $(-1)$-curve$/P_3$,
which is non-exceptional on $S$, does not intersect
$C_2$ on $S$.

Therefore, $D_1$ passes $P_1$,
and, by the log terminal property of $K+D$ and
since $D_1$ intersects $E_1$ on $S^{\rm min}$,
$P_1$ is a simple singularity with
a single ($-l$)-curve $E_1$ on $S^{\rm min}/P_1$.
In addition, $\Supp F=D_1$.
Arguing as above, we can verify that
$P_2$ is a Du Val singularity of type $\mathbb A_l$,
or $S^{\rm min} \to S'$ contracts only ($-2$)-curves
and one ($-1$)-curves -- 
successive blow-ups of $P_2\in S'=\mathbb P^2$.
In addition $F=0$ near $P_2$.
Hence $P_2$ has index $l+1$ for $K+D$,
and $l\ge 2$, because $S$ is not a cone
(except for (A$_2^1$)). 
This gives $d_2=\mult{E_2}{B'}=(l/(l+1))b_2\ge 
(l/(l+1))(6/7)$
by Lemma~\ref{ecomp}.\ref{ineq} in
divisor $E_2/P_2$ on $S^{\rm min}$ intersecting $C_2$.
So, in this case
$0\ge (K+B.C_2)\ge -2+b_1+\mult{D_1}{F}+d_2\ge
-2+6/7+1/2+(l/(l+1))(6/7)$.
Hence $l/(l+1)\le 3/4$ and $l\le 3$.
This gives types (A$_2^4$) and (A$_2^5$)
for $l=2$ and $l=3$ respectively.
The same inequality with $l=2$
gives $\mult{D_1}{F}\le 4/7$
and so $F=(1/2)D_1$.
If $l=3$, $K+(6/7)(C_1+C_2)+(1/2)D_1\equiv 0$.

Finally, we prove that the cases
with $m_1=m_2=0$ are impossible.
First, we verify that each $C_i$
has at least two singularities of $S$.
Otherwise $F\ge (1/2)D_1$ where 
curve $D_1\not= C_1$ and $C_2$ intersects
one of theses curves, say $C_1$, in a non-singular point
and only in this point.
Curve $C_1$ has a singular point $Q_1$ of $S$,
and $C_2$ does so for $P_1$,
because $C^2_i>0$.
Let $F_1/Q_1$ and $E_1/P_1$ be
respectively curves on $S^{\rm min}$
which intersect $C_1$ and $C_2$.
They give different generators $F_1$ and $E_1$
of $S'=\mathbb P^1\times\mathbb P^1$
with multiplicities $f_1=\mult{F_1}{B'}\ge (1/2)b_1$
and $e_1=\mult{E_1}{B'}\ge (1/2)b_1$.
Moreover, if $\Supp{F}$ passes $Q_1$,
equivalently, a component $D_i$ of $F$ passes $Q_1$,
then $f_1\ge (1/2)b_1+(1/2)\mult{D_i}{B}\ge
3/7+1/4$, but
$0\ge (K+B.C_1)\ge -2+b_2+f_1+\mult{D_i}{B}\ge
-2+6/7+3/7+1/4+1/2=1/28$.
By the same reasons, $D_1$ intersects $C_1$
only in one point.
Thus $\Supp{F}=D_1$ does not pass $Q_1$.
If $D_1$ does not intersect $C_2$ in 
a non-singular point of $S$,
then $D_1$ is also a generator $D_1\sim F_1$.
However, $D_1$ intersects transversally
in one point $P$ another generator $E_1$
on $S'=\mathbb P^1\times \mathbb P^1$,
and $\Supp{B'}=D_1+E_1$ near $P$.
On $S^{\rm min}$, $D_1$ and $E_1$ cannot be
disjoint by a chain of rational curves$/P$,
because we can assume as above that $K+B\equiv 0$
and the modifications are crepant.
Then a ($-1$)-curve$/P$ will be
a curve on $S$ non-intersecting $C_1$.
It is impossible.
Hence $D_1$ passes each singular point of $S$ on $C_2$
and we have two of them, or $D_1$ intersects $C_2$
in a non-singular point of $S$.
The former case is impossible:
$0\ge (K+B.C_2)\ge -2+b_2+2\times((1/4)+(1/2)b_2) 
\ge-2+6/7+1/2+6/7=3/4$
as above with $f_1$.
So, $D_1$ intersects $C$ on $S$ in 
two non-singular points: one on each $C_i$.
Then $D_1\sim F_1+E_1$ or has bi-degree $(1,1)$
on $S'=\mathbb P^1\times\mathbb P^1$.
Moreover, $B'=e_1 E_1+f_1F_1+(\mult{D_1}{F} )D_1$,
and $D_1$ passes the intersection point $P$ of
$E_1$ and $F_1$, because 
$D_1$ does not pass singular points of $C$.
Assuming as above that $K+B\equiv 0$ and
the modifications are crepant, we can check
that $S^{\rm min}\to S'$ contracts only
curves over $P$.
Since $\rho(S)=1$, all these curves but one
are contracted on $S$, to singularities
$P_1$, $Q_1$ and maybe one on $D_1$.
However, this is only possible when we contract
$E_1$ and $F_1$ after the first monoidal
transform in $P$ which does not 
produce singularities of $S$ at all.
Otherwise as above, a ($-1$)-curve$/P$ will be
a curve on $S$ non-intersecting some $C_i$.

So, $F$ intersects $C$ only in singular points of $S$.
Hence we have at least four singularities of $S$:
$Q_1,Q_2\in C_1$ and $P_1,P_2\in C_2$.
They give respectively different generators
$F_1\sim F_2$ and $E_1\sim E_2$
with multiplicities $f_i=\mult{F_i}{B'}\ge (1/2)b_1$
and $e_i=\mult{E_i}{B'}\ge (1/2)b_2$.
If for two of these multiplicities, say $e_1$ and $f_1$,
$e_1+f_1<1$, then we cannot disjoint $E_1$ and $F_1$
on $S^{\rm min}$ under assumption $K+B\equiv 0$.
So, we may suppose that $e_1$ and $e_2\ge 1/2$.
Then by Lemma~\ref{ecomp}.\ref{ineq},
$e_1$ and $e_2\ge (2/3)b_1\ge 4/7$.
This gives an equation in the inequality
$0\ge (K+B.C_2)\ge -2+b_1+e_1+e_2\ge
-1+6/7+2\times (4/7)=0$.
Moreover, $e_1=e_2=4/7, B=(6/7)(C_1+C_2)$, $F=0$,
$K+B\equiv 0$ and the modifications are crepant,
Indeed, since $1/4+(1/2)(6/7)=19/28>4/7$,
$C_2$ has exactly two singularities
$P_1$ and $P_2$ of the index $3$ for $K+D$
or $K+C_2$.
Thus they are Du Val of type $\mathbb A_2$.

On the other hand, $S$ has two singularities
$Q_1$ and $Q_2\in C_1$.
We can assume that $Q_2$ is not simple Du Val.
If $Q_1$ is also not simple Du Val,
then again both are Du Val of type $\mathbb A_2$.
So, $B'=(6/7)(C_1+C_2)+(4/7)(E_1+E_2+F_1+F_2)$,
where $F_i$ is curves with $d=4/7$ over $Q_i$.
All curves $C_1\sim E_1\sim E_2$ and
$C_2\sim F_1\sim F_2$ are generators of
corresponding rulings of $S'=\mathbb P^1\times\mathbb P^1$.
On $S^{\rm min}$ there exists a curve $E'$
with $\mult{E'}{B^{\rm min}}=2/7$.
For instance, the second curve in 
a minimal resolution of $P_1$.
Such curves could only be over
the intersection of $P=E_i\cap F_i$,
which is impossible because they have only
one curve $E''/P$ with
$a(E'',K_{S'}+B')\le 1$.
It is a blow-up of $P$ with
$a(E'',K_{S'}+B')=6/7$.

Therefore $Q_1$ is simple Du Val,
$f_1=(1/2)b_1=3/7,f_2=2-b_2-f_1=2-6/7-3/7=5/7$, and
$B'=(6/7)(C_1+C_2)+(4/7)(E_1+E_2)+(3/7)F_1+(5/7)F_2$. 
However, this is possible only for a surface $S$
with $\rho(S)=2$.
Indeed,
to disjoint $F_2$ and $E_1$ we should make
two successive monoidal transformations:
first in $F_2\cap E_1$ which gives $E_3$,
and then in $F_2\cap E_3$ which gives $E_4$.
Both curves $E_1$ and $E_3$ are ($-2$)-curves
on $S^{\rm min}$: over points $F_1\cap E_i$,
the contraction $S^{\rm min}\to S'$ is just
one monoidal transform.
On the other hand, $F_2$ will be a ($-4$)-curve
and $E_4$ is a contractible ($-1$)-curve
passing $Q_2$ and $P_1$.
Hence $\rho (S)\ge 2$.
One can check that $\rho (S)=2$, and
after a contraction of $E_4$ on $S$
we get Case~(I$_2^2$) in \ref{ecomp}.\ref{emainth}.3.

Now we can find minimal complements
in~\ref{ecomp}.\ref{emainth}.3 as in 
\cite[Example 5.2.1]{Sh2},
however we did it only in a few trivial cases.
Indeed, by Proposition~\ref{ecomp}.\ref{logter}, 
we should care only about
the numerical property, not singularities.
\end{proof}

\section{Classification of surface complements}
\label{clasf}

Take $r\in N_2$.
A classification of $r$-complements
or of complements of index $r$ 
means a classification of
surface log canonical pairs $(S/Z,B)$ with
$r(K+B)\sim 0/Z$.
We assume that $S/Z$ is a contraction:
in the local case, according to the Main Theorem, 
such that there exists
a log canonical center$/P$ (cf. Example~\ref{int}.\ref{tcomp}
and Section~\ref{lcomp}).
This classification implies 
a classification of log surfaces $(S/Z,B')$
with such complements.
For instance, the minimal complementary index
$r$ is an important invariant of $(S/Z,B')$.

Moreover, for exceptional $r\in EN_2=\{7,\dots\}$
$K+B$ is to be assumed Kawamata log terminal
and $S$ is complete with $Z=\pt$.
Such cases are bounded.
They have been partially described in Section~\ref{ecomp}.
Here we focus our attention on basic invariants for
regular complements with $r\in RN_2=\{1,2,3,4,6\}$.

For indices $1,2,3,4$ and $6$,
types of complements are denoted 
respectively by $\mathbb A_m^n,\mathbb D_m^n,
\mathbb E3_m^n,\mathbb E4_m^n$ and $\mathbb E6_m^n$,
where 
\begin{itemize}
\item
$n$ is the number of reduced 
(and formally) irreducible components
in $B$ (over a neighborhood of a given $P\in Z$
in the local case), and
\item
$m$ is the number of reduced exceptional
divisors of $B^{\rm min}$ 
on a crepant minimal formally log terminal resolution
(in most cases, this is a minimal resolution
cf. Example~\ref{clasf}.\ref{notation} below) 
$(S^{\rm min},B^{\rm min})\to (S,B)$ 
(over a neighborhood of $P\in Z$),
equivalently, the number of exceptional
divisors on the minimal resolution
with the log discrepancy $0$.
\end{itemize}
We also assume that, 
for the types $\mathbb A_m^n$ and $\mathbb D_m^n$,
$\Supp{B^{\rm min}}$ is a connected and singular curve,
otherwise we denote them respectively
$\mathbb E1_m^n$ and $\mathbb E2_m^n$.
Equivalently, exactly in types $\mathbb E1-2_m^n$
among types $\mathbb A_m^n$ and $\mathbb D_m^n$,
the number of exceptional divisors with 
the log discrepancy $0$ is finite.
The same numerical invariants can be defined
in any dimension up to the LMMP.
But as we see in the next section,
a simplicial space associated to reduced components
is a more important invariant.

\refstepcounter{subsec}
\label{clasfth}
\begin{theorem} 
If $m=n=0$, then $Z=\pt$, and 
$K+B$ is Kawamata log terminal
of a regular index $r\in RN_2$.
It is possible only for types $\mathbb Er_0^0$.
The complements of this type are 
bounded when $B\not=0$ 
or $S$ has non-log terminal point {\rm \cite{Al}}.
Otherwise $B=0$ and $S$ a complete
surface with canonical singularities
and $r K\sim 0$ (their classification
is well known up to a minimal resolution
and can be found in any textbook on
algebraic surfaces, e.~g., {\rm \cite{BPV}} and 
{\rm \cite{Shf}}).

In the other types we suppose that $m+n\ge 1$.
Then they have a non-empty locus of
the log canonical singularities $\LCS{(S,B)}$.
It is the image of $\LCS{(S^{\rm min},B^{\rm min})}=
\Supp{\lfloor B^{\rm min}\rfloor}$ and it has
at most two components: 
Two only
for exceptional types $\mathbb Er_n^m$ with
$n+m=2$ and in the global case.

If $\LCS{(S,B)}$ is connected, it is a point
if and only if $n=0$.
Otherwise it is a connected curve
$C=\Supp{\lfloor B\rfloor}$ with at most nodal
singularities and of the arithmetic genus $g\le 1$
for any sub-curve $C'\subseteq C/P$.
Moreover, if $g=1$ then $C'=C$,
this case is only possible 
for types $\mathbb E1_1^0$ and $\mathbb A_m^n$ 
with $n\ge 1$, whereas $C=C'$ is 
respectively a non-singular curve
of genus $1$ or a Cartesian leaf when $n=1$
and $C$ is a wheel of $n$ rational curves 
when $n\ge 2$.
In all other cases $C$ is a chain of 
$n$ rational curves.

The singularities of $(S,B)$ outside
of $\LCS{(S,B)}$ are log terminal of
index $r$.
In particular they are only canonical when $r=1$;
moreover, only of type $\mathbb A_i$ when $m+n\ge 1$.

If $B=0$ then $n=0$ and $\LCS{(S,B)}$
is the set of elliptic singularities of $S$.

For types $\mathbb A_m^n$ and $\mathbb D_m^n$, 
any natural numbers $m,n$ are possible. 
If $n=0$, then it is a global case
with $B=0$, $K\sim 0$ and $S$ has
a single elliptic singularity. 
  
In the exceptional types $\mathbb Er_m^n$,
$n+m\le 2$ and any $m,n$ under this condition
are possible; $n+m=2$ is only possible
in the global case.
The number of the connected components of $\LCS{(S,B)}$
is $m+n$.
Equivalently, each such component is irreducible.
Moreover, it is a point or a non-singular curve,
respectively, of genus $1$ for types $\mathbb E1_m^n$
and of genus $0$ for the other types,
when $C/P$.
\end{theorem}

\begin{proof}
The most difficult part is
related to connectedness
\cite[Theorem~6.9]{Sh2}.
Other statements follow
an adjunction, except for the statement
on types of canonical singularities
for $1$-complements $(S/Z,B)$.
Essentially it was proved in
Section~\ref{indcomp}. 

So, let $(S/Z,B)$ be a $1$-complement.
After a formal log resolution we can assume that
$\LCS{(S,B)}=B=C$ is a reduced curve
and $S$ non-singular near $C$.
We can also assume that $C$ is minimal, i.~e.,
does not contain $(-1)$-curves.
We verify that the singularities
of $S$ have type $\mathbb A_i$
using an induction on extremal 
contractions.
By the LMMP we have an extremal contraction
$g:S\to T/Z$ with respect to $K$ if
there exists a curve$/P$ not in $C$.
If this is a contraction of a curve $C'$ to a point,
then it is to a non-singular point and 
$S$ has only singularities of type $\mathbb A_i$
near $C'$, because $S$ has only canonical singularities.
If this is a contraction of a fibre type,
then it is a ruling which can have singularities
only when $C$ has an irreducible component $C'$
as a double section.
Then the only possible singularities are
simple double. 
If $T=\pt$, then components of $C$ are ample
and $S=\mathbb P^2$.
Note that if we have no contractions we have
no singularities, because the latter ones
are apart of $C$. 
\end{proof}

Of course, our notation is similar 
to the classical one 
(however, with some twists).

\refstepcounter{subsec}
\label{notation}
\begin{example}
For instance, a singularity $P\in (S/S,H)$
of type $\mathbb A_m$
with the generic hyperplane through $P$
has type $\mathbb A_m^2$
in our notation.
But type $\mathbb D_m$ corresponds to
$\mathbb D_{m-2}^1$ with
some reduced and irreducible $H$.

We have more differences for 
elliptic fiberings $(S/S,E)$.
For instance, the Kodaira type $_m\rm I_b$
is our $\mathbb A_0^b$.

Let $(S/S,L_1+(1/2)(L_2+L_3))$ be
a singularity as in the plane $S=\mathbb P^2$
in the intersection of three lines $L_i$.
Then it has type $\mathbb D_1^1$,
because its minimal log resolution is
a monoidal transform in this point.
\end{example}
 
\refstepcounter{subsec}
\label{toric}
\begin{example}
Each toric variety $X$ has
a natural $1$-complement structure
$(X,D)$ where $D$ is $D=\sum D_i$
with the orbit closures $D_i$.
So, the number of elements in this sum $n$
is the number of edges in the fan.

A toric surface $S$ with $n$ ages has type
$\mathbb A_m^n$.
In addition, $n=\rho(S)+2$.
This characterizes toric surfaces.
\end{example}

\refstepcounter{subsec}
\label{rhob}
\begin{theorem} 
Let $(S/Z,B)$ be with
log canonical $K+B$ and nef$/Z$ divisor $-(K+B)$.
Then $\rho(S/Z)\ge \sum b_i-2$,
where $\rho(S/Z)$ is the rank of Weil
group modulo algebraic equivalence$/Z$, or 
just the Picard number
when the singularities of $S$ are rational.
Moreover, $=$ holds if and only if
$K+B\equiv 0$ and $S/Z$ is formally toric
with $C=\lfloor B\rfloor\subseteq D$. 

In addition, in the case $=$ and reduced $B=C$,
$(S/Z,C)$ is formally toric with $C=D$ 
(see Example~\ref{ecomp}.\ref{toric}).
\end{theorem}

{\it Formally toric\/}$/Z$ means formally
equivalent to a toric contraction,
or locally$/Z$ in analytic topology, when
the base field is $\mathbb C$.

\begin{proof}
First, we can assume that $\LCS{(S,C)}\not=\emptyset$.
In the local case we can do this adding
pull back divisors as in the proof of the General Case
in Theorem~\ref{lcomp}.\ref{lmainth}.
In the global case, after contractions,
we can assume that $\rho(S)=1$.
If $\LCS{(S,C)}=\emptyset$, the inequality
will be improved after contractions.
If $B$ has at most one component $C_i$
with $b_i=\mult{B_i}{B}>0$,
then $\rho(S)=1>b_i-2$.
Otherwise we have at least two curves
$C_i$ and $C_j$ with $b_i$ and $b_j>0$.
We can also assume that $K+B\equiv 0$.
If $\LCS{(S,B)}=\emptyset$ we can change
$b_i$ and $b_j$ in such a way that
$K+B\equiv 0$ and $b_i+b_j$ is not
decreasing.
Indeed $K+B\equiv 0$ gives
a linear equation on $b_i$ and $b_j$. 
Then we get $\LCS{(S,B)}\not=\emptyset$, or
we get $b_i$ or $b_j=0$.
An induction on the number of curves
in $\Supp{B}$ gives the log singularity
or the inequality. 

Second, we can replace $(S,C)$ by its
log terminal minimal resolution 
$(S^{\rm lt},C^{\rm lt})$
over $C=\LCS{(S,B)}\not=0$.
We preserve all the statements.
The contraction will be toric because
it contracts curves of $D$.

If every curve $C'/P$ is in $C$,
we have the local case and
by the adjunction and \cite[Corollary~3.10]{Sh2}, 
we reduce our inequality
to a $1$-dimensional case on $C''/P$.
In addition, for $=$, $S/P$ is non-singular
and toric which is possible to check case by case.
Here we use the monotonicity 
$(m-1)/m+\sum k_i d_i/m\ge d_i$
when $k_i\ge 1$, and even $>d_i$,
when $m\ge 2, k_i\ge 1$ and $1>d_i$
(cf. [ibid]).

Third, we could assume that $K+B\equiv 0$
on $C/P$ after birational contractions.
This improves the inequality.
If we return to a Kawamata log terminal case,
we can find a complement $K+B'\equiv 0/Z$
with $B'\ge B$.
This again improves the inequality.

By \cite[Theorem~6.9]{Sh2}, we assume that
$C/P$ is connected.
Otherwise $Z=\pt$, $C=C_1+C_2$, and after contractions
we can assume that each fractional component, i.~e.,
each component of $B-C$, intersects some $C_i$
(cf. the arguments for the connected case).
Then we reduce the problem to a $1$-dimensional case
on $C_i$.

Therefore we assume that there exists a curve $C'/P$
not in $C$.
Then we have an extremal contraction
$g:S\to T$ which is numerically non-negative for 
a divisor $H$ with $\Supp H=C$ and $g$
is numerically negative for $K+B-\varepsilon H$
with some $\varepsilon>0$.
If $C$ has an exceptional type we take such
an $H$ that is negative on $C$.
Otherwise we assume that $H$ is nef on $C$,
and even ample in the big case.
So, such a $g$ preserves birationally $C$,
whenever birational, or $H$ is ample
and we consider this case as
a contraction to $Z=\pt$ below.
After birational contractions,
$f$ has a fibre type.
If it is to a point, then $Z=\pt$
and $C$ has at most two components
which are intersected by other components of $B$.
We can choose them
and reduce the problem to a $1$-dimensional case
as above.
If $g$ is a ruling we can do similarly
when $C$ has an ample component.
Otherwise $C$ is in a fibre of $g$.
As in the first step, we can assume that we have
at most one other fibre component.
This implies the inequality.
Otherwise we obtain a case when
$C$ is not connected.
Note, that we preserve inequality
only when we contract a curve $C_i/P$
with $(K+B.C)=0$ and $b_i=1$.
Such a transform preserves the formally toric property.
So, $D$ contains $C$ always.
\end{proof}

We hope that in general
$\rho(X/Z)\ge - \dim X+\sum b_i$,
where $\rho$ is the Weil-Picard number,
i.~e., the rank of the Weil divisors
modulo algebraic equivalence.
Moreover, $=$ holds exactly for formally toric varieties
and $\lfloor B\rfloor\subseteq D$.
For instance, this implies 
that locally $\sum b_i\le \dim X$ when
the singularity is $\mathbb Q$-factorial
and $\rho=0$.
If the singularity is not $\mathbb Q$-factorial
we have a stronger inequality
$\sum b_i\le \dim X-1$ for
$B=\sum b_i D_i$ with $\mathbb Q$-Cartier $D_i$
(cf. \cite[Theorem~18.22]{KC}).

\refstepcounter{subsec}
\begin{corollary} 
Let $(S/Z,C)$ be as in Theorem~\ref{clasf}.\ref{rhob}.
The following statements are equivalent.
\begin{itemize}
\item
$(S/Z,C)$ is a surface $1$-complement of 
type $\mathbb A_m^n$ with $n=\rho(S/Z)+2$;
\item
$\rho(S/Z)= n-2$, where $n$ is 
the number of (formally) irreducible components
in $C$ ; and
\item
$(S/Z,C)$ is formally toric.
\end{itemize}
\end{corollary}

For instance, if $\rho(S)=1$ and $Z=\pt$, 
then a $1$-complement $(S,C)$ of type $\mathbb A_m^n$
is toric with $D=C$, if and only if $n=3$.
In other cases $n\le 2$.

Most of the above results work over
non-algebraically closed fields of
the characteristic $0$.

\refstepcounter{subsec}
\label{arith}
\begin{example}
If $C$ is a non-singular curve
of genus $0$, it always has a $1$-complement.
But it has type $\mathbb A_0^2$
only when it has a $k$-point.
Otherwise it has type $\mathbb A_0^1$
and its Fano index is $1$.
\end{example}

A complement with connected $\LCS{(S,B)}$
can be called a {\it monopoly\/}.
Other complements are {\it dipoles\/}.

\refstepcounter{subsec}
\label{dipole}
\begin{theorem}
Any exceptional complement $(S/\pt,B)$
of type $\mathbb Er_0^2$
has a ruling $g:S\to Z$
with a normal curve $Z$ and with two sections
in $\LCS{(S,B)}$;
the genus of $Z$ is $1$ for $(r=1)$-complements  
and $0$ in other cases.
\end{theorem}

\begin{proof} We obtain the ruling
after birational contractions which are
with respect to a curve $C_i$
in $\LCS{(S,B)}$, $C_i^2\le 0$.
Cf. Proof of the Inductive Theorem: Case I.

In addition, if $g(Z)\ge 1$ then $B=0$
and we have a $1$-complement $(S,B)$.
\end{proof}

\refstepcounter{subsec}
\label{dipolec}
\begin{corollary}
Let $(S,B)\to (S',B')$ be
a normalization of 
a connected semi-normal log pair $(S',B')$
with $B'$ under {\rm (M)\/}.
Then $(S',B')$ has an $r$-complement,
if $(S,B)$ has a complement of type $\mathbb Er_m^n$.
\end{corollary}

\begin{proof}
The divisors $B^+$ in the normal part of $S$
belong to fibres  of the contraction $g$
on components of the normalization of $S$
after a log terminal resolution.
The latter has the same $r$ for each component of
$S$ because they are induced from 
curves of non-normal singularities on $S'$.
\end{proof}

For the other dipoles, $\LCS{(S,B)}$ is
a pair of points  or a point and a curve.

\refstepcounter{subsec}
\label{rull}
\begin{remark}
The ruling induces a pencil $\{C_t\}$ of rational curves
through the points.
Similarly, in other cases we can find
\begin{description}
\item{(PEN)} {\it a log proper pencil $\{C_t\}$ 
of log genus $1$ curves\/}, i.~e., $C_t$ does not
intersect $\LCS{(X,B)}$, when its normalization
has genus $1$, and, for the corresponding map 
$g_t:(\mathbb P^1,0+\infty)\to (X,B)$ onto $C_t$,
$g_t(Q)\not\in \LCS{(X,B)}$ when $C_t$ has genus $0$,
$t$ is generic, and $Q\not=\infty$ and $\not= 0$.
\end{description}
This implies easy cases in the Keel-McKernan Theorem
on the log rational covering family \cite[Theorem~1.1]{KM}.
The difficult cases are exceptional and bounded.
Perhaps it can be generalized in a weighted
form for fractional boundaries or m.l.d's.
\end{remark}

\section{Classification of 3-fold log canonical singularities}
\label{csing}

\refstepcounter{subsec}
\label{singth}
\begin{theorem}
Let $(X/Z,B)$ be a birational contraction
$f:X\to Z$ of a log 3-fold $X$
\begin{itemize}
\item 
with boundary $B$ under {\rm (SM)\/} and
\item
nef $-(K+B)$.
\end{itemize}
Then it has an $n$-complement $(X/Z,B^+)$
over a neighborhood of any point $P\in Z$
such that
\begin{itemize}
\item
$n\in N_2$ and
\item
$K+B^+$ is not Kawamata log terminal over $P$.
\end{itemize}
\end{theorem}

\begin{statement1}
We can replace {\rm (SM)\/} by 
\begin{description}
\item{{\rm (M)''}} 
the multiplicities $b_i$ of $B$
are {\rm standard\/}, i.~e., $b_i=(m-1)/m$ 
for a natural number $m$, or
$b_i\ge l/(l+1)$ where $l=\max\{r\in N_2\}$.
\end{description}
\end{statement1}

\refstepcounter{subsec}
\label{singl}
\begin{lemma} Let $(S/Z,B)$ be 
a surface log pair and
$(S^{\rm min},B^{\rm min})$ be
its crepant minimal resolution.
Suppose that $B$ satisfies {\rm (SM)}.
Then
$$K+B \ n-complementary\ \ \Longrightarrow
\ K^{\rm min}+B^{\rm min}\ \ n-complementary.$$
Moreover, we could replace $S^{\rm min}$ by
any resolution $S'\to S$ with subboundary
$B'=B^{S'}$. 

For any $n\in N_2$,
we can replace {\rm (SM)\/} by {\rm (M)''\/}.
\end{lemma}

\begin{proof}
Follows from the Main Lemma~\ref{gcomp}.\ref{mlem3}
or can be done it the same style.

The last statement follows from 
the Monotonicity~\ref{indcomp}.\ref{monl1}.
\end{proof}

\begin{proof}[Sketch Proof of 
Theorem~\ref{csing}.\ref{singth}]
First, we can assume that $K+B$ is
strictly log terminal$/Z$ and $B$ has reduced part
$S=\lfloor B\rfloor\not=\emptyset/P$.
For this we add a multiple of an effective divisor
$D=f^*(H)$ for a hyperplane section $H$ of $Z$
through $P\in Z$.
However, $B+d D$ may contradict (SM) in $D$.
We drop $D$ after a log terminal resolution
of $(S,B+d D)$.
In its turn this can spoil the nef condition
for $-(K+B)$.
This will be preserved when $D$ or
even a divisor $D'\ge D$ below is nef$/P$.
If not we can do this after modifications
in extremal rays on which $D'=B^+-B\ge D$ is negative,
where $B^+$ is a complement
for $B+d D$, i.~e., $K+B^+\equiv 0/P$.
We drop then $D'$.
To do the log flops with respect to $D'$ 
we use the LSEPD trick
\cite[10.5]{Sh2}.
Finally, if $K+B$ is not strictly log terminal,
it holds after a log terminal resolution of $(S,B)$. 

Second, $(S/Z,B_S)$ is a semi-normal and connected
surface over a neighborhood of $P$
where $B_S$ is non-singular in codimension $1$
and again under (SM) or, for \ref{csing}.\ref{singth}.1,
under (M)''.
This follows from the LMMP or \cite[Theorem 17.4]{KC} and
Corollary~\ref{indcomp}.\ref{adjunct}.
Adjunctions $(K+B)|_{S_i}$ on each component
are log terminal \cite[3.2.3]{Sh2}.

We have a complement $(S/Z,B_S^+)$ on $S/Z$
and hence on each $S_i$.
This gives (EC)'.

Third, we can glue complements
from the irreducible components of $S$.
If one of them has no regular complements,
then $S$ is normal and 
there is nothing to verify.
In the other cases we have $r$-complements
with $r\in RN_2$.
Moreover if we have a complement of type
$\mathbb Er_m^n$ on some component $S_i$, 
then $m=0$ by
the log terminality of the adjunction
$(K+B)|_{S_i}$,
and $S$ is a wheel or a chain of its
irreducible components.
Then we can glue complements by 
Corollary~\ref{clasf}.\ref{dipolec}.
Finally, we have the complements of type
$\mathbb D_m^n$ with $r=2$.
They are induced from a $1$-dimensional
non-irreducible case when we always have
a $2$-complement (cf. \cite[Example~5.2.2]{Sh2}).
(However if we have non-standard coefficients $b_i\le 2/3$
we need to use higher complements
(cf. Lemma~\ref{indcomp}.\ref{compl3}).)

Finally, we can act as in the proof of
\cite[Theorem~5.6]{Sh2} (cf. \cite[Theorem 19.6]{KC}).
So, an $r$-complement of $(X/Z,B)$ is induced from
the $r$-complement of $(S/Z, B_S)$.
We can lift the $r$-complement on any log resolution
by Lemma~\ref{csing}.\ref{singl}.
\end{proof}

In particular, we divide $3$-fold birational contractions 
into two types: exceptional and regular.
By Mori's results all
small contractions in the terminal case are 
formally regular.

This time the exceptions are not bounded.
For instance, if $(X,S)$ is a simple
compound Du Val singularity:
$$
x^2+y^2+z^2+w^d=0
$$
with quadratic cone surface $S$
given by $x^2+y^2+z^2=0$.
It is (not formally) an exceptional complement.
Such singularities are not
isomorphic for different $d$ even formally. 
However, they have many common finite invariants:
the m.l.d., the index of complement,
the index of $K$, etc..
They are bounded up to isomorphism of 
a certain degree or order.

\refstepcounter{subsec}
\label{exbound}
\begin{corollary}
Under the assumptions of Theorem~\ref{csing}.\ref{singth},
for any $\varepsilon>0$,
the exceptional contractions $(X/Z,B)$ 
and their complements $(X/Z,B^+)$ are bounded
with respect to the m.l.d. and discrepancies,
when $K+B$ is $\varepsilon$-log terminal:
over a neighborhood of $P\in Z$,
the set of the m.l.d.'s
$a(\eta,B,X)$ for points $\eta$ and 
the discrepancies $a(E,B^+,X)$, 
$a(E,B,X)\le \delta$
for any $\delta>0$ is finite.

It holds also under (M)'', if the contraction
is not divisorial. 
Otherwise we should assume that
the set of $b_i\le 1-\varepsilon$ is finite.
\end{corollary}

\begin{proof}
According to our assumptions, $b_i$ belongs
to a finite set.
Indeed,
if $b_i\le (1-\varepsilon)$.

So, it is enough to verify the finiteness
for discrepancies in exceptional divisors $E$ of $X$.
Indeed, the latter implies that 
we can consider the m.l.d.'s $>3$.
Such does not exist.

Note that $a=a(E,B^+,X)=a(E,B^+,Y)$ form
a finite set because $K+B^+$ has
a finite set of indices $n\in N_2$.

The index $N$ of $S$ is bounded$/P$ as well
because it is bounded locally on $Y$. 
So, we assume that $N$ is the universal index.
To verify the local case, after 
a $\mathbb Q$-factorialization, we can suppose
that $X$ is $\mathbb Q$-factorial.
We assume this below always.
Then the boundedness follows from the boundedness
of quotient singularities on $S$ by the exceptional
property.
In particular, for any point $Q\in S$,
the local fundamental groups of $S\setminus Q$
are bounded.
Along curves in $S$ the index is bounded ($\le 6$)
by \cite[Proposition~3.9]{Sh2}.
After a covering branching over such curves,
$S$ does not pass through codimension $2$
singularities of $Y$.
Then we can argue as in the proof of 
\cite[Corollary~3.7]{Sh2} (cf. \cite[Lemma~1.1.5]{P1}).
The index in such singularities will be bounded 
by orders of cyclic quotients of the fundamental groups.

This also implies that any divisor near $S$ 
has a bounded index.

Now we consider discrepancies $a(E,B,X)$.
Especially, for $E=S/P$.
If it is not exceptional, the finiteness follows
from (SM) and the $\varepsilon$-log terminal property.
The same holds for the other non-exceptional $E$ on $X$.
For exceptional $E=S$, and for any other exceptional $E$,
to compute discrepancies,
we choose an appropriate strictly log terminal model $g:Y\to X$,
on which $S$ is the only exceptional divisor$/X$.
Then $B^Y=g^{-1}B+(1-a(S,B,X))S=
B_Y^+-a(S,B,B)S-D$ where $D=g^{-1}(B^+-B)$ is effective.
Moreover, the multiplicities of $D$ and $B^+-B$
form a finite set.

Take rather generic curve $C\subset S/X/P$.
Then 
\begin{eqnarray*}
0=(K_Y+B^Y.C)&=(K_Y+B_Y^+.C)-a(S,B,X)(S.C)-(D.C)\\
&=-a(S,B,X)(S.C)-(D.C) 
\end{eqnarray*}
implies
$$
a=a(S,B,X)=-(D.C)/(S.C).
$$
Since $a>\varepsilon>0$, and $(D.C)=(D|_S.C)$ belongs
to a finite set,
$(S.C)$ is bounded from below.
On the other hand $(S.C)<0$.
So, $(S.C)$ and $a(S,B,X)$ belong
to a finite set.

This implies the statement for discrepancies
$a(E,B,X)=a(E,B^Y,Y)=a(E,B^+,Y)+
a(S,B,X)\mult{E}{S}+\mult{E}{D}$ with 
centers near $S$ or intersecting $S$.

For the other centers, it is enough to verify
that the index of $D$ is also bounded there.
By \cite[Theorem~3.2]{Sh3}, we need
to verify that the set of exceptional divisors
with discrepancies $a(E,B_Y,Y)<1+1/l$ is bounded
(cf. the proof \cite[Proposition~4.4]{Sh3}).
Since $D$ is effective and $K_Y+B^+$ has index $\le l$,
then $a(E,B_Y^+,Y)\le 1$.
We need to verify that such $E$,
with centers not intersecting $S$, are bounded.
We will see that this bound has 
the form $\le A(1/\varepsilon)$.

Take a terminal resolution $W/Y$ of the above 
exceptional divisors for $K_Y+B_Y^+$.
It does not change the intersection $(S.C)$
for any curve $C\subset S/Z/P$.
(This time $C$ may be not$/X$.)
As above we have inequality:
$$
0\ge (K+B.g(C))=(K_Y+B^Y.C)=-a(S,B,X)(S.C)-(D.C),
$$
or $(S.C)\ge -(D.C)/a(S,B,X)>-(D.C)/\varepsilon$.
So, $(S,C)$ belongs to a finite set even if
we assume that $C$ is ample on $C$ and on 
its other such models of $S$.

Then we apply the LMMP to $K_W+B_W^+-S$.
More precisely, we make flops for $K_W+B_W^+$ with
respect to $-S$, or in extremal rays $R$
with $(S.R)>0$.
So, we decrease $(K_W+B_W^+-S.C)$ and
increase $(S.C)$ and strictly when 
the support $|R|$ has divisorial intersection 
with $S$.
So, the number of exceptional divisors on $W/X$
is bounded, because we contract all of them
during such a LMMP.
\end{proof}

 We have proved more.

\refstepcounter{subsec}
\label{exboundf}
\begin{corollary}
For the exceptional and 
$\varepsilon$-log terminal $(X/Z,B)$,
the fibres  $f^{-1}P$ are bounded.
In particular, if $X/Z$ is small
the number of curves$/P$ is bounded.
\end{corollary}

\refstepcounter{subsec}
\begin{corollary}
For each $\varepsilon>0$,
there exists a finite set $M(\varepsilon)$ 
in $(+\infty,\varepsilon]$
such that $(X/Z,B)$ is not exceptional,
whenever the m.l.d. of $(X,B)/P\ge \varepsilon$ 
and is not in $M(\varepsilon)$.
\end{corollary}

\begin{proof}
$M(\varepsilon)$ is the set of the m.l.d.'s
for the exceptional complements which are $>\varepsilon$.
\end{proof}

\refstepcounter{subsec}
\begin{corollary}
There exists a natural number $n$ such that
any small contraction $X/Z$ of
a $3$-fold $X$ with terminal singularities
has a regular complement, whenever
it has a singularity $Q/P$ in which $K$
has the index $\ge n$. 

A similar bound exists for
the number of curves$/P$
(cf. Corollary~\ref{csing}.\ref{exboundf}).
\end{corollary}

\begin{proof}
Take $n=1/A$, where $A=\min M(1)$.
Note that any terminal {\it singularity\/} of index $n$
has the m.l.d.$=1/n$.
\end{proof}

This result is much weaker than
Mori's on the good element in $|-K|$,
when $X$ is formally $\mathbb Q$-factorial$/Z$,
and has at most $1$ curve$/P$.
Indeed, as was remarked 
to the author by Prokhorov,
then there exists a good element $D\in |-K|$
according to Mori and Kollar.
(It is unknown whether it holds
when the number of curves$/P$ is $>1$ or
when $X/Z$ is not $\mathbb Q$-factorial.) 
So, $(X/Z,D)$ and $(X/Z,0)$ have
a regular complement by \cite[Theorem~5.12]{Sh2}.
In general, the corollary shows that exceptional cases
are the most difficult in combinatorics.
On the other hand we anticipate
a few exceptions (maybe, none) 
among them in the terminal case.

\refstepcounter{subsec}
\begin{example} Let $(X/Z,B)$ be a divisorial
contraction with a surface $B=E/P$, and
let $(K+B)_E=K_E+B_E$ be 
of type (A$_2^6$) or (I$_2^2$)
in~\ref{ecomp}.\ref{emainth}.3.
Then $(X/Z,E)$ has a trivial $7$-complement.
In this case, the singularity $P$ in $Z$ is
(maximal) log canonical of the index $7$,
but it is not log terminal.
\end{example}

\refstepcounter{subsec}
\begin{example}
Let $(X/Z,B)$ be a divisorial
contraction with a surface $B=E=\mathbb P^2/P$, and
let $(K+B)_E=K_E+B_E$ be 
of type (A$_2^1$) in~\ref{ecomp}.\ref{emainth}.3.
Moreover suppose that
$B_E=b_1 L_1+b_2 L_2 +(2/3)L_3+(1/2)L_4$, 
where $L_i$ are straight lines in a general position
and $b_1,b_2=6/7$.
Note that in such a situation 
coefficients $b_i$ are always standard 
\cite[Proposition 3.9]{Sh2} and 
we have a finite choice of them.

Then $(X/Z,E)$ is not regular,
it has a $42$-complement.
Moreover, in this case, $P$ is a purely log terminal
singularity, but {\it not terminal\/} or {\it canonical\/},
except for the case when $X$ has only Du Val
singularities along curves $L_i$ and $K\equiv 0/Z$.
In particular, $P$ has a crepant desingularization.
Here we do not discuss an existence of such
singularities.

Indeed, $E$ has index $42$.
So, for a straight line $L$ in $E$,
$(E.L)<0$ and $42(E.L)$ is an integer.
By the adjunction,
$(K+E.L)=(K_E+B_E.L)=-5/42$ and 
$(E.L)\ge -5/42$, when $P$ is canonical.
Moreover, $(E.L)> -5/42$ in such cases,
except for the above exception with $K\equiv 0/Z$. 
For $(E.L)=-5/42$, 
$X/Z$ is a crepant blow-up with
log canonical singularities along $L_i$.

Hence in the other cases,
$(E.C)=-m/42$
with some integer $m=1,2,3$ or $4$.
Then $(K.L)=-(5-m)/42$ and
discrepancy $d=d(E,0,Z)=(5-m)/m$.
If $P$ is terminal of the index $N\ge 2$,
so it should have at least all $N-1$ discrepancies
$i/N, 1\le i< N$ \cite[Theorem~3.2]{Sh3}.
For instance, if $m=1$ and $X$ has only Du Val
singularities, then $d=4$ and $N\ge 21$.
On the other hand,
making blow-ups over lines $L_i$, 
it is possible to construct
a (minimal log) resolution $Y/X$ with at most
$1+1+2+2\times 6=16$ exceptional divisors.
Then all  other divisors have discrepancies $>1$.
(Moreover the divisors in the resolution give
only a discrepancy $1/7<1$.)
So, it is not a terminal singularity.
Similarly we can exclude other cases.
Therefore $P$ is canonical or worse.

Of course, this approach uses a classification
of terminal singularities.
But it can be replaced by the following arguments.

Let $m=4$, then $d=1/4$ and, 
for any exceptional divisor $E'/P$,
the discrepancy $d(E',0,Z)=d(E',K-(1/4)E,X)=
d(E',K,X)+(1/4)\mult{E'}{E}$.
For instance, $L_4$ is a simple Du Val singularity
and, for $E'/L_4$ on its minimal resolution,
$d(E',0,Z)=d(E',K,X)+(1/4)\mult{E'}{E}=
0+(1/4)(1/2)=1/8$.
But, if $L_3$ is a simple singularity, i.e.,
a divisor $E'/L_3$ is unique on a minimal resolution,
then $d(E',0,Z)=d(E',K,X)+(1/4)\mult{E'}{E}=
-1/3+(1/4)(1/3)=-1/4$.
In this case $P$ is not canonical.
Otherwise, $L_3$ is a Du Val singularity 
of type $\mathbb A_2$.
The same works for other singularities $L_i$,
whenever $d< 1$,
that they are Du Val too.
Otherwise they have a discrepancy $<0$.
It can be checked by induction on 
the number of divisors on a minimal resolution.
But then on a minimal resolution $g:Y\to X$,
\begin{eqnarray*}
-4/42=(E.L)&=(g^*E.g^{-1}L)=
(g^{-1}E.g^{-1}L)+(\mult{L_i}{E}.L)\\
&=I+(1/2)+(2/3)+2\times(6/7)=
I+3-5/42,
\end{eqnarray*}
where $I=(g^{-1}E.g^{-1}L)$ is an integer.
This is impossible.
Similarly, we can make other cases.
\end{example}

\refstepcounter{subsec}
\begin{p-definition}
Let $(X,B)$ be a log terminal pair.
Put $D=\sum D_i=\lfloor B\rfloor$.
Then we can define a
simplicial space $\R{(X,B)}$:
\begin{itemize}
\item
its $l$-simplex is an irreducible component 
$\Delta_l$
in intersection of $l+1$ distinct irreducible
components $D_{i_0},\dots D_{i_l}$:
$$
\Delta_l\subseteq D_{i_0}\cap\dots\cap D_{i_l};
$$
\item
$\Delta_{l'}$ is a face of $\Delta_l$
if $\Delta_{l'}\supseteq \Delta_l$; and
\item 
the intersection of
two simplices $\Delta_l$ and $\Delta_{l'}$ 
consists of a finite
set of simplices $\Delta_{l''}$ such that
$$
\Delta_{l''}\supseteq \Delta_l\cup\Delta_{l'}.
$$
\end{itemize}
The simplices $\Delta_l$ give
a triangulation of $\R{(X,B)}$ 
or a simplicial complex if and only if
we have real global normal crossings in the generic points:
all the intersections $D_{i_0}\cap\dots\cap D_{i_l}$
are irreducible.
The latter can be obtained for
an appropriate log terminal resolution
$(S'/S,B')$.

If $(X,B)$ has a log terminal resolution
$(Y/X,B_Y)$, then the {\it topological type of 
$\R{(Y,B_Y)}$ is independent on such a resolution\/}.
So, we denote it by $\R{(X,B)}$.
The topology of $\R{(X,B)}$ reflects
a complexity of log singularities for $(X,B)$
and in particular of $\LCS(X,B)=\lfloor B\rfloor$
when $(X,B)$ is log terminal.

Put $\reg{(X,B)}=\dim \R{(X,B)}$.

When $X/Z$ is a contraction,
we {\it assume\/} that
components $\Delta_l$ are irreducible
formally or in the analytic topology
on $Z$, i.~e., we consider irreducible
branches over a neighborhood of $P\in Z$.
\end{p-definition}

\begin{proof}
According to Hironaka,
it is enough to verify that 
a monoidal transform in $\Delta_l\subset X$
gives a barycentric triangulation of
$\Delta_l$ in $\R{(X,B)}$.  
\end{proof}

\refstepcounter{subsec}
\label{sabclas}
\begin{example}
If $(S,B)$ is a surface singularity, then
$\R{(S,B)}$ is a graph of $\LCS{(S',B')}$
for a log terminal resolution $(S',B')\to (S,B)$.

Moreover, $\R{(S,B)}$ is homeomorphic to
a circle $S^1$, to a segment $[0,1]$, 
to a point, or to an empty set, when $(S,B)$ is log canonical.
Additionally, the case $S^1$ is
only possible when $B=0$ near the singularity
and it is elliptic with a wheel of rational curves 
for a minimal resolution.

Now let $(S/Z,B)$ an $r$-complement.
Then 
\begin{itemize}
\item
$\reg{(S,B)}=1$ if $(S/Z,B)$ has type
$\mathbb A_m^n$ or $\mathbb D_m^n$,
\item
$\reg{(S,B)}=0$ when $(S/Z,B)$ has type
$\mathbb Er_m^n$ with $(2\ge )m+n\ge 1$, and 
\item
$\reg{(S,B)}=-\infty$ when $(S/Z,B)$ has type
$\mathbb Er_m^n$ with $m+n= 0$, i.~e.,
when $K+B$ is Kawamata log terminal.  
\end{itemize} 

If $(S,B)=(\mathbb P^2, L)$ where
$L=\sum L_i$ with $n$ lines $L_i$ in a generic position,
then $\R{(\mathbb P^2,L)}$ is a complete graph
with $n$ points.
So, it is a manifold with boundary only when
$n\le 3$ or $-(K+L)$ is nef.
\end{example}

We have something similar in dimension $3$.

\refstepcounter{subsec}
\label{topsm}
\begin{p-definition}
Let $(X/Z,B)$ be an $n$-fold contraction.
Then the space $\R{(X,B)}$
has the following property:
\begin{description}
\item {\rm (DIM)}
$\R{(X,B)}$ is a compact topological space
of real dimension
$\reg{(X,B)}\le n-1=\dim X-1$.
\end{description}
Suppose now that
\begin{itemize}
\item
$-(K+B)$ is nef$/Z$, and
\item
$(X,B)$ is log canonical.
\end{itemize}
Then locally$/Z$
\begin{description}
\item {\rm (CN)}
$\R{(X,B)}$ is connected,
whenever $-(K+B)$ is big$/Z$,
or, for $n\le 3$, consists of two points;
moreover, the latter is possible only
if $\dim Z\le \dim X-1= 2$ and there exists
a (birationally unique) conic bundle structure on $X/Z$ with
$2$ reduced and disjoint sections $D_1$ and $D_2$ 
in $B'$
for a log terminal resolution $(X',B')\to (X,B)$;
$\R{(X,B)}=\{D_1,D_2\}$;
\item {\rm (MB)}
$\R{(X,B)}$ is a manifold with boundary;
in addition, it is a manifold when
$(X,B)$ is a $1$-complement and
$B$ is over a given point $P\in Z$.
\end{description}

In particular, we associate to each
birational contraction $(X/Z,B)$
or a singularity, when $X\to Z$ is
an isomorphism, 
a connected manifold $\R{(X,B)}$
locally over a point $P\in Z$,
which will be called a {\it type\/} of $(X/Z,B)$.
The {\it regularity\/} $\reg{(X,B)}$
characterizes its topological difficulty.

If $(X/Z,B)$ is an $n$-complement, that
is not Kawamata log terminal over $P$
as in Theorem~\ref{csing}.\ref{singth},
then $\R{(X,B)}\not=\emptyset$ and
$\reg{(X,B)}\ge 0$.
If $(X/Z,B)$ is an arbitrary log canonical
singularity, we associate with it
a topological manifold of the maximal
dimension (and maximal for inclusions) 
for some complements
$\R{(X,B^+)}$ and its {\it complete\/} 
regularity is $\reg{(X,B^+)}$.  
\end{p-definition}

\begin{proof} (DIM) holds in any dimension
by the very definition.

For $n=3$-folds,
(CN) follows from the LMMP and
proofs of \cite[Connectedness Lemma~5.7
and Theorem~6.9]{Sh2}.
The connectedness, when $-(K+B)$ is big$/Z$,
was proved in \cite[Theorem~17.4]{KC}.

(MB) is a local question on $\R{(X,B)}$
near each point $\Delta_0=D_{i_0}$.
But this is a global question on
$Y=D_{i_0}$ that $\R{(Y,B_Y)}$
satisfies (MB).
More precisely, a neighborhood of $\Delta_0$
is a cone over $\R{(Y,B_Y)}$.
But we may assume (MB) for $\R{(Y,B_Y)}$
by the adjunction and an induction
on $\dim X$.
\end{proof}

\refstepcounter{subsec}
\begin{corollary}
Under the assumptions of Theorem~\ref{csing}.\ref{singth},
let $(X/Z,B^+)$ be an $n$-complement
with the minimal index $n$, then
\begin{itemize}
\item
$\reg{(X,B^+)}=2$ and $\R{(X,B^+)}$ is
a real compact surface with boundary only when $n=1$ or $2$;
$\R{(X,B^+)}$ is closed only for $n=1$;
\item
$\reg{(X,B^+)}=1$ and $\R{(X,B^+)}$ is
a real curve with boundary only when $n=1,2,3,4$ or $6$; and
\item
$\reg{(X,B^+)}=0$ and $\R{(X,B^+)}$ is
a point only when $n\in N_2$;
such complements and contractions are
exceptional.
\end{itemize}
\end{corollary}

\begin{proof}
Follows from the proof of 
Theorem~\ref{csing}.\ref{singth}.
\end{proof}

A topology of log singularities
can be quite difficult when $\reg{(X,B^+)}=2$.
By a real surface below we mean
a connected compact manifold with boundary
of dimension $2$.
It is {\it closed\/} when the boundary is empty. 

\refstepcounter{subsec}
\begin{example}
For any closed real surface $S$,
there exists a $3$-fold $1$-complement $(X,B)$
such that $\R{(X,B)}$ is homeomorphic to
$S$.

First, we take a triangulation $\{\Delta_i\}$ of $S$.

Second, we immerse its dual into $\mathbb P^3$
in such a way that each point $\Delta_0$
is presented by a plane $L_i$.
The plane is in a generic position.

Third, we make a blow-up in such
intersections $L_i\cap L_j$ which does not
correspond to a segment $\Delta_1=L_i\cap L_j$.
Then we get an algebraic surface $B=\sum L_i$
such that $K+B$ is log terminal and
$\R{(X,B)}$ is the triangulation.
So, $\R{(X,B)}$ is homeomorphic to $S$.
Moreover, $K+B$ is numerically trivial on 
the $1$-dimensional skeleton or on each 
curve $C_{i,j}=\Delta_1=L_i\cap L_j$,
because we have exactly two triple-points on each
$\Delta_1$ or each $\Delta_1$ belongs exactly to
two simplices of the triangulation.

Fourth, we contract all $C_{i,j}$ and
something else from $B$ and $X$ to
a point which gives the required singularity.
Indeed, after the LMMP for $K+B$, we can use
a semi-ampleness of $K+B$, when it has
a general type. 
Note that the birational contractions
or flips do not touch $C_{i,j}$.
This can be verified on each $L_i$
by the adjunction.
By the same reason we have no surface contraction on
$L_i$ or flips intersecting $B$, but not in $B$.
However, terminal singularities and flips are possible
outside of $B$ or inside of $B$, which preserves
our assumptions on $B$.
To secure the big property for $K+B$
we can add similarly $B'$ on which
$K+B$ is big.
Then $K+B$ has the log Kodaira dimension $3$
when $B'$ has more than two connected components,
otherwise $(K+B+B')|_{B'}$ will not be big
on some of the components.

The same holds for $2$-complements 
with arbitrary real non-closed surface $S$ with boundary.
We need not to contract birationally some $L_i$ 
and replace them by $(1/2)D$ for generic $D\in |2L_i|$.
This can be done by the above combinatorics,
when $(K+B+B')_{L_i}$ is big. 
\end{example}

\refstepcounter{subsec}
\label{ball}
\begin{corollary}
Let $P\in (X,B)$ be a log canonical singularity
such that 
\begin{itemize}
\item 
$X$ is $\mathbb Q$-factorial, and even formally
or locally in the analytic topology of $X$,
when there exists a non-normal curve in $P$ as 
a center of log canonical discrepancy $0$
for $K+B$; 
\item
$\{B\}\not=0$ in $P$, i.~e., $B$
has a fractional component through $P$; and
\item
$K+\lfloor B\rfloor$ is purely
log terminal near $P$.
\end{itemize}
Then $\R{(X,B)}$ has type $\mathbb B^r$,
where $\mathbb B^r$ is a ball of dimension 
$r=\reg{(X,B)}$.

Moreover, we can drop the conditions when
$r\le 0$, or we can drop the first condition when
\begin{itemize}
\item 
$B$ has a fractional $\mathbb R$-Cartier component $F$,
i.~e., $0<F\le\{B\}$.
\end{itemize}
\end{corollary}

\begin{proof}
According to our assumptions,
$K+\lfloor B\rfloor$ is purely log terminal 
in $P$.
Then $X$ is the only log minimal model
of $(X,\lfloor B\rfloor)$ over $X$ \cite[1.5.7]{Sh2}.

Now we take a log terminal resolution $(Y/X,B_Y)$
and consider formally the LMMP$/X$ for $\lfloor B_Y\rfloor$.
According to the above, the final model
will be $(X,\lfloor B\rfloor)$ 
with $\reg{(X,\lfloor B\rfloor)}\le 1$ and it has
a trivial homotopy type.
On the other hand $\R{(X,B)}=\R{(Y,B_Y)}$.
So, it is enough to check that contractions
and flips preserve the homotopy type.
If the centers of the flip or contraction
are not in $\LCS{(X,B)}=\lfloor B_Y\rfloor$, 
then we even have a homeomorphism.
If we have a divisorial contraction
of a divisor $D_i$ in $\lfloor B_Y\rfloor$, it
induces a fibre contraction $D_i\to Z/X$.
If $Z$ is a curve$/P$, then according to our
conditions this is a curve on another component
$D_j$ in $\lfloor B_Y\rfloor$,
because any contraction of $Y/X$ is divisorial.
Note that each exceptional divisor of $Y/X$
belongs to $\lfloor B_Y\rfloor$, and
$\R{(Y,B_Y)}$ is a gluing of a cone
with vertex $\Delta_0=D_i$ over
$\R{(D_i,B_{D_i})}$ and in the latter.
It is homotopy to $\mathbb B^{r-1}$,
because $-(K_Y+\lfloor B_Y\rfloor)$ is nef$/X$
on all double curves $\Delta_1=D_i\cap D_j$.
The surgery drops this cone.
We have a similar picture when $Z=\pt$
or $Z$ is not over $P$.
In the latter case $Z$ has at most one curves$/P$
formally.
(If non-formally, algebraically, then
at most two double curves$/P$, which 
in addition are connected, by our assumptions.)
Finally, let $Y-\to Y^+/X$ be a flip
in a curve $C_{i,j}=\Delta_1=D_i\cap D_j/P$.
Then according to our assumptions we have
a flip on a third surface $D_k$ in $B_Y$
with $(D_k.C_{i,j})=1$ on $Y$.
The surgery deletes the segment $\Delta_1=C_{i,j}$
and the interior of the triangle 
$\Delta_2=D_i\cap D_j \cap D_k$.
Formally, it looks like a blow-up in $\Delta_1$
which is the barycentric triangulation in $\Delta_1$
and then we {``}contract" the resolution divisor
on a curve on $D_k$.

If $r=1$,
we can also assume that at the start 
$\{B\}\equiv -(K_Y+\lfloor B_Y\rfloor)/Z$
is nef $Z$ \cite[Theorem~5.2]{Sh3}.
Indeed, we may consider the LMMP for
$B_Y+\varepsilon \{B\}$.
Since $-(K_Y+\lfloor B_Y\rfloor)$ is big$/Z$,
we have a birational contraction $Y\to Z/X$;
it contracts all surfaces to a curve $C$
by Theorem~\ref{clasf}.\ref{dipole}.
Then any fractional component of $B$ is positive on $C$,
which contradicts our assumption,
because $Z/X$ is small.

If $B$ has the fractional component $F$,
we can apply an induction on the number
of irreducible curves $C_i/P$ on
a formal $\mathbb Q$-factorialization $X'/X$.
Indeed, if $X'/X$ blows up $1$ such curve $C_i$,
it is irreducible rational, and it has
at most to points in which $K+B$ is 
not log terminal.
So we glue at most two balls $\mathbb B^r$
in $\mathbb B^{r-1}$.
Note that $F$ passes through $C_i$ on $X'$.
Finally, we have no such curves when
$X$ is formally $\mathbb Q$-factorial in $P$.
This drops the first condition.
\end{proof}

In general, for $r=2$ we can verify that 
the formal Weil-Picard number in $P$
is at least $q=h^1(\R{(X,B)},\mathbb R)
=2-\chi(\R{(X,B)})$
(the topological genus),
where $\chi$ is the topological Euler
characteristic.
We may consider $q$ as the {\it genus\/}
of the singularity $P$.

\refstepcounter{subsec}
\label{balls}
\begin{corollary}
Under the assumptions of Theorem~\ref{csing}.\ref{singth},
let $(X/Z,B^+)$ be an $n$-complement
with the minimal index $n$, and 
$P$ be $\mathbb Q$-factorial, formally,
when there exists a non-normal curve in $P$
as a center of log canonical discrepancy $0$
for $K+\lfloor B^+\rfloor$,
and purely log terminal $K+B$, then
$\R{(X,B)}$ has type $\mathbb B^r$
with $r=\reg{(X,B)}$, or $S^2$.
The latter is only possible for $n=1$.
\end{corollary}

If $P$ is an isolated singularity, we may drop
the formal condition for an appropriate
complement with $n\ge 2$.
For $n=1$, we can take a $2$-complement as
in the proof below.
For $n\ge 3$, we can drop even the $\mathbb Q$-factorial
property.

\begin{proof}
We need to consider only the case with $n=1$.
Then we have a reduced component $D$
through $P$ in $B^+$.
If we replace $D$ by a generic $(1/2)D'$,
where $D'\in |2D|$ is rather generic, we obtain
a $2$-complement $B'$, with
$\R{(X,B')}$ homeomorphic to $\mathbb B^2$.
Then $\R{(X,B^+)}$ could be obtained
from this by gluing a cone
over $\R{(D,B_D)}$.
The letter is $[0,1]$ or $S^1$.
That gives respectively $\mathbb B^2$ or $S^2$.
\end{proof}

In the investigation of the m.l.d.'s 
for $3$-folds we can assume that the point $P$ in
$(X,B)$ is 
\begin{description}
\item{(T2)}
$\mathbb Q$-factorial (even formally)
and terminal in codimension $1$ and $2$.
\end{description}
Otherwise after a crepant resolution we can
reduce the problem to (T2) or 
to dimension $2$ or $1$.

\refstepcounter{subsec}
\begin{corollary}
The conjecture on discrepancies 
{\rm \cite[Cojecture 4.2]{Sh3}} holds
for $3$-fold log singularities $(X,B)$
with $b_i \in \Gamma$ under {\rm (M)''} and
{\rm (T2)}, 
when $\reg{(X,B^+)}\le 1$.
Moreover, in such a case 
the only clusters of $A(\Gamma,3)$ are
\begin{description}
\item{\rm (0)} 
$0$, when $\reg{(X,B^+)}=0$, {\rm i.~e.,
in the exceptional cases\/};
\item{\rm (1)}
$0$ and
$$A(\Gamma',2)$$
when $\reg{(X,B^+)}=1$,
where $\Gamma'=\{0,1/2,2/3,3/4,4/5, 5/6\}$, 
and where we consider
only $2$-dimensional singularities $(S,B^+)$
with the boundary multiplicities in $\Gamma'$
and $\reg{(S,B^+)}=0$, {\rm i.~e.,
the exceptional case.}
\end{description}

The cluster points are rational,
our $A(\Gamma,3)$ is closed when $1\in \Gamma$,
and the only cluster of the clusters
is $0$.
\end{corollary}

\begin{proof}[Sketch Proof]
Note that $\Gamma$ satisfies (M)'', but
may not be the d.c.c..
In particular, $\Gamma$ may not be standard.
If $\reg{(X,B)}=0$, then by 
Corollary~\ref{csing}.\ref{exbound},
for any $\varepsilon>0$,
we have a finite subset
$$\{a\in A(\Gamma,3)\mid a\ge \varepsilon\}$$ 
of corresponding $A(\Gamma,3)$.
Of course, the corollary was proved
for the standard $\Gamma$, but 
in our case, according to 
Corollary~\ref{indcomp}.\ref{adjunct}
and because $\reg{(X,B^+)}=0$ all $b_i\le (l-1)/l$,
where $l=\max\{r\in N_2\}$.
So, such $b_i$ is standard due to (M)''.

If $\reg{(X,B)}=1$, we can use the arguments
of \cite{Sh5} and Theorem~\ref{clasf}.\ref{dipole}.
The latter almost reduces our case to
the $2$-dimensional situation.
Almost means except for the edge components
in the chain of the reduced part $S$ of $B^+$
as in the proof of Theorem~\ref{csing}.\ref{singth}
(cf. the proof of Corollary~\ref{csing}.\ref{ball}).
The clusters can be realized as m.l.d.'s for
$(X,B)$ with a reduced part $S=\lfloor B\rfloor$
and the the fractional multiplicities
in $\Gamma'$.
\end{proof}

So, to complete the conjecture on discrepancies for $3$-folds
someone needs to consider the singularities
with $2$-complements of type $\mathbb B^2$ by
Corollary~\ref{csing}.\ref{balls}.
This case is related to $1$- or $2$-complements.
The former are closed to toric complements,
where the conjecture was verified by Borisov \cite{B}.
We also see that $\reg{(X,B^+)}$ conjecturally may have
interpretations in terms of clusters:
the first cluster of $A(\Gamma,n)$ with
$r=\reg{(X,B^+)}$ are $A(\Gamma,n-1)$ with
$\reg{(X,B^+)}=r-1$, in particular, they
are rational when $\Gamma$ is standard.

\refstepcounter{subsec}
\begin{remark}
We anticipate that most of the results in this section 
hold in any dimension and for any regularity
$r=\reg{(X,B)}$ or $\reg{(X,B^+)}$.
\end{remark}

\bigskip
The Johns Hopkins University,
Department of Mathematics,
Baltimore, Maryland, 21218.

shokurov@math.jhu.edu

\end{document}